\documentclass[reqno,10pt,a4paper,twoside]{amsart}
\usepackage{amsfonts,amsmath,amssymb}
\usepackage[dvips]{geometry, graphicx}
\usepackage{color}

\usepackage{hyperref}

\numberwithin{equation}{section}

\newcommand{\beq}{\begin{equation}}
\newcommand{\eeq}{\end{equation}}
\newcommand{\bs}{\begin{split}}
\newcommand{\esplit}{\end{split}}
\newcommand{\begincal}{\begin{eqnarray*}}
\newcommand{\fincal}{\end{eqnarray*}}

\newcommand{\rtwo}{{\mathbb R}^2}

\newcommand{\Om}{\Omega}
\newcommand{\eps}{\varepsilon}

\newcommand{\ds}{\displaystyle}

\newcommand{\xe}{x_\eps}
\newcommand{\xie}{x_{i,\eps}}
\newcommand{\mie}{\mu_{i,\eps}}
\newcommand{\xje}{x_{j,\eps}}
\newcommand{\mje}{\mu_{j,\eps}}
\newcommand{\gie}{\gamma_{i,\eps}}
\newcommand{\gje}{\gamma_{j,\eps}}
\newcommand{\rie}{r_{i,\eps}}
\newcommand{\rje}{r_{j,\eps}}
\newcommand{\xke}{x_{k,\eps}}

\newcommand{\gke}{\gamma_{k,\eps}}
\newcommand{\rke}{r_{k,\eps}}
\newcommand{\sie}{s_{i,\eps}}
\newcommand{\sje}{s_{j,\eps}}
\newcommand{\yie}{y_{i,\eps}}
\newcommand{\nie}{\nu_{i,\eps}}

\newcommand{\die}{d_{i,\eps}}

\newtheorem{thm}{Theorem}[section]
\newtheorem{lemma}{Lemma}[section]

\newtheorem{prop}{Proposition}[section]

\newtheorem{claim}{Claim}[section]

\begin{document}

\title[Multi-bumps analysis for Trudinger-Moser nonlinearities in 2d]
{Multi-bumps analysis for Trudinger-Moser nonlinearities\\ I - Quantification and location of concentration points}

\author{Olivier DRUET}
\address{Olivier Druet et Pierre-Damien Thizy\\
Univ Lyon, Universit\'e Claude Bernard Lyon 1, CNRS UMR 5208, Institut Camille Jordan, 43 blvd. du 11 novembre 1918, F-69622 Villeurbanne cedex, France}
\author{Pierre-Damien THIZY}

\date{October, 23rd, 2017}

\begin{abstract}
In this paper, we investigate carefully the blow-up behaviour of sequences of solutions of some elliptic PDE in dimension two containing a nonlinearity with Trudinger-Moser growth. A quantification result had been obtained by the first author in \cite{DruetDMJ} but many questions were left open. Similar questions were also explicitly asked in subsequent papers, see Del Pino-Musso-Ruf \cite{DPMR1}, Malchiodi-Martinazzi \cite{MalchiodiMartinazzi} or Martinazzi \cite{Martinazzi}. 
We answer all of them, proving in particular that blow up phenomenon is very restrictive because of the strong interaction between bubbles in this equation. This work will have a sequel, giving existence results of critical points of the associated functional at all energy levels via degree theory arguments, in the spirit of what had been done for the Liouville equation in the beautiful work of Chen-Lin \cite{ChenLin-Liouville}.
 
\end{abstract}

\maketitle


\section{Introduction}

We let $\Omega$ be a smooth bounded domain of $\rtwo$ and we consider the equation 
\begin{equation}\label{maineq}
\Delta u = \lambda f u e^{u^2} \hbox{ in }\Omega,\, u>0 \hbox{ in }\Omega,\, u=0\hbox{ on }\partial \Omega\hskip.1cm.
\end{equation}
where $\Delta = -\frac{\partial^2}{\partial x^2}-\frac{\partial^2}{\partial y^2}$, $\lambda>0$ and $f$ is a smooth positive function in $\overline{\Omega}$. 

\smallskip This equation is critical with respect to Trudinger-Moser inequality. Indeed, the nonlinearity in $e^{u^2}$ is the best one can hope to control in dimension $2$ by the $L^2$-norm of the gradient. More precisely, we let $H^1_0\left(\Omega\right)$ be the standard Sobolev space (with zero boundary condition) endowed with the norm ${\ds \left\Vert \nabla u\right\Vert_2^2=\int_\Omega \left\vert \nabla u\right\vert^2\, dx}$. Trudinger proved in \cite{Trudinger} that ${\ds \int_\Omega e^{u^2}\, dx}$ is finite for any function $u$ in $H_0^1\left(\Omega\right)$. Moser was then a little bit more precise in \cite{Moser}, proving that 
\begin{equation}\label{eq-inequality-Moser}
\sup_{u\in H^1_0\left(\Omega\right),\, \left\Vert \nabla u\right\Vert_2=1} \int_\Omega e^{\gamma u^2}\, dx < +\infty \hbox{ if and only if }\gamma\le 4\pi\hskip.1cm.
\end{equation}
Solutions of equation \eqref{maineq} are in fact critical points of the functional 
\begin{equation}\label{eq-functional} 
J(u)= \int_\Omega fe^{u^2}\, dx
\end{equation}
under the constraint ${\ds \int_\Omega \left\vert \nabla u\right\vert^2\, dx =\beta}$ for some $\beta>0$. The $\lambda$ appearing in \eqref{maineq} is then the Euler-Lagrange coefficient. This functional is well-defined on $H^1_0\left(\Omega\right)$ thanks to Trudinger \cite{Trudinger}. It is also easy to find a critical point of $J$ if $\beta<4\pi$ in the constraint thanks to Moser's inequality \eqref{eq-inequality-Moser}~: these critical points may be found as maxima of $J$ under the constraint ${\ds \int_\Omega \left\vert \nabla u\right\vert^2\, dx =\beta<4\pi}$. However, as studied by Adimurth-Prashanth \cite{AdiPrashanth1}, for $\beta=4\pi$, finding critical points is more tricky since a lack of compactness appears in Palais-Smale sequences at this level of energy. Nevertheless, it has been proved by Carleson-Chang \cite{Carleson-Chang} for the unit disk, by Struwe \cite{Struwe-MT} for $\Omega$ close to the disk and by Flucher \cite{Flucher} for a general $\Omega$ that there are extremals in Moser's inequality \eqref{eq-inequality-Moser} for $\gamma=4\pi$, meaning in particular that there are always critical points of $J$ for the critical value $\beta=4\pi$. Note that existence of critical points for $\beta$ slightly larger than $4\pi$ has also been proved by Struwe \cite{Struwe-MT} and Lamm-Robert-Struwe \cite{LRS}. Struwe \cite{Struwe-MT2} also found critical points of higher energy (for some values of $\beta$ between $4\pi$ and $8\pi$) when the domain contains an annulus (in the spirit of Coron \cite{Coron}).  We refer also to the recent Mancini-Martinazzi \cite{ManciniMartinazzi} for an interesting new proof of the existence of extremal functions for Moser's inequality in the disk without using test-functions computations.  

\smallskip In the last decade, tools have been developed to study sequences of solutions of equation \eqref{maineq} and in particular to understand precisely their potential blow-up behaviour. This serie of works started in the minimal energy situation ($\beta$ close to $4\pi$) with Adimurthi-Struwe \cite{AdiStruwe}. Then Adimurthi-Druet \cite{AdiDruet} used this blow-up analysis to obtain an improvement of Moser's inequality (completing the result of Lions \cite{Lions2}). In the radial case (that is in the unit disk with $f\equiv 1$), such a blow-up analysis in the minimal energy case was recently used by Malchiodi-Martinazzi \cite{MalchiodiMartinazzi} to prove that there is a $\beta_0>4\pi$ for which there are solutions of \eqref{maineq} of energy less than or equal to $\beta_0$ but no solutions of energy greater than $\beta_0$. 

\smallskip In order to get solutions of higher energies and to describe precisely the set of solutions for all $\beta$, one needs a fine analysis of blowing-up solutions. The first result in this direction is the quantification result of the first author \cite{DruetDMJ} that we recall here since the questions we adress in the present work come from it~:

\begin{thm}[Druet \cite{DruetDMJ}]\label{thm-DMJ-intro}
Let $\Om$ be a smooth bounded domain in $\rtwo$ and let 
$\left(f_\eps\right)_{\eps>0}$ be a sequence of functions of uniform critical growth in $\Om$. Also let 
$\left(u_\eps\right)_{\eps>0}$ be a sequence of solutions of 
$$\Delta u_\eps = f_\eps\bigl(x,u_\eps(x)\bigr)$$
verifying that $\left\Vert \nabla u_\eps\right\Vert_2^2\to \beta$ as $\eps\to 0$ for some $\beta\in {\mathbb R}$. Then there exists 
a solution $u_0\in C^0\left(\bar{\Om}\right)$ of 
$$\Delta u_0 = f_0\bigl(x,u_0(x)\bigr)\hbox{ in }\Om\hbox{, }u_0= 0\hbox{ on }\partial \Om
\hskip.1cm,$$
and there exists $\tilde{N}\in {\mathbb N}$ such that 
$$\left\Vert \nabla u_\eps\right\Vert_2^2 = \left\Vert \nabla u_0\right\Vert_2^2 + 4\pi \tilde{N}+o\left(1\right)\hskip.1cm.$$
If $\tilde{N}=0$, the convergence of $u_\eps$ to $u_0$ is strong in $H_0^1\left(\Omega\right)$ and actually holds in $C^0\left(\bar{\Om}\right)$. 
\end{thm}

\medskip We do not define here sequences of functions of uniform critical growth in $\Omega$. The only thing we need to know is that they include sequences of the form $f_\eps(x,u) = h_\eps(x)u e^{u^2}$ as soon as $h_\eps>0$ and $h_\eps\to h_0$ in ${\mathcal C}^1\left(\overline{\Omega}\right)$. But they include much more general nonlinearities behaving like $e^{u^2}$ at infinity. Note also that, in the litterature, the nonlinearity is sometimes written as $e^{4\pi u^2}$ (this is for instance the case in \cite{DruetDMJ}), hence the discrepancy of $4\pi$ in some results. 

\medskip This result describes precisely the lack of compactness in the energy space. Note that this result is not true for Palais-Smale sequences, as proved by Costa-Tintarev \cite{CostaTintarev}~: there are Palais-Smale sequences for the above equation which converge to $0$ weakly in $H^1_0\left(\Omega\right)$ and which present a lack of compactness at any level above $4\pi$. This shows that the quantification result of Theorem \ref{thm-DMJ-intro} is specific to sequences of solutions of the equation and require a pointwise analysis as carried out in \cite{DruetDMJ} and below (not only an analysis in energy space).

\medskip Note also that the above result is not empty since Del Pino-Musso-Ruf \cite{DPMR2} constructed, via a Liapunov-Schmidt procedure, multi-pikes sequences of solutions of equation \eqref{maineq} (with $f\equiv 1$) in annuli. These sequences satisfy the hypothesis of the above theorem, converge weakly to $0$ in $H^1_0\left(\Omega\right)$ (that is $u_0\equiv 0$ in the above result) and have an energy converging to $4\pi \tilde{N}$. They can construct such solutions for all $\tilde{N}\ge 1$. This suggests that the topology of the domain plays a crucial role in the existence of solutions of arbitrary energies.

\medskip However, if one wants to push further the existence results, we need to be more precise than Druet \cite{DruetDMJ}. In particular, we need to answer the following natural questions, left open in this work of the first author (see also Del Pino-Musso-Ruf \cite{DPMR1}, Malchiodi-Martinazzi \cite{MalchiodiMartinazzi} or Martinazzi \cite{Martinazzi} where one can find these, or similar, questions)~:

\smallskip 1. Is it possible to have both $u_0\not\equiv 0$ and $\tilde{N}\ge 1$ in the above theorem ? 

\smallskip 2. Are the concentration points appearing when $\tilde{N}\ge 1$ isolated\footnote{By isolated, we mean here that the energy at any concentration point is exactly $4\pi$. In other words, we mean that there are no bubble accumulations and we do not wish to rule out only bubbles towers.} or not ? If yes, where are they ?

\medskip These questions are natural and can be motivated by analogy with Liouville type equations (see among others \cite{Chen-Lin,ChenLin-Liouville,LiShafrir,Linsurvey,Malchiodi,Takahashi}) or Yamabe type equations (see for instance \cite{Brendle,BrendleMarques,DruetJDG,DHRbook,KhuriMarquesSchoen,LiZhang,LiZhang2,LiZhang3, LiZhu,Marques,Marques2,SchoenZhang}). We refer in particular to \cite{DruetENSAIOS} for a survey on this kind of questions. 

\medskip We attack in this paper the questions 1 and 2 above. Our result holds for more general nonlinearities but we restrict, for sake of clearness, to the simplest one. We consider a sequence $\left(u_\eps\right)$ of smooth positive solutions of 
\begin{equation}\label{eq}
\Delta u_\eps = \lambda_\eps f_\eps u_\eps e^{u_\eps^2} \hbox{ in }\Omega,\, u_\eps=0\hbox{ on }\partial \Omega\hskip.1cm,
\end{equation}
for some sequence $\left(\lambda_\eps\right)$ of positive real numbers and some sequence $\left(f_\eps\right)$ of smooth functions in $\overline{\Omega}$ which satisfies 
\begin{equation}\label{eq-cv-heps}
f_\eps\to f_0\hbox{ in }{\mathcal C}^1\left(\overline{\Omega}\right)\hbox{ as }\eps\to 0\hbox{ and }\left\Vert \nabla^2 f_\eps\right\Vert_{L^\infty\left(\Om\right)} = O(1)
\end{equation}
where $f_0>0$ in $\overline{\Omega}$. And we prove the following result~:

\begin{thm}\label{mainthm}
Let $\Omega$ be a smooth bounded domain of $\rtwo$ and let $\left(u_\eps\right)$ be a sequence of smooth solutions of \eqref{eq} which is  bounded in $H^1_0\left(\Omega\right)$. Assume that \eqref{eq-cv-heps} holds. Then, after passing to a subsequence, $\lambda_\eps\to \lambda_0$ as $\eps\to 0$ for some $\lambda_0\in {\mathbb R}$. . 

\smallskip If $\lambda_0\neq 0$, then there exists $u_0\in {\mathcal C}^2\left(\overline{\Omega}\right)$ solution of 
$$\Delta u_0=\lambda_0 f_0 u_0 e^{u_0^2}\hbox{ in }\Omega,\, u_0=0\hbox{ on }\partial \Omega$$
such that $u_\eps\to u_0$ in ${\mathcal C}^2\left(\overline{\Omega}\right)$ as $\eps\to 0$.

\smallskip If $\lambda_0=0$, then $u_\eps\rightharpoonup 0$ weakly in $H^1_0\left(\Omega\right)$. Moreover there exist $N\ge 1$ such that 
$$\int_\Omega \left\vert \nabla u_\eps\right\vert^2\, dx \to 4\pi N\hbox{ as }\eps\to 0$$
and $N$ sequences of points $\left(x_{i,\eps}\right)$ in $\Omega$ such that 

\smallskip a) $x_{i,\eps}\to x_i$ as $\eps\to 0$ with $x_i\in \Omega$ (not on the boundary), all the $x_i$'s being distinct.

\smallskip b) $u_\eps\to 0$ in ${\mathcal C}^2_{loc}\left(\overline{\Omega}\setminus {\mathcal S}\right)$ where ${\mathcal S}=\left\{x_i\right\}_{i=1,\dots,N}$.

\smallskip c) for all $i=1,\dots,N$, we have that $\gie=u_\eps\left(\xie\right)\to +\infty$ as $\eps\to 0$ and that 
$$\gie \left(u_\eps\left(\xie +\mie x\right)-\gie\right)\to U(x)=-\ln \left(1+\frac{1}{4} \vert x\vert^2\right)$$
in $C^2_{loc}\left(\rtwo\right)$ as $\eps\to 0$ where 
$$\mie^{-2} = \lambda_\eps f_\eps\left(\xie\right)\gamma_{i,\eps}^2 e^{\gamma_{i,\eps}^2} \to +\infty \hbox{ as }\eps\to 0\hskip.1cm.$$

\smallskip d) for all $i=1,\dots,N$, there exists $m_i>0$ such that 
$$\sqrt{\lambda_\eps}\gie \to \frac{2}{m_i\sqrt{f_0\left(x_i\right)}}\hbox{ as }\eps\to 0\hskip.1cm.$$

\smallskip e) The points $x_i$ are such that 
$$2m_i \nabla_y {\mathcal H}\left(x_i,x_i\right) +4\pi \sum_{j\neq i} m_j \nabla_y {\mathcal G}\left(x_j,x_i\right) +\frac{1}{2} m_i \frac{\nabla f_0\left(x_i\right)}{f_0\left(x_i\right)}=0$$
and that 
$$4\pi\sum_{j\neq i} m_j{\mathcal G}\left(x_j,x_i\right)+2m_i{\mathcal H}\left(x_i,x_i\right)+m_i\ln \frac{f_0\left(x_i\right)}{m_i^2}+m_i =0$$
for all $i=1,\dots,N$ where 
$${\mathcal G}(x,y) = \frac{1}{2\pi}\left(\ln\frac{1}{\left\vert x-y\right\vert} +{\mathcal H}\left(x,y\right)\right)$$ 
is the Green function of the Laplacian with Dirichlet boundary condition. 
\end{thm}

Note that this theorem proves that, if blow-up occurs, then the weak limit has to be zero so that lack of compactness can occur only at the levels $\beta=4\pi N$ for $N\ge 1$. This is a key information to get general existence result via degree theory from this theorem; this will be the subject of a subsequent paper. We also obtain a precise characterisation of the location of concentration points. This answers in particular by the affirmative to the conjecture of Del Pino-Musso-Ruf \cite{DPMR1} (p. 425) since, in case $f\equiv 1$, the $\left(x_i,m_i\right)$ of Theorem \ref{mainthm} are critical points of the function 
$$\Phi\left(y_i,\alpha_i\right) = 2\pi \sum_{i\neq j} \alpha_i \alpha_j{\mathcal G}\left(y_i,y_j\right)+\sum_{i=1}^N \alpha_i^2 {\mathcal H}\left(y_i,y_i\right) + \sum_{i=1}^N \left(\alpha_i^2 - \alpha_i^2 \ln \alpha_i\right)\hskip.1cm.$$

\medskip The paper is organized as follows. In Section \ref{section-DMJ}, we recall the main results of Druet \cite{DruetDMJ} and set up the proof of the theorem. Section \ref{section-LBA} is devoted to a fine asymptotic analysis in the neighbourhood of a given concentration point while the theorem is proved in Section \ref{section-proofmainthm} which deals with the multi-spikes analysis. At last, we collect some useful estimates concerning the standard bubble and the Green function respectively in appendices A and B.

\section{Previous results and sketch of the proof}\label{section-DMJ}

We set up the proof of Theorem \ref{mainthm} and we recall some results obtained in Druet \cite{DruetDMJ}. We let $\Om$ be a smooth bounded domain of $\rtwo$ and we consider a sequence $\left(u_\eps\right)$ of smooth positive solutions of 
\begin{equation}\label{eqeps}
\Delta u_\eps = \lambda_\eps f_\eps u_\eps e^{u_\eps^2} \hbox{ in }\Omega,\, u_\eps=0\hbox{ on }\partial \Omega
\end{equation}
for some sequence $\left(\lambda_\eps\right)$ of positive real numbers and some sequence $\left(f_\eps\right)$ of smooth functions which satisfies \eqref{eq-cv-heps}. Note that we necessarily have that 
\begin{equation}\label{eq-lambdaepsbounded}
\limsup_{\eps\to 0} \lambda_\eps \le \frac{\lambda_1}{\min_{\overline{\Omega}} f_0}
\end{equation}
where $\lambda_1>0$ is the first eigenvalue of the Laplacian with Dirichlet boundary condition in $\Omega$. Indeed, let $\varphi_1\in {\mathcal C}^\infty\left(\overline{\Omega}\right)$ be a positive (in $\Omega$) eigenfunction associated to $\lambda_1$ and multiply equation \eqref{eqeps} by $\varphi_1$. After integration by parts, we get that 
$$\lambda_1 \int_\Omega u_\eps\varphi_1\, dx = \lambda_\eps \int_\Omega f_\eps u_\eps e^{u_\eps^2}\varphi_1\, dx\hskip.1cm.$$
Since $f_\eps$ becomes positive for $\eps$ small thanks to \eqref{eq-cv-heps} and since $u_\eps$ and $\varphi_1$ are positive, we can write that 
$$\lambda_1 \int_\Omega u_\eps\varphi_1\, dx \ge \lambda_\eps \left(\min_{\overline{\Omega}} f_\eps\right) \int_\Omega u_\eps\varphi_1\, dx\hskip.1cm,$$
which leads to \eqref{eq-lambdaepsbounded}. 

\medskip We assume in the following that there exists $C>0$ such that 
\begin{equation}\label{eq-H1bound}
\int_\Om \left\vert \nabla u_\eps\right\vert^2\, dx \le C\hbox{ for all }\eps>0\hskip.1cm.
\end{equation}
Then we have the following~:

\begin{prop}[Druet \cite{DruetDMJ}]\label{prop-DMJ}
After passing to a subsequence, $\lambda_\eps\to \lambda_0$ as $\eps\to 0$, there exists a smooth solution $u_0$ of the limit equation 
\begin{equation}\label{eqzero}
\Delta u_0 = \lambda_0 f_0 u_0 e^{u_0^2} \hbox{ in }\Omega,\, u_0=0\hbox{ on }\partial \Omega
\end{equation}
and there exist $N\ge 0$ and $N$ sequences $\left(\xie\right)$ of points in $\Omega$ such that the following assertions\footnote{We assume for assertions b) to f) that $N\ge 1$.} hold~:

\smallskip a) $u_\eps \rightharpoonup u_0$ weakly in $H^1_0\left(\Omega\right)$. If $N=0$, the convergence of $u_\eps$ to $u_0$ holds in $C^2\left(\overline{\Omega}\right)$. 

\smallskip b) for any $i\in \left\{1,\dots,N\right\}$, $u_\eps\left(\xie\right)\to +\infty$ as $\eps\to 0$ and $\nabla u_\eps\left(\xie\right)=0$. 

\smallskip c) for any $i,j\in \left\{1,\dots,N\right\}$, $i\neq j$, 
$$\frac{\left\vert \xie-\xje\right\vert}{\mie}\to +\infty \hbox{ as }\eps\to 0$$
where 
$$\mie^{-2} = \lambda_\eps f_\eps\left(\xie\right) u_\eps\left(\xie\right)^2 e^{u_\eps\left(\xie\right)^2}\to +\infty \hbox{ as }\eps\to 0\hskip.1cm.$$

\smallskip d) for any $i\in \left\{1,\dots,N\right\}$, we have that 
$$u_\eps\left(\xie\right)\left(u_\eps\left(\xie+\mie x\right) - u_\eps\left(\xie\right)\right) \to U(x)= -\ln \left(1+\frac{1}{4} \vert x\vert^2\right) $$
in $C^2_{loc}\left(\rtwo\right)$. 

\smallskip e) there exists $C_1>0$ such that 
$$\lambda_\eps \left(\min_{i=1,\dots,N} \left\vert \xie-x\right\vert\right)^2 u_\eps(x)^2 e^{u_\eps(x)^2} \le C_1$$
for all $x\in \Omega$.

\smallskip f) there exists $C_2>0$ such that 
$$\left(\min_{i=1,\dots,N} \left\vert \xie-x\right\vert\right) u_\eps(x)\left\vert \nabla u_\eps(x)\right\vert \le C_2\hskip.1cm.$$

\end{prop}

\medskip {\it Proof} - Even if this result is already contained in \cite{DruetDMJ}, we shall give part of the proof here. The first reason is that it is not exactly stated in this way in \cite{DruetDMJ}. The second reason is that it is proved in greater generality in \cite{DruetDMJ} and we thus give a proof which is in some sense more readable here. 

\smallskip First, it is clear thanks to \eqref{eq-H1bound} that, up to a subsequence, $u_\eps \rightharpoonup u_0$ weakly in $H^1_0\left(\Omega\right)$ where $u_0$ is a solution of \eqref{eqzero}. If $\left\Vert u_\eps\right\Vert_\infty=O(1)$, then, by standard elliptic theory, this convergence holds in $C^2\left(\overline{\Omega}\right)$ and the proposition is true with $N=0$. Let us assume from now on that 
\begin{equation}\label{eq-prop-DMJ-1}
\sup_\Omega u_\eps \to +\infty \hbox{ as }\eps\to 0\hskip.1cm.
\end{equation}
Given $N\ge 1$ and $N$ sequences $\left(\xie\right)$ of points in $\Omega$ which verify that 
\begin{equation}\label{eq-prop-DMJ-2}
\gamma_{i,\eps} = u_\eps\left(\xie\right) \to +\infty \hbox{ as }\eps\to 0 \hbox{ and }\mie^{-2} = \lambda_\eps f_\eps\left(\xie\right)\gamma_{i,\eps}^2 e^{\gamma_{i,\eps}^2}\to +\infty \hbox{ as }\eps\to 0\hskip.1cm,
\end{equation}
we consider the following assertions~:

\smallskip {\bf $\left(P_1^N\right)$} For any $i,j\in \left\{1,\dots,N\right\}$, $i\neq j$, ${\ds \frac{\left\vert \xie-\xje\right\vert}{\mie}\to +\infty}$ as $\eps\to 0$.

\smallskip {\bf $\left(P_2^N\right)$} For any $i\in \left\{1,\dots,N\right\}$, $\nabla u_\eps\left(\xie\right)=0$ and 
$$\gamma_{i,\eps} \left(u_\eps\left(\xie+\mie x\right)-\gamma_{i,\eps}\right)\to U(x)$$
in ${\mathcal C}^2_{loc}\left(\rtwo\right)$ as $\eps \to 0$ where 
$$U(x)=-\ln \left(1+\frac{1}{4} \vert x\vert^2\right)$$
is a solution of $\Delta U = e^{2U}$ in $\rtwo$.

\smallskip {\bf $\left(P_3^N\right)$} There exists $C>0$ such that 
$$\lambda_\eps \left(\min_{i=1,\dots,N} \left\vert \xie-x\right\vert\right)^2 u_\eps(x)^2 e^{u_\eps(x)^2} \le C$$
for all $x\in \Omega$.

\medskip A first obvious remark is that 
\begin{equation}\label{eq-prop-DMJ-3}
\left(P_1^N\right) \hbox{ and }\left(P_2^N\right) \Longrightarrow \int_\Om \left\vert \nabla u_\eps\right\vert^2\, dx \ge 4\pi N +o(1)\hskip.1cm.
\end{equation}
Indeed, one has just to notice that 
$$\int_\Omega \left\vert \nabla u_\eps\right\vert^2\, dx = \lambda_\eps \int_{\Omega}f_\eps u_\eps^2 e^{u_\eps^2}\, dx\hskip.1cm,$$
that ${\mathbb D}_{\xie}\left(R\mie\right)\cap {\mathbb D}_{\xje}\left(R\mje\right)=\emptyset$ for $\eps>0$ small enough thanks to $\left(P_1^N\right)$ and that 
$$\lim_{\eps\to 0}\lambda_\eps \int_{{\mathbb D}_{\xie}\left(R\mie\right)}f_\eps u_\eps^2 e^{u_\eps^2}\, dx =\int_{{\mathbb D}_0(R)}e^{2U}\, dx \to \int_{\rtwo}e^{2U}\, dx = 4\pi\hbox{ as }R\to +\infty$$
thanks to $\left(P_2^N\right)$.

\medskip In the following, we shall say that property ${\mathcal P}_N$ holds if there are $N$ sequences $\left(\xie\right)$ of points in $\Omega$ which verify \eqref{eq-prop-DMJ-2} such that assertions $\left(P_1^N\right)$ and $\left(P_2^N\right)$ hold. 

\medskip {\sc Step 1} - {\it Property ${\mathcal P}_1$ holds.}

\medskip {\sc Proof of Step 1} - Let $x_\eps\in \Omega$ be such that 
$$u_\eps\left(x_\eps\right)=\max_\Omega u_\eps\hskip.1cm.$$
By \eqref{eq-prop-DMJ-1}, we have that 
\begin{equation}\label{eq-prop-DMJ-5}
\gamma_\eps = u_\eps\left(x_\eps\right)\to +\infty\hbox{ as }\eps\to 0\hskip.1cm.
\end{equation}
We just have to check $\left(P_2^1\right)$ since $\left(P_1^1\right)$ is empty. We clearly have that $\nabla u_\eps\left(x_\eps\right)=0$. We set\footnote{The fact that this rescaling is appropriate to understand the blow up behaviour of solutions of equation \eqref{maineq} was first discovered by Adimurthi-Struwe \cite{AdiStruwe}.}
\begin{equation}\label{eq-prop-DMJ-6}
\tilde{u}_\eps(x)= \gamma_\eps\left(u_\eps\left(x_\eps+\mu_\eps x\right)-\gamma_\eps\right)
\end{equation}
for $x\in \Omega_\eps$ where 
$$
\Omega_\eps = \left\{x\in \rtwo \hbox{ s.t. }x_\eps+\mu_\eps x\in \Omega\right\}
$$
and 
\begin{equation}\label{eq-prop-DMJ-8}
\mu_\eps^{-2} = \lambda_\eps f_\eps\left(x_\eps\right) \gamma_\eps^2 e^{\gamma_\eps^2}\hskip.1cm.
\end{equation}
It is clear that 
\begin{equation}\label{eq-prop-DMJ-9}
\mu_\eps\to 0\hbox{ as }\eps\to 0\hskip.1cm.
\end{equation}
Indeed, we can write that 
$$\lambda_\eps f_\eps u_\eps e^{u_\eps^2} \le \lambda_\eps \left(\sup_{\Omega} f_\eps\right)\gamma_\eps e^{\gamma_\eps^2} = \gamma_\eps^{-1} \frac{\sup_{\Omega} f_\eps}{f_\eps\left(x_\eps\right)}\mu_\eps^{-2}=o\left(\mu_\eps^{-2}\right)$$
thanks to \eqref{eq-cv-heps} and \eqref{eq-prop-DMJ-5}. If ever \eqref{eq-prop-DMJ-9} was false, we would have that $\left\Vert \Delta u_\eps\right\Vert_\infty \to 0$ as $\eps\to 0$ which, together with the fact that $u_\eps=0$ on $\partial \Omega$, would contradict \eqref{eq-prop-DMJ-5}. Thus \eqref{eq-prop-DMJ-9} holds.

\noindent Thanks to \eqref{eq-prop-DMJ-9}, we know that, up to a subsequence and up to a harmless rotation, 
\begin{equation}\label{eq-prop-DMJ-10}
\Omega_\eps \to \rtwo \hbox{ or } \Omega_\eps \to {\mathbb R} \times \left(-\infty, d\right)\hbox{ as }\eps\to 0
\end{equation}
where ${\ds d = \lim_{\eps\to 0} \frac{d\left(x_\eps,\partial\Omega\right)}{\mu_\eps}}$. We also have that 
\begin{equation}\label{eq-prop-DMJ-11}
\Delta \tilde{u}_\eps = \frac{f_\eps\left(x_\eps+\mu_\eps x\right)}{f_\eps\left(x_\eps\right)} \frac{u_\eps\left(x_\eps+\mu_\eps x\right)}{\gamma_\eps} e^{u_\eps\left(x_\eps+\mu_\eps x\right)^2-\gamma_\eps^2}
\end{equation}
in $\Omega_\eps$ thanks to \eqref{eqeps} and \eqref{eq-prop-DMJ-8}. Since $0\le u_\eps\le \gamma_\eps$ in $\Omega$ and thanks to \eqref{eq-cv-heps}, this leads to ${\ds \left\Vert \Delta \tilde{u}_\eps\right\Vert_{L^\infty\left(\Omega_\eps\right)}=O(1)}$. Together with the fact that $\tilde{u}_\eps\le 0=\tilde{u}_\eps(0)$ and $\tilde{u}_\eps =-\gamma_\eps^2 \to -\infty$ as $\eps\to 0$ on $\partial \Omega_\eps$, one can check that this implies that 
$$\Omega_\eps \to \rtwo \hbox{ as }\eps\to 0$$
and that
$$\tilde{u}_\eps\to U\hbox{ in }C^1_{loc}\left({\mathbb R}^2\right) \hbox{ as }\eps\to 0$$
after passing to a subsequence. We refer here the reader to \cite{AdiStruwe} or \cite{DruetDMJ} for the details of such an assertion. Moreover, we clearly have that $U\le U(0)=0$ in $\rtwo$. Noting that, as a consequence of the above convergence of $\tilde{u}_\eps$, we have that 
$$u_\eps\left(x_\eps+\mu_\eps x\right)^2-\gamma_\eps^2 \to 2U\hbox{ in }C^0_{loc}\left({\mathbb R}^2\right)\hskip.1cm,$$
one can easily pass to the limit in equation \eqref{eq-prop-DMJ-11} to obtain that 
$$\Delta U = e^{2U}\hbox{ in }\rtwo\hskip.1cm.$$
Moreover, by standard elliptic theory, one has that 
\begin{equation}\label{eq-prop-DMJ-13}
\tilde{u}_\eps\to U\hbox{ in }{\mathcal C}^2_{loc}\left(\rtwo\right)\hbox{ as }\eps\to 0\hskip.1cm.
\end{equation}
In order to apply the classification result of Chen-Li \cite{ChenLi}, we need to check that $e^{2U}\in L^1\left(\rtwo\right)$. Using \eqref{eq-prop-DMJ-6} together with \eqref{eq-cv-heps}, \eqref{eq-prop-DMJ-8} and \eqref{eq-prop-DMJ-13}, we can write that 
$$\lim_{\eps\to 0} \lambda_\eps \int_{{\mathbb D}_{x_\eps}\left(R\mu_\eps\right)} f_\eps u_\eps^2 e^{u_\eps^2}\, dx = \int_{{\mathbb D}_0(R)} e^{2U}\, dx$$
for all $R>0$. Thanks to \eqref{eqeps} and \eqref{eq-H1bound}, we know that 
$$\lambda_\eps \int_{{\mathbb D}_{x_\eps}\left(R\mu_\eps\right)} f_\eps u_\eps^2 e^{u_\eps^2}\, dx \le \lambda_\eps \int_{\Omega} f_\eps u_\eps^2 e^{u_\eps^2}\, dx=\int_\Omega \left\vert \nabla u_\eps\right\vert^2\, dx \le M$$
so that $e^{2U}\in L^1\left(\rtwo\right)$. Remembering that $U\le U(0)=0$, we thus get by \cite{ChenLi} that 
$$U\left(x\right)=-\ln \left(1+\frac{1}{4} \vert x\vert^2\right)\hskip.1cm.$$
This clearly ends the proof of Step 1. \hfill $\spadesuit$

\medskip {\sc Step 2} - {\it Assume that property ${\mathcal P}_N$ holds for some $N\ge 1$. Then either $\left(P_3^N\right)$ holds or ${\mathcal P}_{N+1}$ holds.}

\medskip {\sc Proof of Step 2} - Assume that ${\mathcal P}_N$ holds for some $N\ge 1$ (with associated sequences $\left(\xie\right)$) and that $\left(P_3^N\right)$ does not hold, meaning that 
\begin{equation}\label{eq-prop-DMJ-14}
\lambda_\eps \sup_{x\in \Omega} \left(\min_{i=1,\dots,N} \left\vert \xie-x\right\vert\right)^2 u_\eps(x)^2 e^{u_\eps(x)^2} \to +\infty\hbox{ as }\eps\to 0\hskip.1cm.
\end{equation}
We let then $y_\eps\in \Omega$ be such that 
\begin{equation}\label{eq-prop-DMJ-15}
 \left(\min_{i=1,\dots,N} \left\vert \xie-y_\eps\right\vert\right)^2 u_\eps\left(y_\eps\right)^2 e^{u_\eps\left(y_\eps\right)^2}=\sup_{x\in \Omega} \left(\min_{i=1,\dots,N} \left\vert \xie-x\right\vert\right)^2 u_\eps(x)^2 e^{u_\eps(x)^2}
\end{equation}
and we set 
$$u_\eps\left(y_\eps\right)=\hat{\gamma}_\eps\hskip.1cm.$$
Since $\Omega$ is bounded and $\left(\lambda_\eps\right)$ is bounded, see \eqref{eq-lambdaepsbounded}, we know that 
$$\hat{\gamma}_\eps\to +\infty\hbox{ as }\eps\to 0$$
thanks to \eqref{eq-prop-DMJ-14} and \eqref{eq-prop-DMJ-15}. Thanks to $\left(P_N^2\right)$, \eqref{eq-prop-DMJ-14} and \eqref{eq-prop-DMJ-15}, we also know that 
\begin{equation}\label{eq-prop-DMJ-18}
\frac{\left\vert \xie-y_\eps\right\vert}{\mie}\to +\infty\hbox{ as }\eps\to 0\hbox{ for all }1\le i\le N\hskip.1cm.
\end{equation}
We set 
$$\hat{\mu}_\eps^{-2} =\lambda_\eps f_\eps\left(y_\eps\right) \hat{\gamma}_\eps^2 e^{\hat{\gamma}_\eps^2}$$
so that, with \eqref{eq-cv-heps}, \eqref{eq-prop-DMJ-14} and \eqref{eq-prop-DMJ-15},
$$\hat{\mu}_\eps\to 0\hbox{ as }\eps\to 0$$
and 
\begin{equation}\label{eq-prop-DMJ-21}
\frac{\left\vert \xie-y_\eps\right\vert}{\hat{\mu}_\eps}\to +\infty\hbox{ as }\eps\to 0\hbox{ for all }1\le i\le N\hskip.1cm.
\end{equation}
We set now 
$$\hat{u}_\eps(x)= \hat{\gamma}_\eps\left(u_\eps\left(y_\eps+\hat{\mu}_\eps x\right)-\hat{\gamma}_\eps\right)$$
for $x\in \hat{\Omega}_\eps$ where 
$$\hat{\Omega}_\eps = \left\{x\in \rtwo \hbox{ s.t. }y_\eps+\hat{\mu}_\eps x\in \Omega\right\}\hskip.1cm.$$
We are exactly in the situation of Step 1 except for one thing~: we can not say that $\hat{u}_\eps\le 0$ in $\hat{\Omega}_\eps$. However, combining \eqref{eq-prop-DMJ-15} and \eqref{eq-prop-DMJ-21}, we can say that 
$$\hat{u}_\eps\le o(1) \hbox{ in } K\cap \hat{\Omega}_\eps$$
for all compact subset $K$ of $\rtwo$. This permits to repeat the arguments of Step 1, see \cite{DruetDMJ} for the details, to obtain that 
\begin{equation}\label{eq-prop-DMJ-25}
\hat{u}_\eps\to U\hbox{ in }{\mathcal C}^2_{loc}\left(\rtwo\right)\hbox{ as }\eps\to 0\hskip.1cm.
\end{equation}
Since $U$ has a strict local maximum at $0$, $u_\eps$ must possess, for $\eps>0$ small, a local maximum $x_{N+1,\eps}$ in $\Omega$ such that $\left\vert x_{N+1,\eps}-y_\eps\right\vert = o\left(\hat{\mu}_\eps\right)$. Then $\nabla u_\eps\left(x_{N+1,\eps}\right)=0$ and defining $\gamma_{N+1,\eps}$, $\mu_{N+1,\eps}$ with respect to this point $x_{N+1,\eps}$, it is easily checked that $\left(P_2^{N+1}\right)$ and $\left(P_1^{N+1}\right)$ hold with the sequences $\left(\xie\right)_{i=1,\dots,N+1}$ thanks to \eqref{eq-prop-DMJ-18}, \eqref{eq-prop-DMJ-21} and \eqref{eq-prop-DMJ-25}. This proves that property ${\mathcal P}_{N+1}$ holds and ends the proof of Step 2. \hfill $\spadesuit$

\medskip Starting from Step 1, and applying by induction Step 2, using \eqref{eq-H1bound} and \eqref{eq-prop-DMJ-3} to stop the process, we can easily prove the proposition except for point (f). But this point was the subject of Proposition 2 of \cite{DruetDMJ} and we refer the reader to this paper for the proof. \hfill $\diamondsuit$

\medskip The main result of Druet \cite{DruetDMJ} may be phrased as follows~:

\begin{thm}[Druet \cite{DruetDMJ}]\label{thm-DMJ}
In the framework of Proposition \ref{prop-DMJ}, there exist moreover $M\ge 0$ and $M$ sequences of points $\left(y_{i,\eps}\right)$ in $\Omega$ such that the following assertions hold after passing to a subsequence~:

\smallskip a) For any $i\in \left\{1,\dots,M\right\}$ and any $j\in \left\{1,\dots,N\right\}$,
$$\frac{\left\vert \yie-\xje\right\vert}{\mje}\to +\infty\hbox{ as }\eps\to 0\hskip.1cm.$$

\smallskip b) For any $i\in \left\{1,\dots,M\right\}$,
$$u_\eps\left(\yie\right)\left(u_\eps\left(\yie+\nie x\right) - u_\eps\left(\yie\right)\right) \to U(x)= -\ln \left(1+\frac{1}{4} \vert x\vert^2\right) $$
in $C^2_{loc}\left(\rtwo\setminus {\mathcal S}_i\right)$ where
$$\nie^{-2} = \lambda_\eps f_\eps\left(\yie\right) u_\eps\left(\yie\right)^2e^{u_\eps\left(\yie\right)^2}\to +\infty\hbox{ as }\eps\to 0$$
and 
$${\mathcal S}_i = \left\{\lim_{\eps\to 0} \frac{\xje-\yie}{\nie},\, j=1,\dots, N\right\}\bigcup \left\{\lim_{\eps\to 0} \frac{y_{k,\eps}-\yie}{\nie},\, k=1,\dots, M,\, k\neq i\right\}\hskip.1cm.$$

\smallskip c) The Dirichlet norm of $u_\eps$ is quantified by 
$$\int_\Omega \left\vert \nabla u_\eps\right\vert^2\, dx = \int_\Omega \left\vert \nabla u_0\right\vert^2\, dx + 4\pi\left(N+M\right)+o(1)\hskip.1cm.$$

\end{thm}

\medskip It is the way that the main quantification result of Druet \cite{DruetDMJ} is proved. Proposition 1 in Section 3 of \cite{DruetDMJ} corresponds to Proposition \ref{prop-DMJ} above (at the exception of f)). Then concentration points are added at the end of Section 3 of \cite{DruetDMJ}, point f) of the above proposition is proved in Section 4 of \cite{DruetDMJ} and it is proved during Sections 5 and 6 of \cite{DruetDMJ} that the quantification holds with these concentration points added. 

\medskip Let us comment on this result. First, it is clear that $u_0\not\equiv 0\Rightarrow \lambda_0>0$. Second, if $N=0$, then the convergence of $u_\eps$ to $u_0$ is strong in $H^1_0\left(\Omega\right)$ and in fact even holds in ${\mathcal C}^2\left(\overline{\Omega}\right)$. The two questions left open in this work of the first author were~:

\medskip 1. Is it possible to have $u_0\not\equiv 0$ and $N\ge 1$ together ?

\medskip 2. Are the concentration points $\left(\xie\right)$ isolated or can there be bubbles accumulation ?

\medskip These two questions can be motivated, as explained in the introduction, by the situation in low dimensions for Yamabe type equations, as studied in \cite{DruetJDG} (see also \cite{DruetENSAIOS}).  But they are also crucial in order to understand precisely the number of solutions of equation \eqref{maineq}, a question we shall address in a subsequent paper. 

\medskip Let us briefly sketch the proof of Theorem \ref{mainthm}. We start from the above results of \cite{DruetDMJ}. We shall first give some fine pointwise estimates on the sequence $\left(u_\eps\right)$ in small (but not so small) neighbourhoods of the concentration points. This will be the subject of section \ref{section-LBA}. Then we prove Theorem \ref{mainthm} in section \ref{section-proofmainthm} through a serie of claims proving successively that~: $M=0$ in Theorem \ref{thm-DMJ} above, $\lambda_0=0$ so that $u_0=0$ and, at last, the concentration points are isolated and of comparable size. All Theorem \ref{mainthm} then follows easily.

\section{Local blow up analysis} \label{section-LBA}

In this section, we get some fine estimates on sequences of solutions of equations \eqref{eqeps} in the neighbourhood of one of the concentration points $\left(\xie\right)$ of Theorem \ref{thm-DMJ}. During all this section, $C$ denotes a constant which is independant of $\eps$ or variables $x$, $y$, \dots

\medskip We let $\left(\rho_\eps\right)$ be a bounded sequence of positive real numbers (possibly converging to $0$ as $\eps\to 0$) and we consider a sequence of smooth positive functions $\left(v_\eps\right)$ which are solutions of 
\begin{equation}\label{eq-LBA-1}
\Delta v_\eps = \lambda_\eps f_\eps v_\eps e^{v_\eps^2}\hbox{ in }{\mathbb D}_0\left(\rho_\eps\right)
\end{equation}
where $\left(\lambda_\eps\right)$ is a bounded sequence of positive real numbers, $\left(f_\eps\right)$ is a sequence of smooth positive functions satisfying that there exists $C_0>0$ such that 
\begin{equation}\label{eq-LBA-0}
\frac{1}{C_0}\le f_\eps(0)\le C_0\hbox{, }\left\vert \nabla f_\eps\right\vert \le C_0 \hbox{ and }\left\vert \nabla^2 f_\eps\right\vert \le C_0  \hbox{ in }{\mathbb D}_0\left(\rho_\eps\right)\hskip.1cm.
\end{equation}
Here and in all what follows, ${\mathbb D}_x(r)$ denotes the disk of center $x$ and radius $r$. 
We assume moreover that 
\begin{equation}\label{eq-LBA-2}
\gamma_\eps = v_\eps(0)\to +\infty\hbox{ as }\eps\to 0 \hbox{ and }\nabla v_\eps(0)=0\hskip.1cm,
\end{equation}
that 
\begin{equation}\label{eq-LBA-3}
\mu_\eps^{-2} =\lambda_\eps f_\eps\left(0\right) \gamma_\eps^2 e^{\gamma_\eps^2} \to +\infty\hbox{ as }\eps\to 0\hbox{ with }\frac{\rho_\eps}{\mu_\eps}\to +\infty\hbox{ as }\eps\to 0\hskip.1cm,
\end{equation}
that 
\begin{equation}\label{eq-LBA-4}
\gamma_\eps \left(v_\eps\left(\mu_\eps x\right)-\gamma_\eps\right) \to U(x)=-\ln\left(1+\frac{1}{4}\vert x\vert^2\right)\hbox{ in }{\mathcal C}^2_{loc}\left(\rtwo\right) \hbox{ as }\eps\to 0\hskip.1cm,
\end{equation}
that there exists $C_1>0$ such that 
\begin{equation}\label{eq-LBA-5}
\lambda_\eps \vert x\vert^2 v_\eps^2 e^{v_\eps^2} \le C_1\hbox{ in }{\mathbb D}_0\left(\rho_\eps\right)
\end{equation}
and that there exists $C_2>0$ such that 
\begin{equation}\label{eq-LBA-6}
\vert x\vert \left\vert \nabla v_\eps\right\vert \le \frac{C_2}{\gamma_\eps}\hbox{ in }{\mathbb D}_0\left(\rho_\eps\right)\hskip.1cm.
\end{equation}
The aim of this section will be to compare in a suitable disk the sequence $\left(v_\eps\right)$ with the bubble $B_\eps$ defined as the radial solution in $\rtwo$ of 
\begin{equation}\label{eq-LBA-7}
\Delta B_\eps = \lambda_\eps f_\eps(0)B_\eps e^{B_\eps^2} \hbox{ with } B_\eps(0)=\gamma_\eps\hskip.1cm.
\end{equation}
Thanks to the results of Appendix A, see in particular Claims \ref{claim-appA-4} and \ref{claim-appA-5}, we know that 
\begin{equation}\label{eq-LBA-12}
\left\vert B_\eps(x)-\left(\gamma_\eps-\frac{t_\eps(x)}{\gamma_\eps}-\frac{t_\eps(x)}{\gamma_\eps^3}\right)\right\vert \le C_3\gamma_\eps^{-2} \hbox{ for }x\hbox{ s.t. }t_\eps(x)\le \gamma_\eps^2
\end{equation}
and that 
\begin{equation}\label{eq-LBA-10}
\left\vert \nabla B_\eps(x)-\gamma_\eps^{-1} \frac{2\vec{x}}{\left\vert x\right\vert^2+4\mu_\eps^2} \right\vert \le C_4\gamma_\eps^{-2} \frac{\left\vert x\right\vert}{\left\vert x\right\vert^2+\mu_\eps^2} \hbox{ for }x\hbox{ s.t. }t_\eps(x)\le \gamma_\eps^2
\end{equation}
where $C_3>0$ and $C_4>0$ are some universal constants and  
\begin{equation}\label{eq-LBA-9}
t_\eps(x) = \ln \left(1+\frac{\vert x\vert^2}{4\mu_\eps^2}\right) \hskip.1cm.
\end{equation}

\medskip We prove the following~:

\begin{prop}\label{claim-LBA-1}
We have that~:

\smallskip a) if ${\ds \overline{v}_\eps(r)=\frac{1}{2\pi r}\int_{\partial {\mathbb D}_0(r)} v_\eps\, d\sigma}$,
$$\sup_{0\le r\le \rho_\eps} \left\vert \overline{v}_\eps(r) - B_\eps(r)\right\vert =o\left(\gamma_\eps^{-1}\right)\hskip.1cm.$$
As a consequence, we have that 
$$t_\eps\left(\rho_\eps\right) \le \gamma_\eps^2 - 1 +o\left(1\right)\hskip.1cm.$$

\smallskip b) There exists $C>0$ such that 
$$\left\vert v_\eps-B_\eps\right\vert \le C\gamma_\eps^{-1}$$
and 
$$\left\vert \nabla \left(v_\eps-B_\eps\right)\right\vert \le C\gamma_\eps^{-1}\rho_\eps^{-1}$$
in ${\mathbb D}_0\left(\rho_\eps\right)$.

\smallskip c) After passing to a subsequence, 
$$\gamma_\eps \left(v_\eps\left(\rho_\eps\, \cdot\, \right)- B_\eps\left(\rho_\eps\right)\right)\to 2 \ln \frac{1}{\vert x\vert} + {\mathcal H}$$
as $\eps\to 0$ in ${\mathcal C}^1_{loc}\left({\mathbb D}_0(1)\setminus \left\{0\right\}\right)$ where ${\mathcal H}$ is some harmonic function in the unit disk satisfying 
$${\mathcal H}(0)=0\hbox{ and }  \nabla {\mathcal H}(0) = -\frac{1}{2}\lim_{\eps\to 0} \frac{\rho_\eps\nabla f_\eps(0)}{f_\eps(0)}\hskip.1cm.$$ 
\end{prop}

\medskip {\it Proof of Proposition \ref{claim-LBA-1}} - Let us first remark that we may assume without loss of generality that 
\begin{equation}\label{eq-LBA-12bis}
t_\eps\left(\rho_\eps\right)\le \gamma_\eps^2\hskip.1cm.
\end{equation}
Indeed, up to reduce $\rho_\eps$, this is the case and once a) is proved, we know that $t_\eps\left(\rho_\eps\right) \le \gamma_\eps^2 -1 +o\left(1\right)$. This will easily permit to prove that, for the original $\rho_\eps$, we had $t_\eps\left(\rho_\eps\right)\le \gamma_\eps^2$ since $t_\eps(r)\le \gamma_\eps^2 - \frac{1}{2}$ as long as $t_\eps(r)\le \gamma_\eps^2$. 

\smallskip Fix $0<\eta<1$ and let 
\begin{equation}\label{eq-LBA-13}
r_\eps=\sup\left\{r\in \left(0,\rho_\eps\right)\hbox{ s.t. } \left\vert \overline{v}_\eps(s)-B_\eps(s)\right\vert\le \frac{\eta}{\gamma_\eps}\hbox{ for all }0\le s\le r\right\}
\end{equation}
where 
$$\overline{v}_\eps(r)=\frac{1}{2\pi r} \int_{\partial {\mathbb D}_0(r)} v_\eps\, d\sigma\hskip.1cm.$$
Note that we know thanks to \eqref{eq-LBA-3} and \eqref{eq-LBA-4} that 
\begin{equation}\label{eq-LBA-14}
\frac{r_\eps}{\mu_\eps}\to +\infty \hbox{ as }\eps\to 0\hskip.1cm.
\end{equation}
We have that 
\begin{equation}\label{eq-claim-LBA-1-1}
\left\vert \overline{v}_\eps(r)-B_\eps(r)\right\vert\le \frac{\eta}{\gamma_\eps}\hbox{ for all }0\le r\le r_\eps
\end{equation}
and that 
\begin{equation}\label{eq-claim-LBA-1-2}
\left\vert \overline{v}_\eps\left(r_\eps\right)-B_\eps\left(r_\eps\right)\right\vert= \frac{\eta}{\gamma_\eps}\hbox{ if }r_\eps<\rho_\eps\hskip.1cm.
\end{equation}
We set 
\begin{equation}\label{eq-claim-LBA-1-3}
v_\eps= B_\eps+w_\eps
\end{equation}
in ${\mathbb D}_0\left(\rho_\eps\right)$. Thanks to \eqref{eq-LBA-6} and \eqref{eq-claim-LBA-1-1}, we know that 
\begin{equation}\label{eq-claim-LBA-1-4}
\left\vert w_\eps\right\vert \le \frac{\eta+\pi C_2}{\gamma_\eps}\hbox{ in }{\mathbb D}_0\left(r_\eps\right)\hskip.1cm.
\end{equation}
This clearly implies since $\vert B_\eps\vert\le \gamma_\eps$ that 
\begin{equation}\label{eq-claim-LBA-1-5}
\left\vert v_\eps^2-B_\eps^2\right\vert \le 3\left(\eta+\pi C_2\right) \hbox{ in }{\mathbb D}_0\left(r_\eps\right)\hskip.1cm.
\end{equation}
Thanks to \eqref{eq-LBA-1}, we can write that
\begincal
\Delta w_\eps &=& \lambda_\eps f_\eps v_\eps e^{v_\eps^2}- \lambda_\eps f_\eps(0) B_\eps e^{B_\eps^2}\\
&=& \lambda_\eps e^{B_\eps^2} \left(f_\eps v_\eps e^{v_\eps^2-B_\eps^2} - f_\eps(0) B_\eps\right)\\
&=& \lambda_\eps e^{B_\eps^2} \left(f_\eps w_\eps e^{v_\eps^2-B_\eps^2} +f_\eps B_\eps e^{v_\eps^2-B_\eps^2}  - f_\eps(0) B_\eps\right)\\
\fincal
in ${\mathbb D}_0\left(r_\eps\right)$ so that, using \eqref{eq-LBA-0}, \eqref{eq-claim-LBA-1-4} and \eqref{eq-claim-LBA-1-5} but also \eqref{eq-LBA-12}, we get the existence of some $C>0$ such that 
\begin{equation}\label{eq-claim-LBA-1-6}
\left\vert \Delta w_\eps\right\vert \le C \lambda_\eps f_\eps(0)\left(1+B_\eps^2\right)e^{B_\eps^2} \left\vert w_\eps\right\vert + C \lambda_\eps \vert x\vert \left(\frac{2}{\gamma_\eps}+B_\eps\right)e^{B_\eps^2} \hbox{ in }{\mathbb D}_0\left(r_\eps\right)\hskip.1cm.
\end{equation}
We let $\varphi_\eps$ be such that 
\begin{equation}\label{eq-claim-LBA-1-7}
\Delta \varphi_\eps =0 \hbox{ in }{\mathbb D}_0\left(r_\eps\right)\hbox{ and }\varphi_\eps= w_\eps \hbox{ on }\partial{\mathbb D}_0\left(r_\eps\right)\hskip.1cm.
\end{equation}
Using \eqref{eq-LBA-6} and \eqref{eq-LBA-10}, we know that 
$$\left\vert \nabla w_\eps\right\vert \le C\gamma_\eps^{-1}r_\eps^{-1} \hbox{ on }\partial{\mathbb D}_0\left(r_\eps\right)$$
for some $C>0$ so that 
\begin{equation}\label{eq-claim-LBA-1-8}
\left\Vert \nabla \varphi_\eps\right\Vert_{L^\infty\left({\mathbb D}_0\left(r_\eps\right)\right)} = O\left(\frac{1}{\gamma_\eps r_\eps}\right)\hskip.1cm.
\end{equation}
Note also that, up to a subsequence, 
\begin{equation}\label{eq-claim-LBA-1-8bis}
\gamma_\eps\varphi_\eps\left(r_\eps\,\cdot\,\right)\to\varphi_0\hbox{ in }{\mathcal C}^2_{loc}\left({\mathbb D}_0(1)\right)\hbox{ as }\eps\to 0
\end{equation}
since $\left\vert \overline{\varphi}_\eps\left(r_\eps\right)\right\vert \le \eta\gamma_\eps^{-1}$ thanks to \eqref{eq-claim-LBA-1-1} and \eqref{eq-claim-LBA-1-3}. It follows from standard elliptic theory thanks to \eqref{eq-claim-LBA-1-7}.

\medskip {\sc Step 1} - {\it There exists $C>0$ such that 
$$\left\vert \nabla \left(w_\eps-\varphi_\eps\right)(y)\right\vert \le C \left(\left\Vert \nabla w_\eps\right\Vert_{L^\infty\left({\mathbb D}_0\left(r_\eps\right)\right)}+\gamma_\eps^{-1}\right)\left( \frac{\mu_\eps}{\mu_\eps +\vert y\vert}+\gamma_\eps^{-2}\right) + C \gamma_\eps^{-2}+Cr_\eps^{-1}\gamma_\eps^{-3}$$
for all $y\in {\mathbb D}_0\left(r_\eps\right)$. }

\medskip {\sc Proof of Step 1} - Let $y_\eps\in {\mathbb D}_0\left(r_\eps\right)$. Using the Green representation formula and \eqref{eq-claim-LBA-1-6}, we can write that 
\begin{eqnarray}\label{eq-claim-LBA-1-9}
\left\vert \nabla \left(w_\eps-\varphi_\eps\right)\left(y_\eps\right)\right\vert &\le & C\lambda_\eps f_\eps(0)\int_{{\mathbb D}_0\left(r_\eps\right)} \frac{1}{\vert x-y_\eps\vert}\left(1+B_\eps(x)^2\right)e^{B_\eps(x)^2} \left\vert w_\eps(x)\right\vert\, dx\\
&&+ C \lambda_\eps \int_{{\mathbb D}_0\left(r_\eps\right)}\frac{1}{\vert x-y_\eps\vert}\vert x\vert \left(\frac{2}{\gamma_\eps}+ B_\eps(x)\right) e^{B_\eps(x)^2}\, dx \hskip.1cm.\nonumber
\end{eqnarray}
We let in the following 
\begin{equation}\label{eq-claim-LBA-1-10}
t_{1,\eps} = \frac{1}{4}\gamma_\eps^2\hbox{ and }t_{2,\eps} = \gamma_\eps^2 -\gamma_\eps \end{equation}
and we let 
\begin{eqnarray}\label{eq-claim-LBA-1-11}
\Omega_{0,\eps} &=& {\mathbb D}_0\left(r_\eps\right)\cap \left\{t_\eps(x)\le t_{1,\eps}\right\}, \nonumber\\
\Omega_{1,\eps} &= &{\mathbb D}_0\left(r_\eps\right)\cap \left\{t_{1,\eps}\le t_\eps(x)\le t_{2,\eps}\right\}\hbox{ and }\\
\Omega_{2,\eps} &=& {\mathbb D}_0\left(r_\eps\right)\cap \left\{t_\eps(x)\ge t_{2,\eps}\right\}\hskip.1cm.\nonumber
\end{eqnarray}
We also set, for $i=0,1,2$, 
\begin{equation}\label{eq-claim-LBA-1-12}
I_{i,\eps} = \lambda_\eps f_\eps(0)\int_{\Omega_{i,\eps}}\frac{1}{\vert x-y_\eps\vert}\left(1+B_\eps(x)^2\right)e^{B_\eps(x)^2} \left\vert w_\eps(x)\right\vert\, dx
\end{equation}
and 
\begin{equation}\label{eq-claim-LBA-1-13}
J_{i,\eps} = \lambda_\eps \int_{\Omega_{i,\eps}}\frac{1}{\vert x-y_\eps\vert}\vert x\vert \left(\frac{2}{\gamma_\eps}+ B_\eps(x)\right) e^{B_\eps(x)^2}\, dx\hskip.1cm.
\end{equation}

\medskip {\sc Case 1} - We assume first that $\vert y_\eps\vert = O\left(\mu_\eps\right)$. Since $w_\eps(0)=0$ and using \eqref{eq-LBA-12}, we can write that 
$$I_{0,\eps}\le C\lambda_\eps f_\eps(0) \gamma_\eps^2 e^{\gamma_\eps^2}\left\Vert \nabla w_\eps\right\Vert_{L^\infty\left(\Omega_{0,\eps}\right)}
\int_{\Omega_{0,\eps}} \frac{\vert x\vert}{\vert x-y_\eps\vert}e^{\frac{t_\eps(x)^2}{\gamma_\eps^2}-2t_\eps(x)}\, dx\hskip.1cm.$$
Thanks to \eqref{eq-LBA-3}, we can rewrite this as 
$$I_{0,\eps}\le C\mu_\eps^{-2}\left\Vert \nabla w_\eps\right\Vert_{L^\infty\left(\Omega_{0,\eps}\right)}
\int_{\Omega_{0,\eps}} \frac{\vert x\vert}{\vert x-y_\eps\vert}e^{\frac{t_\eps(x)^2}{\gamma_\eps^2}-2t_\eps(x)}\, dx\hskip.1cm.$$
Since 
$$\frac{t_\eps(x)^2}{\gamma_\eps^2}-2t_\eps(x)\le -\frac{7}{4}t_\eps(x)\hbox{ in }\Omega_{0,\eps}\hskip.1cm,$$
this leads to 
\begincal
I_{0,\eps}&\le & C\mu_\eps^{-2}\left\Vert \nabla w_\eps\right\Vert_{L^\infty\left(\Omega_{0,\eps}\right)}
\int_{\Omega_{0,\eps}} \frac{\vert x\vert}{\vert x-y_\eps\vert} \left(1+\frac{\vert x\vert^2}{4\mu_\eps^2}\right)^{-\frac{7}{4}}\, dx\\
&\le & C\left\Vert \nabla w_\eps\right\Vert_{L^\infty\left(\Omega_{0,\eps}\right)}\int_{{\mathbb R}^2} 
\frac{\vert x\vert}{\left\vert x-\frac{y_\eps}{\mu_\eps}\right\vert} \left(1+\frac{\vert x\vert^2}{4}\right)^{-\frac{7}{4}}\, dx\hskip.1cm.
\fincal
Since $\left\vert y_\eps\right\vert =O\left(\mu_\eps\right)$, we obtain by the dominated convergence theorem that
\begin{equation}\label{eq-claim-LBA-1-14}
I_{0,\eps}=O\left(\left\Vert \nabla w_\eps\right\Vert_{L^\infty\left(\Omega_{0,\eps}\right)}\right)\hskip.1cm.
\end{equation}
In $\Omega_{1,\eps}$, we have that $\vert x\vert \le \left(1+o(1)\right)\left\vert x-y_\eps\right\vert$ since $\vert y_\eps\vert = O\left(\mu_\eps\right)$ so that, we can write, as above
\begincal
I_{1,\eps} &\le & C\mu_\eps^{-2}\left\Vert \nabla w_\eps\right\Vert_{L^\infty\left(\Omega_{1,\eps}\right)}
\int_{\Omega_{1,\eps}}e^{\frac{t_\eps(x)^2}{\gamma_\eps^2}-2t_\eps(x)}\, dx\\
&\le & C\left\Vert \nabla w_\eps\right\Vert_{L^\infty\left(\Omega_{1,\eps}\right)}\int_{t_{1,\eps}}^{t_{2,\eps}} e^{\frac{t^2}{\gamma_\eps^2}-t}\, dt
\fincal
by the change of variables $t=\ln\left(1+\frac{\vert x\vert^2}{4\mu_\eps^2}\right)$. Since 
$$\frac{t^2}{\gamma_\eps^2}-t \le -\frac{t}{\gamma_\eps}\le -\frac{1}{4}\gamma_\eps$$
for $\frac{1}{4}\gamma_\eps^2 = t_{1,\eps}\le t\le t_{2,\eps}=\gamma_\eps^2-\gamma_\eps$, we immediately get that 
\begin{equation}\label{eq-claim-LBA-1-15}
I_{1,\eps} \le C \left\Vert \nabla w_\eps\right\Vert_{L^\infty\left(\Omega_{1,\eps}\right)} \gamma_\eps^2 e^{-\frac{1}{4}\gamma_\eps}\hskip.1cm.
\end{equation}
In $\Omega_{2,\eps}$, we have that $B_\eps = O\left(1\right)$ thanks to \eqref{eq-LBA-12} so that, using \eqref{eq-LBA-0} and \eqref{eq-claim-LBA-1-4} , we can write that 
$$I_{2,\eps}\le C\lambda_\eps\gamma_\eps^{-1} \int_{{\mathbb D}_0\left(r_\eps\right)} \frac{1}{\left\vert x-y_\eps\right\vert}\, dx $$
so that 
\begin{equation}\label{eq-claim-LBA-1-16}
I_{2,\eps}\le C \lambda_\eps r_\eps\gamma_\eps^{-1}\hskip.1cm.
\end{equation}
Now we notice that $t_\eps\left(r_\eps\right)\le \gamma_\eps^2$ implies that 
$$\frac{r_\eps^2}{\mu_\eps^2}\le 4e^{\gamma_\eps^2}\hskip.1cm.$$
Using \eqref{eq-LBA-0} and \eqref{eq-LBA-3}, this gives that 
\begin{equation}\label{eq-claim-LBA-1-18bis}
\lambda_\eps r_\eps^2 \le \frac{C}{\gamma_\eps^2}\hskip.1cm.
\end{equation}
Thus we get that 
\begin{equation}\label{eq-claim-LBA-1-16bis}
I_{2,\eps}\le C r_\eps^{-1}\gamma_\eps^{-3}\hskip.1cm.
\end{equation}

For the second set of integrals, things are similar and easier. We write that 
$$J_{0,\eps} \le C\mu_\eps^{-2}\gamma_{\eps}^{-1} \int_{\Omega_{0,\eps}} \frac{\vert x\vert}{\vert x-y_\eps\vert}e^{\frac{t_\eps(x)^2}{\gamma_\eps^2}-2t_\eps(x)}\, dx$$
so that, see above,
\begin{equation}\label{eq-claim-LBA-1-17}
J_{0,\eps} \le C\gamma_\eps^{-1}\hskip.1cm.
\end{equation}
We also have that
\begin{equation}\label{eq-claim-LBA-1-18}
J_{1,\eps} \le C  \gamma_\eps e^{-\frac{1}{4}\gamma_\eps}
\end{equation}
in the same way than above. At last, for $J_{2,\eps}$, we write that
$$J_{2,\eps} \le C\lambda_\eps r_\eps\int_{{\mathbb D}_0\left(r_\eps\right)} \frac{1}{\left\vert x-y_\eps\right\vert}\, dx \le C\lambda_\eps r_\eps^2\hskip.1cm.$$
Thus we have thanks to \eqref{eq-claim-LBA-1-18bis}  that 
\begin{equation}\label{eq-claim-LBA-1-19}
J_{2,\eps} \le C\gamma_\eps^{-2}\hskip.1cm.
\end{equation}

Summarizing, we obtain in this first case, coming back to \eqref{eq-claim-LBA-1-9} with \eqref{eq-claim-LBA-1-14}, \eqref{eq-claim-LBA-1-15}, \eqref{eq-claim-LBA-1-16bis}, \eqref{eq-claim-LBA-1-17}, \eqref{eq-claim-LBA-1-18} and \eqref{eq-claim-LBA-1-19}, that 
\begin{equation}\label{eq-claim-LBA-1-20}
\left\vert \nabla \left(w_\eps-\varphi_\eps\right)\left(y_\eps\right)\right\vert \le C\left\Vert \nabla w_\eps\right\Vert_{L^\infty\left({\mathbb D}_0\left(r_\eps\right)\right)} + C\gamma_\eps^{-1}+Cr_\eps^{-1}\gamma_\eps^{-3}\hskip.1cm.
\end{equation}

\medskip {\sc Case 2} - We assume now that $\frac{\left\vert y_\eps\right\vert}{\mu_\eps}\to +\infty$ as $\eps\to 0$. 

\medskip We follow the lines of the first case to estimate most of the integrals. Thus we only emphasize on the changes. First, we write that 
\begincal
I_{0,\eps} &\le &  C\mu_\eps^{-2}\left\Vert \nabla w_\eps\right\Vert_{L^\infty\left(\Omega_{0,\eps}\right)}
\int_{\Omega_{0,\eps}} \frac{\vert x\vert}{\vert x-y_\eps\vert}\left(1+\frac{\vert x\vert^2}{4\mu_\eps^2}\right)^{-\frac{7}{4}}\, dx\\
&\le & C\left\Vert \nabla w_\eps\right\Vert_{L^\infty\left(\Omega_{0,\eps}\right)}\int_{{\mathbb R}^2} \frac{\vert x\vert}{\vert x-\mu_\eps^{-1}y_\eps\vert}\left(1+\frac{\vert x\vert^2}{4}\right)^{-\frac{7}{4}}\, dx \hskip.1cm.
\fincal
Now we can write that 
\begincal
&&\int_{{\mathbb R}^2} \frac{\vert x\vert}{\vert x-\mu_\eps^{-1}y_\eps\vert}\left(1+\frac{\vert x\vert^2}{4}\right)^{-\frac{7}{4}}\, dx \\
&&\quad = \left(\frac{\left\vert y_\eps\right\vert}{\mu_\eps}\right)^{-\frac{3}{2}}\int_{{\mathbb R}^2}\frac{\vert x\vert}{\vert x-\left\vert y_\eps\right\vert^{-1}y_\eps\vert}\left(\frac{\mu_\eps^2}{\left\vert y_\eps\right\vert^2}+\frac{\vert x\vert^2}{4}\right)^{-\frac{7}{4}}\, dx\\
&&\quad\le  C\left(\frac{\left\vert y_\eps\right\vert}{\mu_\eps}\right)^{-\frac{3}{2}} 
+ 2\left(\frac{\left\vert y_\eps\right\vert}{\mu_\eps}\right)^{-\frac{3}{2}}\int_{{\mathbb D}_0\left(\frac{1}{2}\right)}\vert x\vert\left(\frac{\mu_\eps^2}{\left\vert y_\eps\right\vert^2}+\frac{\vert x\vert^2}{4}\right)^{-\frac{7}{4}}\, dx\\
&&\quad\le  C\left(\frac{\left\vert y_\eps\right\vert}{\mu_\eps}\right)^{-\frac{3}{2}} 
+ 2 \frac{\mu_\eps}{\left\vert y_\eps\right\vert}\int_{{\mathbb R}^2}\vert x\vert\left(1+\frac{\vert x\vert^2}{4}\right)^{-\frac{7}{4}}\, dx
\fincal
so that 
\begin{equation}\label{eq-claim-LBA-1-21}
I_{0,\eps} \le C\left\Vert \nabla w_\eps\right\Vert_{L^\infty\left(\Omega_{0,\eps}\right)} \frac{\mu_\eps}{\left\vert y_\eps\right\vert}\hskip.1cm.
\end{equation}
Let us write once again that 
$$I_{1,\eps} \le C\mu_\eps^{-2}\left\Vert \nabla w_\eps\right\Vert_{L^\infty\left(\Omega_{1,\eps}\right)}
\int_{\Omega_{1,\eps}}\frac{\vert x\vert}{\left\vert x-y_\eps\right\vert}e^{\frac{t_\eps(x)^2}{\gamma_\eps^2}-2t_\eps(x)}\, dx\hskip.1cm.$$
Let us split this integral into two parts. First, 
\begincal
\mu_\eps^{-2}\int_{\Omega_{1,\eps}\setminus {\mathbb D}_{y_\eps}\left(\frac{1}{2}\left\vert y_\eps\right\vert\right)} \frac{\vert x\vert}{\left\vert x-y_\eps\right\vert}e^{\frac{t_\eps(x)^2}{\gamma_\eps^2}-2t_\eps(x)}\, dx &\le& 
3 \mu_\eps^{-2}\int_{\Omega_{1,\eps}}e^{\frac{t_\eps(x)^2}{\gamma_\eps^2}-2t_\eps(x)}\, dx\\
&\le &C  \int_{t_{1,\eps}}^{t_{2,\eps}} e^{\frac{t^2}{\gamma_\eps^2}-t}\, dt \\
 &\le & C \gamma_\eps^2 e^{-\frac{1}{4}\gamma_\eps}
\fincal
as in Case 1. Second, 
\begincal
\mu_\eps^{-2}\int_{\Omega_{1,\eps}\cap {\mathbb D}_{y_\eps}\left(\frac{1}{2}\left\vert y_\eps\right\vert\right)} \frac{\vert x\vert}{\left\vert x-y_\eps\right\vert}e^{\frac{t_\eps(x)^2}{\gamma_\eps^2}-2t_\eps(x)}\, dx
&\le & \frac{3}{2}\mu_\eps^{-2} \left\vert y_\eps\right\vert e^{\frac{s_\eps^2}{\gamma_\eps^2}-2s_\eps}\int_{ {\mathbb D}_{y_\eps}\left(\frac{1}{2}\left\vert y_\eps\right\vert\right)} \frac{1}{\left\vert x-y_\eps\right\vert}\, dx
\fincal
where 
$$s_\eps=t_\eps\left(\frac{y_\eps}{2}\right)\hskip.1cm.$$
Thus we have that 
\begincal
\mu_\eps^{-2}\int_{\Omega_{1,\eps}\cap {\mathbb D}_{y_\eps}\left(\frac{1}{2}\left\vert y_\eps\right\vert\right)} \frac{\vert x\vert}{\left\vert x-y_\eps\right\vert}e^{\frac{t_\eps(x)^2}{\gamma_\eps^2}-2t_\eps(x)}\, dx 
&\le & C \frac{\left\vert y_\eps\right\vert^2}{\mu_\eps^2}e^{\frac{s_\eps^2}{\gamma_\eps^2}-2s_\eps}
\fincal
Note that ${\ds \Omega_{1,\eps}\cap {\mathbb D}_{y_\eps}\left(\frac{1}{2}\left\vert y_\eps\right\vert\right) =\emptyset}$ if 
$$t_\eps\left(\frac{3}{2}\left\vert y_\eps\right\vert\right)=\ln \left(1+\frac{9\left\vert y_\eps\right\vert^2}{16\mu_\eps^2}\right)\le t_{1,\eps}=\frac{1}{4}\gamma_\eps^2$$ 
so that we may assume that 
$$\ln \left(1+\frac{9\left\vert y_\eps\right\vert^2}{16\mu_\eps^2}\right)> \frac{1}{4}\gamma_\eps^2\hskip.1cm.$$
Thus 
$$s_\eps = \ln \left(1+\frac{\left\vert y_\eps\right\vert^2}{16\mu_\eps^2}\right)\ge \frac{1}{4}\gamma_\eps^2 -\ln 9\hskip.1cm.$$
It is also clear that if ${\ds \Omega_{1,\eps}\cap {\mathbb D}_{y_\eps}\left(\frac{1}{2}\left\vert y_\eps\right\vert\right) \neq\emptyset}$, $s_\eps\le t_{2,\eps}=\gamma_\eps^2-\gamma_\eps$. Thus we have that  
$$\frac{s_\eps^2}{\gamma_\eps^2}-s_\eps\le -\frac{1}{4}\gamma_\eps +O(1)\hskip.1cm.$$
We deduce that, if not zero, 
\begincal
\mu_\eps^{-2}\int_{\Omega_{1,\eps}\cap {\mathbb D}_{y_\eps}\left(\frac{1}{2}\left\vert y_\eps\right\vert\right)} \frac{\vert x\vert}{\left\vert x-y_\eps\right\vert}e^{\frac{t_\eps(x)^2}{\gamma_\eps^2}-2t_\eps(x)}\, dx &\le& C \frac{\left\vert y_\eps\right\vert^2}{\mu_\eps^2}e^{-\frac{1}{4}\gamma_\eps} e^{-s_\eps}\\
&\le & C \frac{\left\vert y_\eps\right\vert^2}{\mu_\eps^2}e^{-\frac{1}{4}\gamma_\eps} \left(1+\frac{\left\vert y_\eps\right\vert^2}{16\mu_\eps^2}\right)^{-1}\\
&\le & Ce^{-\frac{1}{4}\gamma_\eps}\hskip.1cm.
\fincal 
Thus we arrive to 
\begin{equation}\label{eq-claim-LBA-1-22}
I_{1,\eps} \le C\gamma_\eps^2 e^{-\frac{1}{4}\gamma_\eps}\left\Vert \nabla w_\eps\right\Vert_{L^\infty\left(\Omega_{1,\eps}\right)}\hskip.1cm.
\end{equation}
At last, for $I_{2,\eps}$, we have nothing to change to get that 
\begin{equation}\label{eq-claim-LBA-1-23}
I_{2,\eps} \le Cr_\eps^{-1}\gamma_\eps^{-3}\hskip.1cm.
\end{equation}
For $J_{0,\eps}$, $J_{1,\eps}$ and $J_{2,\eps}$, we proceed as above or as in Case 1 to get that 
$$J_{0,\eps}\le C \gamma_\eps^{-1} \frac{\mu_\eps}{\left\vert y_\eps\right\vert}, \,\, J_{1,\eps}\le C \gamma_\eps e^{-\frac{1}{4}\gamma_\eps} \hbox{ and }J_{2,\eps}\le \frac{C}{\gamma_\eps^2}\hskip.1cm.$$
Thus, in this second case, we obtain coming back to \eqref{eq-claim-LBA-1-9} with \eqref{eq-claim-LBA-1-21}, \eqref{eq-claim-LBA-1-22}, \eqref{eq-claim-LBA-1-23} and these last estimates that 
\begin{equation}\label{eq-claim-LBA-1-24}
\left\vert \nabla \left(w_\eps-\varphi_\eps\right)\left(y_\eps\right)\right\vert \le C \left(\left\Vert \nabla w_\eps\right\Vert_{L^\infty\left({\mathbb D}_0\left(r_\eps\right)\right)}+\gamma_\eps^{-1}\right) \left(\frac{\mu_\eps}{\left\vert y_\eps\right\vert}+\gamma_\eps^2e^{-\frac{1}{4}\gamma_\eps}\right) + C\gamma_\eps^{-2}+Cr_\eps^{-1}\gamma_\eps^{-3}\hskip.1cm.
\end{equation}

\medskip The study of these two cases clearly permits to conclude Step 1. \hfill $\spadesuit$

\medskip {\sc Step 2} - {\it We have that 
$$\left\Vert \nabla \left(w_\eps-\varphi_\eps\right)\right\Vert_{L^\infty\left({\mathbb D}_0\left(r_\eps\right)\right)} = o\left(\gamma_\eps^{-1}r_\eps^{-1}\right)+O\left(\gamma_\eps^{-1}\right)$$
and that 
$$\left\Vert w_\eps-\varphi_\eps\right\Vert_{L^\infty\left({\mathbb D}_0\left(r_\eps\right)\right)} = o\left(\gamma_\eps^{-1}\right)\hskip.1cm.$$
Moreover, if $r_\eps\not\to 0$ as $\eps\to 0$, we have that 
$$\lim_{\eps\to 0} \frac{\nabla f_\eps(0)}{f_\eps(0)}= -2\left(\lim_{\eps\to 0}\frac{1}{r_\eps}\right) \nabla \varphi_0(0)\hskip.1cm.$$
}

\medskip {\sc Proof of Step 2} - Let $y_\eps\in {\mathbb D}_0\left(r_\eps\right)$ be such that 
\begin{equation}\label{eq-claim-LBA-1-25}
\left\vert \nabla \left(w_\eps-\varphi_\eps\right)\left(y_\eps\right)\right\vert = \left\Vert \nabla \left(w_\eps-\varphi_\eps\right)\right\Vert_{L^\infty\left({\mathbb D}_0\left(r_\eps\right)\right)}
\end{equation}
and let us assume that 
\begin{equation}\label{eq-claim-LBA-1-26}
\alpha_\eps = \left\vert \nabla \left(w_\eps-\varphi_\eps\right)\left(y_\eps\right)\right\vert \ge \frac{\delta}{r_\eps\gamma_\eps}+\frac{1}{\delta\gamma_\eps}
\end{equation}
for some $\delta>0$. Thanks to \eqref{eq-claim-LBA-1-8}, we have that 
\begin{equation}\label{eq-claim-LBA-1-26bis}
\left\Vert \nabla w_\eps\right\Vert_{L^\infty\left({\mathbb D}_0\left(r_\eps\right)\right)} \le \alpha_\eps +Cr_\eps^{-1}\gamma_\eps^{-1}\le \alpha_\eps \left(1+\frac{C}{\delta}\right)\hskip.1cm.
\end{equation}
Applying Step 1 to this sequence $\left(y_\eps\right)$, we get thanks to \eqref{eq-claim-LBA-1-25}, \eqref{eq-claim-LBA-1-26} and \eqref{eq-claim-LBA-1-26bis} that 
$$\left(\frac{1}{\delta}+\frac{\delta}{r_\eps}\right)\gamma_\eps^{-1}\le \alpha_\eps = \left\vert \nabla \left(w_\eps-\varphi_\eps\right)\left(y_\eps\right)\right\vert \le C_\delta \alpha_\eps\left( \frac{\mu_\eps}{\mu_\eps +\left\vert y_\eps\right\vert}+\gamma_\eps^{-2}\right) + C \gamma_\eps^{-2}+Cr_\eps^{-1}\gamma_\eps^{-3}\hskip.1cm.$$
This proves that 
\begin{equation}\label{eq-claim-LBA-1-27} 
\frac{y_\eps}{\mu_\eps}\to y_0\in {\mathbb R}^2 \hbox{ as }\eps\to 0
\end{equation}
after passing to a subsequence and, thanks to Step 1 and \eqref{eq-claim-LBA-1-26}, that 
\begin{equation}\label{eq-claim-LBA-1-28}
\left\vert \nabla \left(\tilde{w}_\eps-\tilde{\varphi}_\eps\right)(x)\right\vert \le \frac{C_\delta}{1+\vert x\vert}+o(1)\hbox{ for all }x\in {\mathbb R}^2
\end{equation}
where $C_\delta$ depends only on $\delta$ and 
\begin{equation}\label{eq-claim-LBA-1-29}
\tilde{w}_\eps(x) = \frac{1}{\mu_\eps \alpha_\eps} w_\eps\left(\mu_\eps x\right), \, \tilde{\varphi}_\eps(x)=\frac{1}{\mu_\eps \alpha_\eps} \varphi_\eps\left(\mu_\eps x\right)\hskip.1cm.
\end{equation}
We know that 
\begin{equation}\label{eq-claim-LBA-1-30}
\tilde{w}_\eps(0)=0,\, \nabla \tilde{w}_\eps(0)=0\hbox{ and }\left\vert \nabla\left(\tilde{w}_\eps-\tilde{\varphi}_\eps\right)\left(\frac{y_\eps}{\mu_\eps}\right)\right\vert =1\hskip.1cm.
\end{equation}
We also know thanks to \eqref{eq-claim-LBA-1-8bis} and \eqref{eq-claim-LBA-1-26} that, after passing to a subsequence, 
\begin{equation}\label{eq-claim-LBA-1-30bis}
\nabla \tilde{\varphi}_\eps(x)\to \left(\lim_{\eps\to 0}\frac{1}{\gamma_\eps r_\eps\alpha_\eps}\right) \nabla \varphi_0(0)=\vec{A}\hbox{ in }{\mathcal C}^1_{loc}\left({\mathbb R}^2\right)\hbox{ as }\eps\to 0\hskip.1cm.
\end{equation}
Using \eqref{eq-claim-LBA-1-6}, we can write that 
$$\left\vert \Delta\tilde{w}_\eps\right\vert \le C\lambda_\eps \mu_\eps^2 f_\eps(0) \left(1+B_\eps\left(\mu_\eps x\right)^2\right) e^{B_\eps\left(\mu_\eps x\right)^2} \left\vert \tilde{w}_\eps\right\vert + C \lambda_\eps \alpha_\eps^{-1}\mu_\eps^2 \vert x\vert \left(\frac{2}{\gamma_\eps}+B_\eps\left(\mu_\eps x\right) \right)e^{B_\eps\left(\mu_\eps x\right)^2}\hskip.1cm.$$
Noting thanks to \eqref{eq-claim-LBA-1-28}, \eqref{eq-claim-LBA-1-30} and \eqref{eq-claim-LBA-1-30bis} that 
$$\left\vert \tilde{w}_\eps(x)\right\vert \le C_\delta\ln \left(1+\vert x\vert\right)+\left\vert \vec{A}\right\vert \vert x\vert +o\left(\vert x\vert\right)$$
and is thus uniformly bounded on any compact subset of ${\mathbb R}^2$, we easily deduce from the above estimate together with the definition \eqref{eq-LBA-3} of $\mu_\eps$ and \eqref{eq-claim-LBA-1-26} that 
$\left( \Delta\tilde{w}_\eps\right)$ is uniformly bounded in any compact subset of ${\mathbb R}^2$.
Thus, by standard elliptic theory, we have that, after passing to a subsequence,  
\begin{equation}\label{eq-claim-LBA-1-31}
\tilde{w}_\eps\to w_0\hbox{ in }C^{1,\eta}_{loc}\left({\mathbb R}^2\right) \hbox{ as }\eps\to 0\hskip.1cm.
\end{equation}
Moreover, we have thanks to \eqref{eq-claim-LBA-1-27}, \eqref{eq-claim-LBA-1-28}, \eqref{eq-claim-LBA-1-30} and \eqref{eq-claim-LBA-1-31} that 
\begin{equation}\label{eq-claim-LBA-1-32}
w_0(0)=0, \, \nabla w_0(0)=0,\, \left\vert\nabla w_0\left(y_0\right)-\vec{A}\right\vert =1\hbox{ and }\left\vert \nabla w_0(x)-\vec{A}\right\vert \le \frac{C_\delta}{1+\vert x\vert}\hbox{ in }{\mathbb R}^2\hskip.1cm.
\end{equation}
Thus $w_0\not\equiv 0$. Since we know that $\gamma_\eps w_\eps\left(\mu_\eps x\right)\to 0$ in ${\mathcal C}^1_{loc}\left({\mathbb R}^2\right)$ as $\eps\to 0$ thanks to \eqref{eq-LBA-4}, we deduce that 
\begin{equation}\label{eq-claim-LBA-1-33}
\gamma_\eps \mu_\eps \alpha_\eps \to 0 \hbox{ as }\eps\to 0\hskip.1cm.
\end{equation}
Thanks to \eqref{eq-LBA-1}, \eqref{eq-LBA-3}, \eqref{eq-LBA-7}, \eqref{eq-claim-LBA-1-3} and \eqref{eq-claim-LBA-1-29}, we can write that 
\begincal
\Delta \tilde{w}_\eps (x) &=& \frac{1}{\alpha_\eps} \mu_\eps \lambda_\eps \left(f_\eps\left(\mu_\eps x\right) \left(B_\eps\left(\mu_\eps x\right)+w_\eps\left(\mu_\eps x\right)\right) e^{\left(B_\eps\left(\mu_\eps x\right)+w_\eps\left(\mu_\eps x\right)\right)^2}\right.\\
&&\left. - f_\eps(0) B_\eps\left(\mu_\eps x\right) e^{B_\eps\left(\mu_\eps x\right)^2}\right)\hskip.1cm.\\
&=& \frac{B_\eps\left(\mu_\eps x\right)}{\gamma_\eps} e^{B_\eps\left(\mu_\eps x\right)^2-\gamma_\eps^2} \frac{1}{\alpha_\eps \mu_\eps \gamma_\eps} 
\left(\frac{f_\eps\left(\mu_\eps x\right)}{f_\eps(0)} e^{2B_\eps\left(\mu_\eps x\right) w_\eps\left(\mu_\eps x\right)+w_\eps\left(\mu_\eps x\right)^2}-1\right)\\
&& + \gamma_\eps^{-2} \frac{f_\eps\left(\mu_\eps x\right)}{f_\eps(0)} \tilde{w}_\eps(x)e^{B_\eps\left(\mu_\eps x\right)^2-\gamma_\eps^2+2B_\eps\left(\mu_\eps x\right) w_\eps\left(\mu_\eps x\right)+w_\eps\left(\mu_\eps x\right)^2}\hskip.1cm.
\fincal
Let us write now that 
$$\gamma_\eps^{-1} B_\eps\left(\mu_\eps x\right) \to 1 \hbox{ in }{\mathcal C}^0_{loc}\left({\mathbb R}^2\right)\hbox{ as }\eps\to 0\hskip.1cm,$$
that 
$$B_\eps\left(\mu_\eps x\right)^2-\gamma_\eps^2 \to 2 U(x) \hbox{ in }{\mathcal C}^0_{loc}\left({\mathbb R}^2\right)\hbox{ as }\eps\to 0$$
where 
$$U(x)=-\ln \left(1+\frac{\vert x\vert^2}{4}\right)$$
thanks to \eqref{eq-LBA-12} and \eqref{eq-LBA-9}. We can also write that 
$$\frac{f_\eps\left(\mu_\eps x\right)}{f_\eps(0)} = 1 + f_\eps(0)^{-1}\mu_\eps x^i\partial_i f_\eps(0) + O\left(\mu_\eps^2 \vert x\vert^2\right)$$
thanks to \eqref{eq-LBA-0} and that 
$$2B_\eps\left(\mu_\eps x\right) w_\eps\left(\mu_\eps x\right)+w_\eps\left(\mu_\eps x\right)^2 = 2\mu_\eps \alpha_\eps\gamma_\eps \left(w_0+o(1)\right) = o(1)$$
thanks to \eqref{eq-claim-LBA-1-31} and \eqref{eq-claim-LBA-1-33}. Thus we can write that 
\begincal 
&&\frac{f_\eps\left(\mu_\eps x\right)}{f_\eps(0)} e^{2B_\eps\left(\mu_\eps x\right) w_\eps\left(\mu_\eps x\right)+w_\eps\left(\mu_\eps x\right)^2}-1\\
&&\quad = 2\mu_\eps \alpha_\eps\gamma_\eps w_0 + \mu_\eps f_\eps(0)^{-1} x^i\partial_i f_\eps(0) + o\left(\mu_\eps\alpha_\eps\gamma_\eps\right)
\fincal
Thus we obtain that 
$$\Delta \tilde{w}_\eps (x) =  e^{2U(x)}\left( 2w_0(x) + \frac{1}{\alpha_\eps\gamma_\eps} f_\eps(0)^{-1} x^i\partial_i f_\eps(0)\right)+o(1)\hskip.1cm.$$
Thanks to \eqref{eq-claim-LBA-1-26}, we know that, after passing to a subsequence,
\begin{equation}\label{eq-claim-LBA-1-34}
\frac{1}{\alpha_\eps\gamma_\eps}\frac{\partial_i f_\eps(0)}{f_\eps(0)} \to X_i \hbox{ as }\eps\to 0\hskip.1cm.
\end{equation}
Note that we have, again thanks to \eqref{eq-claim-LBA-1-26}, that 
\begin{equation}\label{eq-claim-LBA-1-35}
\vec{X} = 0 \hbox{ if }r_\eps\to 0\hbox{ as }\eps\to 0\hskip.1cm.
\end{equation}
Then we can write that
$$\Delta \tilde{w}_\eps (x) =  e^{2U(x)}\left( 2w_0(x) + X_i x^i\right)+o(1)$$
so that 
\begin{equation}\label{eq-claim-LBA-1-36}
\Delta w_0 = e^{2U}\left(2w_0 + X_i x^i\right)\hbox{ in }{\mathbb R}^2\hskip.1cm.
\end{equation}
Now, thanks to \cite{Chen-Lin}, lemma 2.3 or \cite{Laurain}, lemma C.1, we know that the only solution of this equation satisfying \eqref{eq-claim-LBA-1-32} is 
\begin{equation}\label{eq-claim-LBA-1-37}
w_0(x)=\frac{\vert x\vert^2}{4+\vert x\vert^2} A^i x_i
\end{equation}
and, moreover, we must have 
\begin{equation}\label{eq-claim-LBA-1-38}
\vec{A} = -\frac{1}{2}\vec{X} \hskip.1cm.
\end{equation}
Since $w_0\not\equiv 0$, we must have $\vec{A}\neq 0$ and thus $\vec{X}\neq 0$. 

\smallskip This permits to prove the step. Indeed, if $r_\eps\to 0$, then we have that $\vec{X}=0$ by \eqref{eq-claim-LBA-1-35}, which is a contradiction. Thus, if $r_\eps\to 0$, we get that \eqref{eq-claim-LBA-1-26} is impossible so that $\alpha_\eps = o\left(\gamma_\eps^{-1}r_\eps^{-1}\right)$ in this case. This proves the first estimate of the step in the case $r_\eps\to 0$ as $\eps\to 0$. If $r_\eps\not\to 0$ as $\eps\to 0$, we know thanks to the fact that $\vec{X}\neq 0$ and to \eqref{eq-claim-LBA-1-34} that $\alpha_\eps =O\left(\gamma_\eps^{-1}\right)$ if \eqref{eq-claim-LBA-1-26} holds and if it does not hold, we again have that $\alpha_\eps =O\left(\gamma_\eps^{-1}\right)$. Thus we also have that the first estimate of the step holds if $r_\eps\not\to 0$ as $\eps\to 0$. Moreover, in this second case, we know that 
$$\lim_{\eps\to 0} \frac{\partial_i f_\eps(0)}{f_\eps(0)}=-2\left(\lim_{\eps\to 0}\frac{1}{r_\eps}\right) \nabla \varphi_0(0)$$
thanks to \eqref{eq-claim-LBA-1-30bis}, \eqref{eq-claim-LBA-1-34} and \eqref{eq-claim-LBA-1-38}. This proves the last part of the step. 

\smallskip It remains to notice that the second estimate of the step is a simple consequence of the first. Indeed, coming back to the estimate of Step 1 with the estimate on the gradient just proved, we have that 
$$\left\vert \nabla \left(w_\eps-\varphi_\eps\right)(y)\right\vert \le C\gamma_\eps^{-1} \left(1+O\left(r_\eps^{-1}\right)\right)\left( \frac{\mu_\eps}{\mu_\eps +\vert y\vert}+\gamma_\eps^{-2}\right) + C \gamma_\eps^{-2}+Cr_\eps^{-1}\gamma_\eps^{-3}$$
Since $w_\eps-\varphi_\eps=0$ on $\partial {\mathbb D}_0\left(r_\eps\right)$, this leads after integration to 
$$\left\vert w_\eps(y)-\varphi_\eps(y)\right\vert \le C\gamma_\eps^{-1}\left(1+O\left(r_\eps^{-1}\right)\right) \mu_\eps \ln \frac{\mu_\eps+r_\eps}{\mu_\eps+\vert y\vert} + O\left(\gamma_\eps^{-2}\right) $$
for all $y\in {\mathbb D}_0\left(r_\eps\right)$. This leads to
$$\left\Vert w_\eps-\varphi_\eps\right\Vert_{L^\infty\left({\mathbb D}_0\left(r_\eps\right)\right)} 
= O\left(\gamma_\eps^{-1}\mu_\eps \ln \left(1+\frac{r_\eps}{\mu_\eps}\right)\right) 
+ O\left(\gamma_\eps^{-1}\frac{\mu_\eps}{r_\eps} \ln \left(1+\frac{r_\eps}{\mu_\eps}\right)\right) +O\left(\gamma_\eps^{-2}\right)=o\left(\gamma_\eps^{-1}\right)$$
thanks to \eqref{eq-LBA-14}.

\smallskip This ends the proof of Step 2. \hfill $\spadesuit$

\medskip We are now in position to prove Proposition \ref{claim-LBA-1}. First, since $w_\eps(0)=0$ and $\nabla w_\eps(0)=0$, we get with Step 2 that 
\begin{equation}\label{eq-claim-LBA-1-39}
\varphi_\eps(0)=o\left(\gamma_\eps^{-1}\right)
\end{equation}
and that 
\begin{equation}\label{eq-claim-LBA-1-40}
\left\vert \nabla\varphi_\eps(0)\right\vert =o\left(\gamma_\eps^{-1}r_\eps^{-1}\right)+O\left(\gamma_\eps^{-1}\right)\hskip.1cm.
\end{equation}
Since $\varphi_\eps$ is harmonic, \eqref{eq-claim-LBA-1-39} gives that 
$$\gamma_\eps \varphi_\eps(0)=\frac{\gamma_\eps}{2\pi r_\eps} \int_{\partial {\mathbb D}_0\left(r_\eps\right)}\varphi_\eps\, d\sigma \to 0\hbox{ as }\eps\to 0\hskip.1cm.$$
Since $v_\eps-B_\eps=w_\eps=\varphi_\eps$ on $\partial {\mathbb D}_0\left(r_\eps\right)$, this leads to 
$$\gamma_\eps \left\vert \overline{v}_\eps\left(r_\eps\right) -B_\eps\left(r_\eps\right)\right\vert \to 0\hbox{ as }\eps\to 0\hskip.1cm,$$
which is impossible if $r_\eps<\rho_\eps$ thanks to \eqref{eq-claim-LBA-1-2}. Thus we have proved that 
\begin{equation}\label{eq-claim-LBA-1-41}
r_\eps=\rho_\eps\hskip.1cm,
\end{equation}
for any choice of $\eta\in \left(0,1\right)$. This proves the first part of a). The second part of a) is then just a consequence of \eqref{eq-LBA-12}. Indeed, $\gamma_\eps \left\vert \overline{v}_\eps\left(\rho_\eps\right) -B_\eps\left(\rho_\eps\right)\right\vert=o(1)$ implies that $\gamma_\eps B_\eps\left(\rho_\eps\right)\ge \gamma_\eps \overline{v}_\eps\left(\rho_\eps\right)+o(1)$. And \eqref{eq-LBA-12} gives that 
$$\gamma_\eps^2 - t_\eps\left(\rho_\eps\right)-\gamma_\eps^{-2}t_\eps\left(\rho_\eps\right) \ge \gamma_\eps \overline{v}_\eps\left(\rho_\eps\right)+o(1)\hskip.1cm.$$
which leads to $t_\eps\left(\rho_\eps\right)\le \gamma_\eps^2-1+o\left(1\right)$ since $ \overline{v}_\eps\left(\rho_\eps\right)\ge 0$. 
Point b) of the proposition is a consequence of Step 2 together with \eqref{eq-claim-LBA-1-8}. It remains to prove c). Let us write that 
$$\gamma_\eps \left(v_\eps\left(\rho_\eps x\right)-B_\eps\left(\rho_\eps\right)\right) = \gamma_\eps w_\eps\left(\rho_\eps x\right) + \gamma_\eps\left(B_\eps\left(\rho_\eps x\right)-B_\eps\left(\rho_\eps\right)\right)\hskip.1cm.$$
We write that 
$$\gamma_\eps \left(B_\eps\left(\rho_\eps x\right)-B_\eps\left(\rho_\eps\right)\right) \to 2\ln \frac{1}{\vert x\vert}\hbox{ in }{\mathcal C}^1_{loc}\left({\mathbb D}_0(1)\setminus \left\{0\right\}\right)\hbox{ as }\eps\to 0$$
thanks to \eqref{eq-LBA-12} and \eqref{eq-LBA-10}. Moreover, thanks to Step 2, we know that 
$$ \gamma_\eps \left\Vert w_\eps\left(\rho_\eps x\right) -\varphi_\eps\left(\rho_\eps x\right)\right\Vert_{L^\infty\left({\mathbb D}_0(1)\right)} = o(1)$$
and, combining Steps 1 and 2, that 
\begincal
&&\gamma_\eps\rho_\eps \left\vert  \nabla w_\eps\left(\rho_\eps x\right)-\nabla \varphi_\eps\left(\rho_\eps x\right)\right\vert \\
&&\quad \le C\left(\frac{\mu_\eps}{\mu_\eps +\rho_\eps \vert x\vert} +\gamma_\eps^{-2}\right) + C\rho_\eps \gamma_\eps^{-1} + C\gamma_\eps^{-2}
\fincal
in ${\mathbb D}_0(1)$. Thus we have that 
$$\gamma_\eps w_\eps\left(\rho_\eps x\right) \to \varphi_0 \hbox{ in }{\mathcal C}^1_{loc}\left({\mathbb D}_0(1)\setminus\left\{0\right\}\right)\hbox{ as }\eps\to 0$$
thanks to \eqref{eq-claim-LBA-1-8bis}. We thus have obtained that 
$$\gamma_\eps \left(v_\eps\left(\rho_\eps\, \cdot\, \right)- B_\eps\left(\rho_\eps\right)\right)\to 2 \ln \frac{1}{\vert x\vert} +\varphi_0$$
in ${\mathcal C}^1_{loc}\left({\mathbb D}_0(1)\setminus \left\{0\right\}\right)$ as $\eps\to 0$. Moreover, we have thanks to Step 2 that 
$$\varphi_0(0)=0\hbox{ and } \nabla \varphi_0(0) = -\frac{1}{2}\lim_{\eps\to 0} \frac{\rho_\eps\nabla f_\eps(0)}{f_\eps(0)}\hskip.1cm.$$
This ends the proof of the proposition. \hfill $\diamondsuit$

\section{Proof of Theorem \ref{mainthm}}\label{section-proofmainthm}

Let $\left(u_\eps\right)$ be a sequence of smooth positive solutions of 
\begin{equation}\label{eq-MT-1}
\Delta u_\eps = \lambda_\eps f_\eps u_\eps e^{u_\eps^2} \hbox{ in }\Omega,\, u_\eps=0\hbox{ on }\partial \Omega
\end{equation}
for some sequence $\left(\lambda_\eps\right)$ of positive real numbers and some sequence $\left(f_\eps\right)$ of functions in ${\mathcal C}^1\left(\overline{\Omega}\right)$ which satisfies \eqref{eq-cv-heps}. We assume that there exists $C>0$ such that 
\begin{equation}\label{eq-MT-2}
\int_\Om \left\vert \nabla u_\eps\right\vert^2\, dx \le C\hskip.1cm.
\end{equation}
We consider the concentration points $\left(\xie\right)_{i=1,\dots,N}$ given by Proposition \ref{prop-DMJ} together with the $\gamma_{i,\eps}$'s and $\mie$'s. For any $i\in \left\{1,\dots,N\right\}$, we let 
\begin{equation}\label{eq-MT-3}
r_{i,\eps} = \frac{1}{2}\min\left\{\min_{j\in\left\{1,\dots,N\right\}, j\neq i} \left\vert \xie-\xje\right\vert, \, d\left(\xie,\partial \Omega\right)\right\}\hskip.1cm.
\end{equation}
Note that we have
\begin{equation}\label{eq-MT-4}
\lambda_\eps \left\vert x-\xie\right\vert^2 u_\eps(x)^2 e^{u_\eps(x)^2} \le C_1\hbox{ in }{\mathbb D}_{\xie}\left(\rie\right)
\end{equation}
and 
\begin{equation}\label{eq-MT-5}
\left\vert x-\xie\right\vert u_\eps(x) \left\vert \nabla u_\eps(x)\right\vert \le C_2\hbox{ in }{\mathbb D}_{\xie}\left(\rie\right)
\end{equation}
thanks to assertions e) and f) of Proposition \ref{prop-DMJ}. 

We let, for $i\in\left\{1,\dots,N\right\}$, $B_{i,\eps}$ be the radial solution, studied in Appendix A, of  
$$\Delta B_{i,\eps} = \lambda_\eps f_\eps\left(\xie\right) B_{i,\eps} e^{B_{i,\eps}^2}\hbox{ and }B_{i,\eps}(0)=\gamma_{i,\eps}$$
and we shall write, by an obvious and not misleading abuse of notation, 
\begin{equation}\label{eq-MT-6}
B_{i,\eps}(x)= B_{i,\eps}\left(\left\vert x-\xie\right\vert\right)\hskip.1cm.
\end{equation}
We let also 
\begin{equation}\label{eq-MT-6bis}
t_{i,\eps}\left(r\right) =\ln\left(1+\frac{r^2}{4\mie^2}\right)\hbox{ and }t_{i,\eps}(x)=t_{i,\eps}\left(\left\vert x-\xie\right\vert\right)\hskip.1cm.
\end{equation}
At last, we define for $i=1,\dots,N$
\begin{equation}\label{eq-MT-8}
d_{i,\eps}=d\left(\xie,\partial\Omega\right)\hskip.1cm.
\end{equation}
Let us first state a claim which explains how we shall use the results of Section \ref{section-LBA} for the multibumps analysis~:

\begin{claim}\label{claim-MT-LBA}
Assume that $\left(u_\eps\right)$ satisfies equation \eqref{eq-MT-1} with $\left(f_\eps\right)$ satisfying \eqref{eq-cv-heps}. Assume also that \eqref{eq-MT-2} holds so that we have concentration points $\left(\xie\right)$ satisfying \eqref{eq-MT-4} and \eqref{eq-MT-5}. Let $0\le r_\eps\le \rie$ be such that there exists $C_3>0$ such that 
$$\left\vert x-\xie\right\vert \left\vert \nabla u_\eps(x)\right\vert \le C_3\gie^{-1}\hbox{ in }{\mathbb D}_{\xie}\left(r_\eps\right)\hskip.1cm.$$
Then we have that~:

\smallskip a) $t_{i,\eps}\left(r_\eps\right) \le \gie^2 - 1 +o\left(1\right)$ and 
$$\frac{1}{2\pi r_\eps}\int_{\partial {\mathbb D}_{\xie}\left(r_\eps\right)} u_\eps\, d\sigma = B_{i,\eps}\left(r_\eps\right)+o\left(\gie^{-1}\right)\hskip.1cm.$$

\smallskip b) There exists $C>0$ such that 
$$\left\vert u_\eps-B_{i,\eps}\right\vert \le C\gie^{-1}$$
and 
$$\left\vert \nabla \left(u_\eps-B_{i,\eps}\right)\right\vert \le C\gie^{-1}r_\eps^{-1}$$
in ${\mathbb D}_{\xie}\left(r_\eps\right)$.

\smallskip c) If $r_\eps=\rie$, after passing to a subsequence, 
$$\gie \left(u_\eps\left(\xie+\rie \, \cdot\, \right)- B_{i,\eps}\left(\rie\right)\right)\to 2 \ln \frac{1}{\vert x\vert} + {\mathcal H}_i$$
as $\eps\to 0$ in ${\mathcal C}^1_{loc}\left({\mathbb D}_0(1)\setminus \left\{0\right\}\right)$ where ${\mathcal H}_i$ is some harmonic function in the unit disk satisfying 
$${\mathcal H}_i(0)=0\hbox{ and } \nabla {\mathcal H}_i(0)=-\frac{1}{2}\lim_{\eps\to 0} \frac{r_{i,\eps}\nabla f_\eps\left(\xie\right)}{f_\eps\left(\xie\right)}\hskip.1cm.$$
\end{claim}

Let us start with a simple consequence of this claim~:

\begin{claim}\label{claim-MT-0}
For any $i\in \left\{1,\dots,N\right\}$ and any sequence $\left(r_\eps\right)$ of positive real numbers such that ${\mathbb D}_{\xie}\left(r_\eps\right)\subset \Omega$,  we have that~:

\smallskip a) If $r_\eps\le r_{i,\eps}$ and $B_{i,\eps}\left(r_\eps\right)\ge \delta \gamma_{i,\eps}$ for some $\delta >0$, there exists $C>0$ such that 
$$\left\vert u_\eps-B_{i,\eps}\right\vert \le \frac{C}{\gamma_{i,\eps}}\hbox{ in }{\mathbb D}_{\xie}\left(r_\eps\right)\hskip.1cm.$$
Moreover, we have that 
$$\frac{1}{2\pi r_\eps}\int_{\partial {\mathbb D}_{\xie}\left(r_\eps\right)} u_\eps\, d\sigma = B_{i,\eps}\left(r_\eps\right)+o\left(\gamma_{i,\eps}^{-1}\right)\hskip.1cm.$$

\smallskip b) If ${\ds \limsup_{\eps\to 0}\gie^{-1} B_{i,\eps}\left(r_\eps\right) \le 0}$ and ${\ds \limsup_{\eps\to 0}\gie^{-1} B_{i,\eps}\left(r_{i,\eps}\right) \le 0}$, then we have that 
$$\inf_{\partial{\mathbb D}_{\xie}\left(r_\eps\right)} u_\eps \le B_{i,\eps}\left(r_\eps\right) + o\left(\gamma_{i,\eps}^{-1}\right)\hskip.1cm.$$

\smallskip c) If ${\ds \limsup_{\eps\to 0}\gie^{-1} B_{i,\eps}\left(r_{i,\eps}\right) \le 0}$, we have that $t_{i,\eps}\left(\die\right)\le \gie^2$ for $\eps>0$ small enough. In other words, we have that 
$$\lambda_\eps f_\eps\left(\xie\right)\gie^2 d_{i,\eps}^2 \le 4$$
for $\eps$ small enough. Here, $\die$ is as in \eqref{eq-MT-8}.
\end{claim}

\medskip {\it Proof of Claim \ref{claim-MT-0}} - We first prove a). We assume that $B_{i,\eps}\left(r_\eps\right)\ge \delta\gamma_{i,\eps}$ for some $\delta>0$. Define 
$0\le s_\eps\le r_{\eps}$ as
$$s_\eps=\max\left\{0\le s\le r_\eps\hbox{ s.t. }u_\eps\ge \frac{1}{2}\delta \gamma_{i,\eps} \hbox{ in }{\mathbb D}_{\xie}(s)\right\}\hskip.1cm.$$
Thanks to \eqref{eq-MT-5}, we have that 
$$\left\vert x-\xie\right\vert \left\vert \nabla u_\eps\right\vert \le C\gie^{-1}\hbox{ in }{\mathbb D}_{\xie}\left(s_\eps\right)$$
for some $C>0$ so that we can apply Claim \ref{claim-MT-LBA}. Assertion b) of this claim gives that 
$$\left\vert u_\eps - B_{i,\eps}\right\vert \le C\gamma_{i,\eps}^{-1}\hbox{ in }{\mathbb D}_{\xie}\left(s_\eps\right)$$
for some $C>0$. Since $B_{i,\eps}\left(s_\eps\right)\ge B_{i,\eps}\left(r_\eps\right)\ge \delta \gamma_{i,\eps}$, we obtain in particular that $s_\eps=r_\eps$. Indeed, if $s_\eps<r_\eps$, there would exist some $x_\eps\in \partial{\mathbb D}_{\xie}\left(s_\eps\right)$ such that $u_\eps\left(x_\eps\right)=\frac{\delta}{2}\gamma_{i,\eps}$, which is impossible by what we just proved. Thus a) is clearly proved, applying again Claim \ref{claim-MT-LBA} this time with $r_\eps$. 

Let us now prove b). Let us assume first that $1+\frac{r_\eps^2}{4\mie^2}\le e^{\gie^2}$, that ${\ds \limsup_{\eps\to 0}\gie^{-1} B_{i,\eps}\left(r_\eps\right) \le 0}$ and that ${\ds \limsup_{\eps\to 0}\gie^{-1} B_{i,\eps}\left(r_{i,\eps}\right) \le 0}$ and assume by contradiction that there exists $0<\eta<1$ such that 
\begin{equation}\label{eq-claim-MT-0-1}
\inf_{{\mathbb D}_{\xie}\left(r_\eps\right)} u_\eps \ge B_{i,\eps}\left(r_\eps\right) + \eta\gamma_{i,\eps}^{-1}\hskip.1cm.
\end{equation}
We claim that 
\begin{equation}\label{eq-claim-MT-0-2}
u_\eps\ge \gie +\frac{1}{\gamma_{i,\eps}} \ln \frac{4\mie^2}{\left\vert \xie-x\right\vert^2} -\frac{1-\eta}{\gie} +o\left(\gie^{-1}\right)\hbox{ in }{\mathbb D}_{\xie}\left(r_\eps\right)\setminus {\mathbb D}_{\xie}\left(R_0\mie\right)
\end{equation}
where 
$$R_0=\frac{4}{\sqrt{e^{1-\eta}-1}}\hskip.1cm.$$
The right-hand side of \eqref{eq-claim-MT-0-2} being harmonic and $u_\eps$ being super-harmonic, it is sufficient to check the inequality on $\partial {\mathbb D}_{\xie}\left(r_\eps\right)$ and on $\partial {\mathbb D}_{\xie}\left(R_0\mie\right)$. For that purpose, let us write that 
$$B_{i,\eps}\left(r_\eps\right) = \gie - \gie^{-1} t_{i,\eps}\left(r_\eps\right)-\gie^{-3} t_{i,\eps}\left(r_\eps\right)+O\left(\gie^{-2}\right)$$
 as proved in Appendix A, Claim \ref{claim-appA-4}, since we assumed for the moment that $t_{i,\eps}\left(r_\eps\right)\le \gamma_{i,\eps}^2$. Since we assumed that ${\ds \limsup_{\eps\to 0}\gie^{-1} B_{i,\eps}\left(r_\eps\right) \le 0}$, this gives that ${\ds \frac{t_{i,\eps}\left(r_\eps\right)}{\gie^2}\to 1}$ as $\eps\to 0$ so that 
\begincal
B_{i,\eps}\left(r_\eps\right) &=& \gie-\gie^{-1}\ln \left(1+\frac{r_\eps^2}{4\mie^2}\right) - \gie^{-1} +o\left(\gie^{-1}\right)\\
&=& \gie +\frac{1}{\gamma_{i,\eps}} \ln \frac{4\mie^2}{r_\eps^2} -\gie^{-1} +o\left(\gie^{-1}\right)\hskip.1cm.
\fincal
This implies with \eqref{eq-claim-MT-0-1} that 
\begin{equation}\label{eq-claim-MT-0-3}
u_\eps\ge \gie +\frac{1}{\gamma_{i,\eps}} \ln \frac{4\mie^2}{\left\vert \xie-x\right\vert^2} -\frac{1-\eta}{\gie} +o\left(\gie^{-1}\right)\hbox{ on }\partial {\mathbb D}_{\xie}\left(r_\eps\right)\hskip.1cm.
\end{equation}
Let us write now that 
$$u_\eps-B_{i,\eps} = o\left(\gie^{-1}\right) \hbox{ on }\partial {\mathbb D}_{\xie}\left(R_0\mie\right)$$
thanks to d) of Proposition \ref{prop-DMJ}. Since 
$$B_{i,\eps}\left(R_0\mie\right)= \gie - \gie^{-1}\ln \left(1+\frac{R_0^2}{4}\right) +o\left(\gie^{-1}\right)\hskip.1cm,$$
we obtain that 
\begin{equation}\label{eq-claim-MT-0-4}
u_\eps\ge \gie +\frac{1}{\gamma_{i,\eps}} \ln \frac{4\mie^2}{\left\vert \xie-x\right\vert^2} -\frac{1-\eta}{\gie} +o\left(\gie^{-1}\right)\hbox{ on }\partial {\mathbb D}_{\xie}\left(R_0\mie\right)
\end{equation}
provided that 
$$ \ln\left(1+\frac{4}{R_0^2}\right)< 1-\eta \hskip.1cm,$$
which is the case with our choice of $R_0$. Thus \eqref{eq-claim-MT-0-2} is proved.

Now there exists $R_0\mie\le s_\eps\le \min\left\{r_\eps, r_{i,\eps}\right\}$ such that $B_{i,\eps}\left(s_\eps\right)=\frac{\eta}{2} \gamma_{i,\eps}$ since ${\ds \limsup_{\eps\to 0}\gie^{-1} B_{i,\eps}\left(r_\eps\right) \le 0}$ and ${\ds \limsup_{\eps\to 0}\gie^{-1} B_{i,\eps}\left(r_{i,\eps}\right) \le 0}$. We can apply a) of the claim to get that 
$$\frac{1}{2\pi s_\eps} \int_{\partial {\mathbb D}_{\xie}\left(s_\eps\right)} u_\eps \, d\sigma =\frac{\eta}{2}\gie + o\left(\gie^{-1}\right)\hskip.1cm.$$
Applying \eqref{eq-claim-MT-0-2}, this leads to 
\begin{equation}\label{eq-claim-MT-0-5}
\gie +\frac{1}{\gamma_{i,\eps}} \ln \frac{4\mie^2}{s_\eps^2} -\frac{1-\eta}{\gie} +o\left(\gie^{-1}\right)\le \frac{\eta}{2}\gie + o\left(\gie^{-1}\right)\hskip.1cm.
\end{equation}
Since $B_{i,\eps}\left(s_\eps\right)=\frac{\eta}{2} \gie$, it is not difficult to check thanks to Claim \ref{claim-appA-4} of Appendix A that 
$$t_{i,\eps}\left(s_\eps\right)=\left(1-\frac{\eta}{2}\right)\left(\gie^2-1\right)+O\left(\gie^{-1}\right)$$
so that, since $\frac{s_\eps}{\mie}\to +\infty$ as $\eps\to 0$, 
$$\ln \frac{4\mie^2}{s_\eps^2} =-\left(\gie^2-1\right)\left(1-\frac{\eta}{2}\right) + o\left(1\right)\hskip.1cm.$$
Coming back to \eqref{eq-claim-MT-0-5} with this leads to a contradiction. This proves that \eqref{eq-claim-MT-0-1} is absurd for any $0<\eta<1$. Thus we have proved assertion b) as long as $t_{i,\eps}\left(r_\eps\right)\le \gie^2$. 

\smallskip We shall now prove c), which will by the way prove that b) holds since the condition $t_{i,\eps}\left(r_\eps\right)\le \gie^2$ will always be satisfied.  Let us assume by contradiction that $t_{i,\eps}\left(\die\right)\ge \gie^2$. Then ${\mathbb D}_{\xie}\left(r_\eps\right)\subset \Omega$ for $\eps>0$ small where
$$1+\frac{r_\eps^2}{4\mu_\eps^2} = e^{\gamma_\eps^2-\frac{1}{2}}\hskip.1cm.$$
We can apply b) in this case since $t_{i,\eps}\left(r_\eps\right)\le \gie^2$ and 
$$B_{i,\eps}\left(r_\eps\right)=  -\frac{1}{2}\gie^{-1}+O\left(\gie^{-2}\right)$$
by Claim \ref{claim-appA-4} of Appendix A. This leads to a contradiction since $u_\eps\ge 0$ in $\Omega$. Thus c) is proved thanks to the definition of $\mie$ and b) is also proved. This ends the proof of this claim. \hfill $\diamondsuit$

\begin{claim}\label{claim-MT-1}
For any $i\in \left\{1,\dots,N\right\}$, we have that 
$$\limsup_{\eps\to 0} \gamma_{i,\eps}^{-1}B_{i,\eps}\left(\rie\right) \le 0\hskip.1cm.$$
\end{claim}

\medskip {\it Proof of Claim \ref{claim-MT-1}} - Let us reorder for this proof the concentration points in such a way that 
\begin{equation}\label{eq-claim-MT-1-1}
r_{1,\eps}\le r_{2,\eps}\le \dots \le r_{N,\eps}\hskip.1cm.
\end{equation}
We prove the assertion by induction on $i$. Let $i\in \left\{1,\dots,N\right\}$ and let us assume that 
\begin{equation}\label{eq-claim-MT-1-2}
\limsup_{\eps\to 0} \gamma_{j,\eps}^{-1}B_{j,\eps}\left(\rje\right) \le 0\hbox{ for all }1\le j\le i-1\hskip.1cm.
\end{equation}
Note that we do not assume anything if $i=1$. We proceed by contradiction, assuming that, after passing to a subsequence,
\begin{equation}\label{eq-claim-MT-1-3}
\gamma_{i,\eps}^{-1}B_{i,\eps}\left(\rie\right) \ge 2\eps_0
\end{equation}
for some $\eps_0>0$. 

\medskip {\sc Step 1} - {\it If \eqref{eq-claim-MT-1-3} holds, then $\frac{d\left(\xie,\partial \Omega\right)}{\rie}\to +\infty$ as $\eps\to 0$. In particular, this implies that $r_{i,\eps}\to 0$ as $\eps\to 0$.}

\medskip {\sc Proof of Step 1} - For any $\eta>0$ small enough, there exists a path of length less than or equal to $Cd\left(\xie,\partial\Omega\right)$ joining the boundary of $\Omega$ and the boundary of the disk ${\mathbb D}_{\xie}\left(\eta d\left(\xie,\partial\Omega\right)\right)$, and avoiding all the disks ${\mathbb D}_{\xje}\left(\eta d\left(\xie,\partial\Omega\right)\right)$ for $j=1,\dots, N$. Using f) of Proposition \ref{prop-DMJ}, we deduce that, for any $\eta>0$, there exists $C>0$ such that 
$$u_\eps \le C \hbox{ on }\partial {\mathbb D}_{\xie}\left(\eta d\left(\xie,\partial\Omega\right)\right)\hskip.1cm.$$
If $d\left(\xie,\partial\Omega\right)=O\left(\rie\right)$, we can find $\eta>0$ small enough such that $\eta d\left(\xie,\partial\Omega\right)\le \rie$. Then the above estimate would clearly contradict a) of Claim \ref{claim-MT-0} together with \eqref{eq-claim-MT-1-3}. Thus Step 1 is proved. \hfill $\spadesuit$

\medskip Thanks to Step 1, we know that, if  \eqref{eq-claim-MT-1-3} holds, then 
\begin{equation}\label{eq-claim-MT-1-4}
{\mathcal D}_i = \left\{j\in\left\{1,\dots,N\right\},\, j\neq i \hbox{ s.t. } \left\vert \xje-\xie\right\vert =O\left(\rie\right)\right\} \neq \emptyset\hskip.1cm.
\end{equation}
There exists $0<\delta<1$ such that, for any $j\in {\mathcal D}_i$, any point of $\partial {\mathbb D}_{\xje}\left(\delta\rie\right)$ can be joined to a point of $\partial {\mathbb D}_{\xie}\left(\delta\rie\right)$ by a path $\gamma_\eps:[0,1]\to \Omega$ such that $\left\vert \gamma_\eps(t)-\xke\right\vert\ge \delta\rie$ for all $k=1,\dots,N$ and all $0\le t\le 1$ and such that $\left\vert \gamma_\eps'(t)\right\vert \le \delta^{-1}\rie$. Thanks to assertion f) of Proposition \ref{prop-DMJ}, the existence of such paths give that 
$$\inf_{\partial {\mathbb D}_{\xje}\left(\delta\rie\right)} u_\eps^2 \ge \inf_{\partial {\mathbb D}_{\xie}\left(\delta\rie\right)} u_\eps^2 - 2C_2\delta^{-2}\hbox{ for all }j\in {\mathcal D}_i\hskip.1cm.$$
Thanks to \eqref{eq-claim-MT-1-3}, we can apply a) of Claim \ref{claim-MT-0} to obtain also that 
$$u_\eps \ge B_{i,\eps}\left(\delta\rie\right)-C\gie^{-1}\hbox{ on }\partial {\mathbb D}_{\xie}\left(\delta\rie\right)$$
for some $C>0$. Since 
$$B_{i,\eps}\left(\delta\rie\right)=B_{i,\eps}\left(\rie\right)+O\left(\gie^{-1}\right)\hskip.1cm,$$
the two previous estimates , together with \eqref{eq-claim-MT-1-3}, lead to the existence of some $C>0$ such that 
\begin{equation}\label{eq-claim-MT-1-5}
u_\eps\ge B_{i,\eps}\left(r_{i,\eps}\right) -C\gie^{-1} \hbox{ on }\partial {\mathbb D}_{\xje}\left(\delta\rie\right)\hbox{ for all }j\in{\mathcal D}_i\hskip.1cm.
\end{equation}

\medskip {\sc Step 2} - {\it If \eqref{eq-claim-MT-1-3} holds, then for any $j\in {\mathcal D}_i$, we have that 
$$\liminf_{\eps\to 0} \gje^{-1} B_{j,\eps}\left(\rje\right) >0\hskip.1cm.$$
In particular, we have that $j\ge i+1$.}

\medskip {\sc Proof of Step 2} - Assume on the contrary that there exists $j\in {\mathcal D}_i$ such that, after passing to a subsequence, 
\begin{equation}\label{eq-claim-MT-1-6}
\limsup_{\eps\to 0} \gje^{-1} B_{j,\eps}\left(r_{j,\eps}\right)\le 0\hskip.1cm.
\end{equation}
Since $j\in {\mathcal D}_i$, we also know that 
\begin{equation}\label{eq-claim-MT-1-7}
r_{j,\eps}\le \frac{1}{2}\left\vert \xie-\xje\right\vert\le C\rie\hskip.1cm.
\end{equation}
Thus we also have that 
\begin{equation}\label{eq-claim-MT-1-8}
\limsup_{\eps\to 0} \gje^{-1} B_{j,\eps}\left(\delta r_{i,\eps}\right)\le 0\hskip.1cm.
\end{equation}
We can apply b) of Claim \ref{claim-MT-0} with $r_\eps=\delta\rie$ to obtain that 
\begin{equation}\label{eq-claim-MT-1-9}
B_{i,\eps}\left(r_{i,\eps}\right) -C\gie^{-1} \le B_{j,\eps}\left(\delta \rie\right) + o\left(\gje^{-1}\right)
\end{equation}
thanks to \eqref{eq-claim-MT-1-5}. Combining \eqref{eq-claim-MT-1-3} and \eqref{eq-claim-MT-1-8}, we get that 
\begin{equation}\label{eq-claim-MT-1-10}
\gie=o\left(\gje\right)\hskip.1cm.
\end{equation}
Thus we also have that $\mu_{j,\eps}\le\mie$. Let us write now thanks to Claim \ref{claim-appA-4} of Appendix A that 
$$B_{j,\eps}\left(\delta \rie\right) = -\gje^{-1} \ln \left(\frac{\rie^2}{4}\right) - \gje^{-1}\ln \left(\lambda_\eps \gje^2\right) + O\left(\gje^{-1}\right)$$
and that 
$$B_{i,\eps}\left(\rie\right) = -\gie^{-1}\ln \left(\frac{\rie^2}{4}\right) - \gie^{-1}\ln \left(\lambda_\eps \gie^2\right) + O\left(\gie^{-1}\right)$$
to obtain that 
$$B_{j,\eps}\left(\delta \rie\right) = \frac{\gamma_{i,\eps}}{\gamma_{j,\eps}} B_{i,\eps}\left(\rie\right)+ \gje^{-1} \ln \left(\frac{\gie^2}{\gje^2}\right) + O\left(\gje^{-1}\right)\hskip.1cm.$$
Coming back to \eqref{eq-claim-MT-1-9} with this, \eqref{eq-claim-MT-1-3} and \eqref{eq-claim-MT-1-10}, we obtain that 
$$\left(2\eps_0+o(1)\right) \gie\le \gje^{-1} \ln \left(\frac{\gie^2}{\gje^2}\right) +O\left(\gamma_{i,\eps}^{-1}\right)\le O\left(\gie^{-1}\right)\hskip.1cm,$$
which is a clear contradiction. Step 2 is proved. \hfill $\spadesuit$

\medskip We can now conclude the proof of the claim by proving that \eqref{eq-claim-MT-1-3} is absurd if \eqref{eq-claim-MT-1-2} holds. Continue to assume that \eqref{eq-claim-MT-1-3} holds. Then we know thanks to Step 2 that for any $j\in {\mathcal D}_i$, $j\ge i+1$ so that $\rje\ge \rie$. We set, for $j\in {\mathcal D}_i$, and up to a subsequence, 
\begin{equation}\label{eq-claim-MT-1-11}
\hat{x}_j = \lim_{\eps\to 0} \frac{\xje-\xie}{\rie}
\end{equation}
and we let 
\begin{equation}\label{eq-claim-MT-1-12}
\hat{\mathcal S} = \left\{\hat{x}_j,\, j\in {\mathcal D}_i\right\}\hskip.1cm.
\end{equation}
We know thanks to Step 1 that there exists $j\in {\mathcal D}_i$ such that 
\begin{equation}\label{eq-claim_MT-1-13}
\left\vert \hat{x}_j\right\vert = 2
\end{equation}
and that 
\begin{equation}\label{eq-claim-MT-1-14}
\left\vert \hat{x}_k-\hat{x}_l\right\vert \ge 2 \hbox{ for all }k,l\in {\mathcal D}_i,\, k\neq l\hskip.1cm.
\end{equation}
Since $\rje$ and $\rie$ are comparable, we also have thanks to Step 2 that 
\begin{equation}\label{eq-claim-MT-1-15}
\liminf_{\eps\to 0} \gje^{-1} B_{j,\eps}\left(\rje\right) >0\hskip.1cm.
\end{equation}
Let $K$ be a compact subset of ${\mathbb R}^2\setminus \hat {\mathcal S}$. We can use assertion f) of Proposition \ref{prop-DMJ} to write\footnote{see the argument between Steps 1 and 2.} that 
\begin{equation}\label{eq-claim-MT-1-16}
\gie \left\vert u_\eps\left(\xie+\rie x\right) -B_{i,\eps}\left(\rie\right)\right\vert \le C_K\hbox{ in }K\hskip.1cm.
\end{equation}
Thanks to \eqref{eq-MT-1}, we can write that 
$$\Delta \hat{u}_\eps = \lambda_\eps \rie^2 \gie f_\eps\left(\xie+\rie x\right) u_\eps\left(\xie+\rie x\right)^2 e^{u_\eps\left(\xie+\rie x\right)^2}$$
where 
$$\hat{u}_\eps = \gie \left(u_\eps\left(\xie+\rie x\right) -B_{i,\eps}\left(\rie\right)\right)\hskip.1cm.$$
Using \eqref{eq-claim-MT-1-16}, we can write that 
$$\left\vert \Delta \hat{u}_\eps\right\vert \le C_K \mie^{-2} \rie^2  e^{B_{i,\eps}\left(\rie\right)^2-\gie^2}\hbox{ in }K$$
for any compact subset $K$ of ${\mathbb R}^2\setminus \hat {\mathcal S}$. Thanks to \eqref{eq-claim-MT-1-3}, we have that 
$$e^{B_{i,\eps}\left(\rie\right)^2-\gie^2}\le C\left(1+\frac{\rie^2}{4\mie^2}\right)^{-1-2\eps_0}$$
so that 
$$\left\vert \Delta \hat{u}_\eps\right\vert \le C_K \left(\frac{\mie}{\rie}\right)^{4\eps_0}\to 0 \hbox{ uniformly} \hbox{ in }K\hskip.1cm.$$
By standard elliptic theory, we thus have that 
\begin{equation}\label{eq-claim-MT-1-17}
\hat{u}_\eps = \gie \left(u_\eps\left(\xie+\rie x\right) -B_{i,\eps}\left(\rie\right)\right) \to \hat{u}_0\hbox{ in }{\mathcal C}^1_{loc}\left({\mathbb R}^2\setminus \hat{\mathcal S}\right)\hbox{ as }\eps\to 0
\end{equation}
where 
\begin{equation}\label{eq-claim-MT-1-18}
\Delta \hat{u}_0=0\hbox{ in }{\mathbb R}^2\setminus \hat{\mathcal S}\hskip.1cm.
\end{equation}
Since $\rje\ge \rie$ for $j\in {\mathcal D}_i$, \eqref{eq-claim-MT-1-15} permits to apply a) of Claim \ref{claim-MT-0}, which in turn implies thanks to \eqref{eq-MT-5} that we can apply Claim \ref{claim-MT-LBA} for all $j\in{\mathcal D}_i$ with $r_\eps=\rie$. Assertion c) of this claim gives that 
\begin{equation}\label{eq-claim-MT-1-19}
\gje\left(u_\eps\left(\xje+\rie x\right)-B_{j,\eps}\left(\rie\right)\right) \to 2\ln \frac{1}{\vert x\vert} + {\mathcal H}_j
\end{equation}
in ${\mathcal C}^1_{loc}\left({\mathbb D}_0(1)\setminus\left\{0\right\}\right)$ as $\eps\to 0$ 
where ${\mathcal H}_j$ is harmonic in the unit disk and satifies ${\mathcal H}_j(0)=0$ and $\nabla {\mathcal H}_j(0)=0$ (note here that we know thanks to Step 1 that $\rie\to 0$ as $\eps\to 0$). This gives that 
\begin{equation}\label{eq-claim-MT-1-20}
\frac{\gje}{\gie} \hat{u}_\eps +\gje\left(B_{i,\eps}\left(\rie\right)-B_{j,\eps}\left(\rie\right) \right)\to 2\ln \frac{1}{\vert x-\hat{x}_j \vert} + \hat{\mathcal H}_j
\end{equation}
in ${\mathcal C}^1_{loc}\left({\mathbb D}_{\hat{x}_j}(1)\setminus \left\{\hat{x}_j\right\}\right)$ as $\eps\to 0$ for all $j\in {\mathcal D}_i$ (and also for $j=i$ if we set $\hat{x}_i=0$). It remains to write thanks to Claim \ref{claim-appA-4} of Appendix A and to \eqref{eq-cv-heps} that
$$B_{i,\eps}\left(\rie\right)-B_{j,\eps}\left(\rie\right) = \left(1-\frac{\gie}{\gje}\right)B_{i,\eps}\left(\rie\right)  +\gje^{-1}\ln \frac{\gje^2}{\gie^2}+O\left(\gje^{-1}\right)$$
to deduce from \eqref{eq-claim-MT-1-3}, \eqref{eq-claim-MT-1-17} and \eqref{eq-claim-MT-1-20} that 
\begin{equation}\label{eq-claim-MT-1-21}
\frac{\gje}{\gie}\to 1\hbox{ as }\eps\to 0 \hbox{ and }\gie\left\vert \gie-\gje\right\vert = O\left(1\right)\hskip.1cm.
\end{equation}
Then  \eqref{eq-claim-MT-1-17} and \eqref{eq-claim-MT-1-20} just lead to 
\begin{equation}\label{eq-claim-MT-1-22}
\hat{u}_0 = 2\ln \frac{1}{\left\vert x-\hat{x}_j\right\vert} + \varphi_j
\end{equation}
in ${\mathbb D}_{\hat{x}_j}\left(1\right)\setminus \left\{\hat{x}_j\right\}$ where $\varphi_j$ is smooth and harmonic and satisfies $\nabla \varphi_j\left(\hat{x}_j\right)=0$. Thus we can write that 
\begin{equation}\label{eq-claim-MT-1-23}
\hat{u}_0=2\ln\frac{1}{\vert x\vert} + 2\sum_{j\in {\mathcal D}_i} \ln\frac{1}{\left\vert x -\hat{x}_j\right\vert} + \varphi
\end{equation}
where $\varphi$ is a smooth harmonic function in ${\mathbb R}^2$. Thanks to assertion f) of Proposition \ref{prop-DMJ}, we also know that 
$$\left\vert \nabla \varphi(x)\right\vert \le \frac{C}{1+\vert x\vert}\hbox{ in }{\mathbb R}^2$$
for some $C>0$ so that $\varphi\equiv Cst$. Now this gives that for any $k\in {\mathcal D}_i$,
$$\nabla \left(\ln\frac{1}{\vert x\vert}+\sum_{j\in {\mathcal D}_i, j\neq k} \ln\frac{1}{\left\vert x- \hat{x}_j\right\vert}\right)\left(\hat{x}_k\right)=0\hskip.1cm.$$
Let $k\in {\mathcal D}_i$ be such that $\left\vert \hat{x}_k\right\vert \ge \left\vert \hat{x}_j\right\vert$ for all $j\in {\mathcal D}_i$. Then 
 \begincal
\langle \nabla \left(\ln\frac{1}{\vert x\vert}+\sum_{j\in {\mathcal D}_i, \, j\neq k} \ln\frac{1}{\left\vert x -\hat{x}_j\right\vert}\right)\left(\hat{x}_k\right), \hat{x}_k\rangle &=& -\left\vert \hat{x}_k\right\vert - \sum_{j\in {\mathcal D}_i,\, j\neq k} \frac{\left\vert \hat{x}_k\right\vert^2-\langle \hat{x}_k,\hat{x}_j\rangle}{\left\vert \hat{x}_k-\hat{x}_j\right\vert}\\
&\le &-\left\vert \hat{x}_k\right\vert <0\hskip.1cm,
\fincal
which gives the desired contradiction. This proves that \eqref{eq-claim-MT-1-3} is absurd as soon as \eqref{eq-claim-MT-1-2} holds. And this ends the proof of the claim by an induction argument. \hfill $\diamondsuit$

\begin{claim}\label{claim-MT-2}
For any $i=1,\dots,N$, we have that 
$$\lambda_\eps f_\eps\left(\xie\right)\gie^2 d_{i,\eps}^2 \le 4$$
for $\eps$ small enough.
\end{claim}

\medskip {\it Proof of Claim \ref{claim-MT-2}} - It is a direct consequence of c) of Claim \ref{claim-MT-0} together with Claim \ref{claim-MT-1}. \hfill $\diamondsuit$

\begin{claim}\label{claim-MT-3}
We have that 
$$\int_\Omega \left\vert \nabla u_\eps\right\vert^2\, dx = \int_{\Omega}\left\vert \nabla u_0\right\vert^2\, dx + 4\pi N +o(1)\hskip.1cm.$$
In other words, $M=0$ in Theorem \ref{thm-DMJ}.
\end{claim}

\medskip {\it Proof of Claim \ref{claim-MT-3}} - We prove that $M=0$ in Theorem \ref{thm-DMJ}. Assume on the contrary that there exists some sequence $\left(y_{1,\eps}\right)$ such that the assertion b) of Theorem \ref{thm-DMJ} holds. This means that 
$$\nu_{1,\eps}^{-2} = \lambda_\eps f_\eps\left(y_{1,\eps}\right) u_\eps\left(y_{1,\eps}\right)^2 e^{u_\eps\left(y_{1,\eps}\right)^2} \to +\infty\hbox{ as }\eps\to 0\hskip.1cm.$$
By e) of Proposition \ref{prop-DMJ}, we know that 
\begincal
\nu_{1,\eps}^{-2} \left(\min_{i=1,\dots, N} \left\vert \xie-y_{1,\eps}\right\vert\right)^2 &=& \left(\min_{i=1,\dots, N} \left\vert \xie-y_{1,\eps}\right\vert\right)^2 \lambda_\eps f_\eps\left(y_{1,\eps}\right) u_\eps\left(y_{1,\eps}\right)^2 e^{u_\eps\left(y_{1,\eps}\right)^2} \\
&\le & C_1 f_\eps\left(y_{1,\eps}\right)\hskip.1cm.
\fincal
This proves that there exists $i\in \left\{1,\dots,N\right\}$ such that 
$$\left\vert \xie-y_{1,\eps}\right\vert =O\left(\nu_{1,\eps}\right)\hskip.1cm.$$
Since 
$$\frac{\left\vert \xie-y_{1,\eps}\right\vert }{\mie}\to +\infty\hbox{ as }\eps\to 0$$
by a) of Theorem \ref{thm-DMJ}, we have that 
$$\frac{\nu_{1,\eps}}{\mie}\to +\infty\hbox{ as }\eps\to 0\hskip.1cm.$$
Thanks to the definition of $\nu_{1,\eps}$ and $\mie$, this leads to 
$$e^{\gie^2 - u_\eps\left(y_{1,\eps}\right)^2}\frac{\gie^2}{u_\eps\left(y_{1,\eps}\right)^2}\to +\infty\hbox{ as }\eps\to 0\hskip.1cm,$$
which implies that 
\begin{equation}\label{eq-claim-MT-3-1}
\gie^2-u_\eps\left(y_{1,\eps}\right)^2\to +\infty\hbox{ as }\eps\to 0\hskip.1cm.
\end{equation}
Now, by the convergence of b) of Theorem \ref{thm-DMJ}, we know that
$$u_\eps\ge u_\eps\left(y_{1,\eps}\right) - C u_\eps\left(y_{1,\eps}\right)^{-1} \hbox{ on }\partial {\mathbb D}_{\xie}\left(R\nu_{1,\eps}\right)$$
for some $R>0$ and $C>0$. Thanks to Claim \ref{claim-MT-1}, we can use assertion b) of Claim \ref{claim-MT-0} to deduce that 
$$u_\eps\left(y_{1,\eps}\right) - C u_\eps\left(y_{1,\eps}\right)^{-1}\le B_{i,\eps}\left(R\nu_{1,\eps}\right) + o\left(\gie^{-1}\right)\hskip.1cm.$$
This leads after some simple computations, using Claim \ref{claim-appA-4} of Appendix A, to 
$$u_\eps\left(y_{1,\eps}\right) - C u_\eps\left(y_{1,\eps}\right)^{-1}\le \gie^{-1} \ln \left(\frac{u_\eps\left(y_{1,\eps}\right)^2}{\gie^2}\right) + \frac{u_\eps\left(y_{1,\eps}\right)^2}{\gie} + C\gie^{-1}$$
so that, thanks to \eqref{eq-claim-MT-3-1}, 
$$u_\eps\left(y_{1,\eps}\right)^2 \left(1-\frac{u_\eps\left(y_{1,\eps}\right)}{\gie}\right) \le C\hskip.1cm.$$
This clearly implies that 
$$\frac{u_\eps\left(y_{1,\eps}\right)}{\gie} \to 1\hbox{ as }\eps\to 0$$
and then that 
$$u_\eps\left(y_{1,\eps}\right)\ge \gamma_{i,\eps} - C\gie^{-1}$$
for some $C>0$. This contradicts \eqref{eq-claim-MT-3-1}. Thus we have proved that $M=0$ in Theorem \ref{thm-DMJ} and the claim follows. \hfill $\diamondsuit$ 

\medskip For any $i\in \left\{1,\dots,N\right\}$, thanks to Claim \ref{claim-MT-1} and a) of Claim \eqref{claim-MT-0}, there exists $0\le \sie\le \rie$ such that 
\begin{equation}\label{eq-MT-9}
\limsup_{\eps\to 0} \gie^{-1} B_{i,\eps}\left(\sie\right)\le 0 \hbox{ and }\left\vert u_\eps-B_{i,\eps}\right\vert \le \frac{D_i}{\gie}\hbox{ in }{\mathbb D}_{\xie}\left(\sie\right)
\end{equation}
for some $D_i>0$. 

\begin{claim}\label{claim-MT-4}
We have that 
$$\liminf_{\eps\to 0} \int_{{\mathbb D}_{\xie}\left(\sie\right)} \left\vert \nabla u_\eps\right\vert^2\, dx \ge 4\pi\hskip.1cm.$$
\end{claim}

\medskip {\it Proof of Claim \ref{claim-MT-4}} - Let $\delta>0$. Let us write thanks to \eqref{eq-MT-9} that 
$$ \int_{{\mathbb D}_{\xie}\left(\sie\right)} \left\vert \nabla u_\eps\right\vert^2\, dx\ge \int_{{\mathbb D}_{\xie}\left(\sie\right)} \left\vert \nabla \left(u_\eps-\delta\gie\right)^+\right\vert^2\, dx=\int_{{\mathbb D}_{\xie}\left(\sie\right)} \left(u_\eps-\delta\gie\right)^+\Delta u_\eps\, dx\hskip.1cm.$$
Thanks to \eqref{eq-MT-1}, this leads to 
\begincal
\int_{{\mathbb D}_{\xie}\left(\sie\right)} \left\vert \nabla u_\eps\right\vert^2\, dx&\ge& \lambda_\eps \int_{{\mathbb D}_{\xie}\left(\sie\right)} f_\eps \left(u_\eps-\delta\gie\right)^+ u_\eps e^{u_\eps^2}\, dx\\
&\ge & \lambda_\eps \int_{{\mathbb D}_{\xie}\left(R\mie\right)} f_\eps \left(u_\eps-\delta\gie\right)^+ u_\eps e^{u_\eps^2}\, dx
\fincal
for all $R>0$. Now we have that 
$$\lim_{\eps\to 0}\lambda_\eps \int_{{\mathbb D}_{\xie}\left(R\mie\right)} f_\eps \left(u_\eps-\delta\gie\right)^+ u_\eps e^{u_\eps^2}\, dx = \left(1-\delta\right) \int_{{\mathbb D}_0(R)} e^{2U}\, dx$$
thanks to d) of Proposition \ref{prop-DMJ}. Since 
$$\int_{{\mathbb R}^2} e^{2U}\, dx=4\pi\hskip.1cm,$$
the result follows by letting $R$ go to $+\infty$ and $\delta$ go to $0$. \hfill $\diamondsuit$

\medskip Let us set now 
\begin{equation}\label{eq-MT-10}
\Omega_\eps = \Omega\setminus \bigcup_{i=1}^N {\mathbb D}_{\xie}\left(\sie\right)
\end{equation}
where $\sie$ is as in \eqref{eq-MT-9} and 
\begin{equation}\label{eq-MT-11}
w_\eps = \left\{\begin{array}{ll}
{\ds u_\eps}&{\ds \hbox{ in }\Omega_\eps}\\
{\ds \min\left\{u_\eps, B_{i,\eps}\left(\sie\right)+2\frac{D_i}{\gie}\right\}}&{\ds \hbox{ in }{\mathbb D}_{\xie}\left(\sie\right)\hbox{ for }i=1,\dots,N}
\end{array}\right.
\end{equation}

\begin{claim}\label{claim-MT-5}
We have that 
$$\int_\Omega \left\vert \nabla\left(w_\eps-u_0\right)\right\vert^2\, dx \to 0\hbox{ as }\eps\to 0\hskip.1cm.$$
\end{claim}

\medskip {\it Proof of Claim \ref{claim-MT-5}} - Let us write that 
\begincal
\int_\Omega \left\vert \nabla\left(w_\eps-u_0\right)\right\vert^2\, dx &=& \int_\Omega \left\vert \nabla w_\eps\right\vert^2\, dx - 2\int_\Omega \langle \nabla w_\eps,\nabla u_0\rangle\, dx + \int_\Omega \left\vert \nabla u_0\right\vert^2\, dx\\
&=& \int_\Omega \left\vert \nabla u_\eps\right\vert^2\, dx - 2\int_\Omega \langle \nabla u_\eps,\nabla u_0\rangle\, dx +\int_\Omega \left\vert \nabla u_0\right\vert^2\, dx\\
&&\quad + \int_\Omega \langle \nabla\left(w_\eps-u_\eps\right), \nabla u_\eps+\nabla w_\eps-2\nabla u_0 \rangle\, dx\\
&=& 4\pi N +o(1)+ \int_\Omega \langle \nabla\left(w_\eps-u_\eps\right), \nabla u_\eps+\nabla w_\eps-2\nabla u_0 \rangle\, dx
\fincal
thanks to the weak convergence of $u_\eps$ to $u_0$ in $H^1$ and to Claim \ref{claim-MT-3}. Let us remark now that $\nabla\left(w_\eps-u_\eps\right)\equiv 0$ in $\Omega_\eps$ and that $\langle \nabla\left(w_\eps-u_\eps\right),\nabla w_\eps\rangle=0$ a.e. Thus we can write that 
\begin{equation}\label{eq-claim-MT-5-1}
\int_\Omega \left\vert \nabla\left(w_\eps-u_0\right)\right\vert^2\, dx  = 4\pi N + o(1)+\sum_{i=1}^N \int_{{\mathbb D}_{\xie}\left(\sie\right)} \langle \nabla\left(w_\eps-u_\eps\right), \nabla u_\eps-2\nabla u_0 \rangle\, dx\hskip.1cm.
\end{equation}
Since $w_\eps-u_\eps$ is null on the boundary of ${\mathbb D}_{\xie}\left(\sie\right)$, we can proceed as in the proof of Claim \ref{claim-MT-4} to get that 
\begincal
\int_{{\mathbb D}_{\xie}\left(\sie\right)} \langle \nabla\left(w_\eps-u_\eps\right), \nabla u_\eps-2\nabla u_0 \rangle\, dx &=& \int_{{\mathbb D}_{\xie}\left(\sie\right)} \left(w_\eps-u_\eps\right)\left(\Delta u_\eps-2\Delta u_0\right)\, dx\\
&\le & -4\pi +o(1) + O\left(\gie \int_{{\mathbb D}_{\xie}\left(\sie\right)}  \left\vert \Delta u_0\right\vert\, dx\right)\hskip.1cm.
\fincal
Here we used the fact that $w_\eps\le u_\eps$ and $\left\vert w_\eps\right\vert =o\left(\gie\right)$ in ${\mathbb D}_{\xie}\left(\sie\right)$. If $u_0\equiv 0$, the last term disappears. If $u_0\not\equiv 0$, then $\lambda_\eps\to \lambda_0$ with $\lambda_0>0$ and Claim \ref{claim-MT-2} gives that $\gie s_{i,\eps}^2=o(1)$. Thus, in any case, we have that 
$$\int_{{\mathbb D}_{\xie}\left(\sie\right)} \langle \nabla\left(w_\eps-u_\eps\right), \nabla u_\eps-2\nabla u_0 \rangle\, dx \le -4\pi +o(1)\hskip.1cm.$$
Coming back to \eqref{eq-claim-MT-5-1} with this proves the claim. \hfill $\diamondsuit$

\medskip The next two claims are devoted to obtaining good pointwise estimates on $u_\eps$ and $\nabla u_\eps$.

\begin{claim}\label{claim-MT-6}
For any sequence $\left(x_\eps\right)$ of points in $\Omega$ such that 
$$\frac{\left\vert \xe-\xie\right\vert}{\mie}\to +\infty\hbox{ as }\eps\to 0\hbox{ for }i=1,\dots,N\hskip.1cm,$$
we have that 

\smallskip a) if $d_\eps=d\left(x_\eps,\partial\Omega\right)\not\rightarrow 0$ as $\eps\to 0$, then 
\begincal
u_\eps\left(x_\eps\right)&=&\psi_\eps\left(x_\eps\right)+ \sum_{i=1}^N \frac{4\pi+o(1)}{\gie} {\mathcal G}\left(\xie,x_\eps\right) \\
&&+O\left(\sum_{i=1}^N \gie^{-1}\left(\frac{\mie}{\left\vert \xie-x_\eps\right\vert} +\gie^{-1} \ln\left(\frac{\sie}{\left\vert \xie-x_\eps\right\vert}+2\right)\right)\right)\hskip.1cm,
\fincal

\smallskip b) if $d_\eps\to 0$ as $\eps\to 0$, then
\begincal
u_\eps\left(x_\eps\right)&=&\psi_\eps\left(x_\eps\right)+ \sum_{i=1}^N \frac{4\pi+o(1)}{\gie} {\mathcal G}\left(\xie,x_\eps\right)   \\
&&+O\left( \sum_{i\in {\mathcal A}} \gie^{-1}\left(\frac{\mie}{\left\vert \xie-x_\eps\right\vert} +\gie^{-1} \ln\left(\frac{\sie}{\left\vert \xie-x_\eps\right\vert}+2\right)\right)\right)\\
&&+O\left(\sum_{i\in {\mathcal B}} \frac{d_\eps}{d_\eps+d_{i,\eps}} \gie^{-1}\left(\mie + \gie^{-1}\sie\right)\right)
\fincal

\smallskip\noindent where ${\mathcal G}$ is the Green function of the Laplacian with Dirichlet boundary condition in $\Omega$ and $\psi_\eps$ is a solution of 
$$\Delta \psi_\eps= \lambda_\eps f_\eps w_\eps e^{w_\eps^2}\hbox{ in }\Omega \hbox{ and }\psi_\eps=0\hbox{ on }\partial \Omega\hskip.1cm.$$
In b), ${\mathcal A}$ is defined as the set of $i\in \left\{1,\dots,N\right\}$ such that $\left\vert \xie-x_\eps\right\vert \le \sie +o\left(d_\eps\right)$ and ${\mathcal B}$ as its complementary. 
\end{claim}

\medskip {\it Proof of Claim \ref{claim-MT-6}} - We let ${\mathcal G}$ be the Green function of the Laplacian with Dirichlet boundary condition in $\Omega$. We let $\left(x_\eps\right)$ be a sequence of points in $\Omega$ such that 
\begin{equation}\label{eq-claim-MT-6-0}
\frac{\left\vert \xe-\xie\right\vert}{\mie}\to +\infty\hbox{ as }\eps\to 0\hbox{ for }i=1,\dots,N\hskip.1cm.
\end{equation}
Then we have thanks to \eqref{eq-MT-1} and to the definition of $\psi_\eps$ that 
$$u_\eps\left(\xe\right)-\psi_\eps\left(x_\eps\right)=\lambda_\eps \int_\Omega {\mathcal G}\left(\xe,x\right) f_\eps(x) \left(u_\eps(x)e^{u_\eps(x)^2}-w_\eps(x)e^{w_\eps(x)^2}\right)\, dx\hskip.1cm.$$
Using the definition \eqref{eq-MT-11} of $w_\eps$, this gives that 
\begin{eqnarray}\label{eq-claim-MT-6-1}
&&u_\eps\left(\xe\right)-\psi_\eps\left(x_\eps\right) \\
&&\quad =\sum_{i=1}^N \lambda_\eps {\mathcal G}\left(x_\eps,\xie\right)\int_{{\mathbb D}_{\xie}\left(\sie\right)}f_\eps(x) \left(u_\eps(x)e^{u_\eps(x)^2}-w_\eps(x)e^{w_\eps(x)^2}\right)\, dx+ \sum_{i=1}^N A_{i,\eps}\nonumber
\end{eqnarray}
where 
\begin{equation}\label{eq-claim-MT-6-2}
A_{i,\eps}= \lambda_\eps\int_{{\mathbb D}_{\xie}\left(\sie\right)}\left({\mathcal G}\left(\xe,x\right)-{\mathcal G}\left(x_\eps,\xie\right)\right) f_\eps(x) \left(u_\eps(x)e^{u_\eps(x)^2}-w_\eps(x)e^{w_\eps(x)^2}\right)\, dx\hskip.1cm.
\end{equation}
We fix $i\in \left\{1,\dots,N\right\}$ in the following and we let 
\begin{eqnarray}\label{eq-claim-MT-6-3}
\Omega_{0,\eps} &=& {\mathbb D}_{\xie}\left(\sie\right)\cap \left\{t_\eps(x)\le t_{1,\eps}\right\}, \nonumber\\
\Omega_{1,\eps} &= &{\mathbb D}_{\xie}\left(\sie\right)\cap \left\{t_{1,\eps}\le t_\eps(x)\le t_{2,\eps}\right\}\hbox{ and }\\
\Omega_{2,\eps} &=& {\mathbb D}_{\xie}\left(\sie\right)\cap \left\{t_\eps(x)\ge t_{2,\eps}\right\}\hskip.1cm.\nonumber
\end{eqnarray}
where $t_\eps(x)=\ln \left(1+\frac{\left\vert x-\xie\right\vert^2}{4\mie^2}\right)$, $t_{1,\eps} = \frac{1}{4}\gamma_{i,\eps}^2$ and $t_{2,\eps} = \gamma_{i,\eps}^2 -\gamma_{i,\eps}$.

\medskip {\sc Step 1} - {\it We have that 
$$\lambda_\eps \int_{{\mathbb D}_{\xie}\left(\sie\right)}f_\eps(x) \left(u_\eps(x)e^{u_\eps(x)^2}-w_\eps(x)e^{w_\eps(x)^2}\right)\, dx=4\pi \gie^{-1}+o\left(\gie^{-1}\right)\hskip.1cm.$$
}

\medskip {\sc Proof of Step 1} - We write that 
\begincal
&&\lambda_\eps \int_{{\mathbb D}_{\xie}\left(\sie\right)}f_\eps(x) \left(u_\eps(x)e^{u_\eps(x)^2}-w_\eps(x)e^{w_\eps(x)^2}\right)\, dx \\
&&\quad = \lambda_\eps \int_{{\mathbb D}_{\xie}\left(R_\eps\mie\right)}f_\eps(x) \left(u_\eps(x)e^{u_\eps(x)^2}-w_\eps(x)e^{w_\eps(x)^2}\right)\, dx\\
&&\qquad+\lambda_\eps \int_{\Omega_{0,\eps}\setminus {\mathbb D}_{\xie}\left(R_\eps\mie\right)}f_\eps(x) \left(u_\eps(x)e^{u_\eps(x)^2}-w_\eps(x)e^{w_\eps(x)^2}\right)\, dx\\
&&\qquad+\lambda_\eps \int_{\Omega_{1,\eps}}f_\eps(x) \left(u_\eps(x)e^{u_\eps(x)^2}-w_\eps(x)e^{w_\eps(x)^2}\right)\, dx\\
&&\qquad+\lambda_\eps \int_{\Omega_{2,\eps}}f_\eps(x) \left(u_\eps(x)e^{u_\eps(x)^2}-w_\eps(x)e^{w_\eps(x)^2}\right)\, dx
\fincal
where $R_\eps\to +\infty$ is such that $\left\vert u_\eps-B_{i,\eps}\right\vert =o\left(\gie^{-1}\right)$ and $\gie^{-1}B_{i,\eps}(x)=1+o(1)$ in ${\mathbb D}_{\xie}\left(R_\eps\mie\right)$. Such a $R_\eps$ does exist thanks to d) of Proposition \ref{prop-DMJ}. Then we have, using also \eqref{eq-MT-9}, that 
\begincal
&&\lambda_\eps \int_{{\mathbb D}_{\xie}\left(R_\eps\mie\right)}f_\eps(x) \left(u_\eps(x)e^{u_\eps(x)^2}-w_\eps(x)e^{w_\eps(x)^2}\right)\, dx\\
&&\quad = \lambda_\eps \left(f_\eps\left(\xie\right)+o(1)\right) \int_{{\mathbb D}_{\xie}\left(R_\eps\mie\right)} B_{i,\eps}(x) e^{B_{i,\eps}(x)^2}\, dx \\
&&\qquad + O\left(\lambda_\eps B_{i,\eps}\left(\sie\right) e^{B_{i,\eps}\left(\sie\right)^2} R_\eps^2 \mie^2\right)\\
&&\quad = \gie^{-1}\left(\int_{{\mathbb D}_0\left(R_\eps\right)} e^{2U}\, dx+o(1)\right)+o\left(\gie^{-1}\right)\\
&&\quad = 4\pi \gie^{-1}+o\left(\gie^{-1}\right)\hskip.1cm.
\fincal
In $\Omega_{0,\eps}$, we write that 
$$
B_{i,\eps}\left(x\right)^2 = \gie^2 - 2 t_\eps\left(x\right)+\frac{t_\eps(x)^2}{\gie^2} +O(1)
\le  \gie^2 - \frac{7}{4}t_\eps(x) + O(1)$$
so that 
$$e^{B_{i,\eps}(x)^2} \le e^{\gie^2} \left(1+\frac{\left\vert x-\xie\right\vert^2}{4\mie^2}\right)^{-\frac{7}{4}}\hskip.1cm.$$
Thus we can write that 
\begincal
0\le &&\lambda_\eps \int_{\Omega_{0,\eps}\setminus {\mathbb D}_{\xie}\left(R_\eps\mie\right)}f_\eps(x) \left(u_\eps(x)e^{u_\eps(x)^2}-w_\eps(x)e^{w_\eps(x)^2}\right)\, dx\\
&& \qquad \le C\gie^{-1} \mie^{-2}\int_{\Omega\setminus {\mathbb D}_{\xie}\left(R_\eps\mie\right)}\left(1+\frac{\left\vert x-\xie\right\vert^2}{4\mie^2}\right)^{-\frac{7}{4}}\, dx = o\left(\gie^{-1}\right) \hskip.1cm.
\fincal
In $\Omega_{1,\eps}$, we write that 
$$e^{B_{i,\eps}\left(x\right)^2} \le e^{\gie^2-\frac{1}{4}\gie} e^{-t_\eps(x)}$$
so that 
\begincal
0&\le &\lambda_\eps \int_{\Omega_{1,\eps}}f_\eps(x) \left(u_\eps(x)e^{u_\eps(x)^2}-w_\eps(x)e^{w_\eps(x)^2}\right)\, dx\\
&&\qquad \le C\mie^{-2} e^{-\frac{1}{4}\gie}\gie^{-1} \int_{\Omega_{1,\eps}}\left(1+\frac{\left\vert x-\xie\right\vert^2}{4\mie^2}\right)^{-1}\, dx\\
&&\qquad \le C e^{-\frac{1}{4}\gie}\gie^{-1}\ln \frac{\sie}{\mie} = o\left(\gie^{-2}\right)
\fincal 
since $2\ln \frac{\sie}{\mie}= \gie^2 +O(1)$ thanks to Claim \ref{claim-MT-2}. At last, in $\Omega_{2,\eps}$, we have that $B_{i,\eps}=O(1)$ so that 
$$
0\le \lambda_\eps \int_{\Omega_{2,\eps}}f_\eps(x) \left(u_\eps(x)e^{u_\eps(x)^2}-w_\eps(x)e^{w_\eps(x)^2}\right)\, dx\le C\lambda_\eps \sie^2 = O\left(\gie^{-2}\right)$$
thanks to Claim \ref{claim-MT-2}. Combining all these estimates clearly proves Step 1. \hfill $\spadesuit$

\medskip We shall now estimate the $A_i$'s involved in \eqref{eq-claim-MT-6-1} and defined in \eqref{eq-claim-MT-6-2}. We write since $u_\eps\ge w_\eps$ and thanks to \eqref{eq-MT-9} that 
\begin{equation}\label{eq-claim-MT-6-4}
\left\vert A_{i,\eps}\right\vert \le C\lambda_\eps \int_{{\mathbb D}_{\xie}\left(\sie\right)}\left\vert{\mathcal G}\left(\xe,x\right)-{\mathcal G}\left(x_\eps,\xie\right)\right\vert \left(B_{i,\eps} +C\gie^{-1}\right)e^{B_{i,\eps}(x)^2}\, dx\hskip.1cm.
\end{equation}

\medskip {\sc Step 2} - {\it Assume that $d_\eps=d\left(x_\eps,\partial\Omega\right)\ge d$ for some $d>0$. Then we have that 
$$\left\vert A_{i,\eps}\right\vert \le C\gie^{-1}\left(\frac{\mie}{\left\vert \xie-x_\eps\right\vert} +\gie^{-1}\ln\left(\frac{\sie}{\left\vert \xie-x_\eps\right\vert}+2\right)\right)\hskip.1cm.$$
}

\medskip {\sc Proof of Step 2} - We use \eqref{eq-appB-1} to write that 
$$\left\vert{\mathcal G}\left(\xe,x\right)-{\mathcal G}\left(x_\eps,\xie\right)\right\vert\le \frac{1}{2\pi}\left\vert \ln \frac{\left\vert x_{i,\eps}-x_\eps\right\vert}{\left\vert x_\eps-x\right\vert}\right\vert + C\left\vert x-\xie\right\vert\hskip.1cm.$$
Thus we have thanks to \eqref{eq-claim-MT-6-4} that 
$$\left\vert A_{i,\eps}\right\vert \le C\lambda_\eps \int_{{\mathbb D}_{\xie}\left(\sie\right)}\left(\left\vert\ln \frac{\left\vert x_{i,\eps}-x_\eps\right\vert}{\left\vert x_\eps-x\right\vert}\right\vert +\left\vert \xie-x\right\vert\right)\left(B_{i,\eps} +C\gie^{-1}\right)e^{B_{i,\eps}(x)^2}\, dx\hskip.1cm.$$
In $\Omega_{0,\eps}$, we have that 
$$B_{i,\eps}(x)\le \gie^2 - \frac{7}{4}\ln\left(1+\frac{\left\vert \xie-x\right\vert^2}{4\mie^2}\right)$$
so that 
\begincal 
&& \lambda_\eps \int_{\Omega_{0,\eps}}\left(\left\vert\ln \frac{\left\vert x_{i,\eps}-x_\eps\right\vert}{\left\vert x_\eps-x\right\vert}\right\vert +\left\vert \xie-x\right\vert\right)\left(B_{i,\eps} +C\gie^{-1}\right)e^{B_{i,\eps}(x)^2}\, dx\\
&&\quad \le C \mie^{-2}\gie^{-1}\int_{\Omega_{0,\eps}}\left(\left\vert\ln \frac{\left\vert x_{i,\eps}-x_\eps\right\vert}{\left\vert x_\eps-x\right\vert}\right\vert +\left\vert \xie-x\right\vert\right)\left(1+\frac{\left\vert \xie-x\right\vert^2}{4\mie^2}\right)^{-\frac{7}{4}}\, dx\hskip.1cm.
\fincal
This leads after simple computations, since $\frac{\left\vert \xie-x_\eps\right\vert}{\mie}\to +\infty$, as $\eps\to 0$ to 
$$\lambda_\eps \int_{\Omega_{0,\eps}}\left(\left\vert\ln \frac{\left\vert x_{i,\eps}-x_\eps\right\vert}{\left\vert x_\eps-x\right\vert}\right\vert +\left\vert \xie-x\right\vert\right)\left(B_{i,\eps} +C\gie^{-1}\right)e^{B_{i,\eps}(x)^2}\, dx \le C\gie^{-1}\frac{\mie}{\left\vert \xie-x_\eps\right\vert}\hskip.1cm.$$
In $\Omega_{1,\eps}$, we can write that 
$$e^{B_{i,\eps}(x)^2} \le e^{\gie^2-\frac{1}{4}\gie} \left(1+\frac{\left\vert x-\xie\right\vert^2}{4\mie^2}\right)^{-1}$$
so that 
\begincal
&&\lambda_\eps \int_{\Omega_{1,\eps}}\left(\left\vert\ln \frac{\left\vert x_{i,\eps}-x_\eps\right\vert}{\left\vert x_\eps-x\right\vert}\right\vert +\left\vert \xie-x\right\vert\right)\left(B_{i,\eps} +C\gie^{-1}\right)e^{B_{i,\eps}(x)^2}\, dx\\
&&\quad \le C\gie^{-1} e^{-\frac{1}{4}\gie} \int_{\Omega_{1,\eps}}\left(\left\vert\ln \frac{\left\vert x_{i,\eps}-x_\eps\right\vert}{\left\vert x_\eps-x\right\vert}\right\vert +\left\vert \xie-x\right\vert\right)\left\vert x-\xie\right\vert^{-2}\, dx\\
&&\quad \le C\gie^{-1}e^{-\frac{1}{4}\gie} \left(\ln \frac{r_{2,\eps}}{\mu_{i,\eps}}\right)^2
\fincal
where $t_{i,\eps}\left(r_{2,\eps}\right)=t_{2,\eps}$. We have that 
$$\ln \frac{r_{2,\eps}}{\mu_{i,\eps}} \le C\gamma_\eps^2$$
so that 
$$\lambda_\eps \int_{\Omega_{1,\eps}}\left(\left\vert\ln \frac{\left\vert x_{i,\eps}-x_\eps\right\vert}{\left\vert x_\eps-x\right\vert}\right\vert +\left\vert \xie-x\right\vert\right)\left(B_{i,\eps} +C\gie^{-1}\right)e^{B_{i,\eps}(x)^2}\, dx=O\left(\gie^3 e^{-\frac{1}{4}\gie}\right)=o\left(\gie^{-2}\right)\hskip.1cm.$$
At last, in $\Omega_{2,\eps}$, we have that $B_{i,\eps}=O(1)$ so that 
\begincal
&&\lambda_\eps \int_{\Omega_{2,\eps}}\left(\left\vert\ln \frac{\left\vert x_{i,\eps}-x_\eps\right\vert}{\left\vert x_\eps-x\right\vert}\right\vert +\left\vert \xie-x\right\vert\right)\left(B_{i,\eps} +C\gie^{-1}\right)e^{B_{i,\eps}(x)^2}\, dx\\
&&\quad \le \lambda_\eps  \int_{\Omega_{2,\eps}}\left(\left\vert\ln \frac{\left\vert x_{i,\eps}-x_\eps\right\vert}{\left\vert x_\eps-x\right\vert}\right\vert +\left\vert \xie-x\right\vert\right)\, dx\\
&&\quad \le \lambda_\eps \sie^2 \ln \left(\frac{\sie}{\left\vert \xie-x_\eps\right\vert}+2\right)\le C\gie^{-2}\ln \left(\frac{\sie}{\left\vert \xie-x_\eps\right\vert}+2\right)
\fincal
by direct computations and Claim \ref{claim-MT-2}. Combining the above estimates gives Step 2. \hfill $\spadesuit$

\medskip {\sc Step 3} - {\it Assume now that $d_\eps=d\left(x_\eps,\partial\Omega\right)\to 0$ as $\eps\to 0$ and that $\left\vert \xie-x_\eps\right\vert \ge \sie+\delta d_\eps$ for some $\delta>0$. Then we have that 
$$\left\vert A_{i,\eps}\right\vert \le C\frac{d_\eps}{d_\eps+d_{i,\eps}} \left(\gie^{-1}\mie + \gie^{-2}\sie\right)\hskip.1cm.$$
}

\medskip {\sc Proof of Step 3} - In this setting, we can apply \eqref{eq-appB-11} to write that 
$$\left\vert {\mathcal G}\left(x_\eps,x\right)- {\mathcal G}\left(x_\eps,\xie\right)\right\vert \le C\frac{d_\eps}{d_\eps+d_{i,\eps}}\left\vert x-\xie\right\vert$$
so that 
$$\left\vert A_{i,\eps}\right\vert \le C \lambda_\eps \frac{d_\eps}{d_\eps+d_{i,\eps}}\int_{{\mathbb D}_{\xie}\left(\sie\right)} \left\vert x-\xie\right\vert \left(B_{i,\eps}+C\gie^{-1}\right) e^{B_{i,\eps}^2}\, dx$$
thanks to \eqref{eq-claim-MT-6-2}. In $\Omega_{0,\eps}$, we have that 
$$\left(B_{i,\eps}+C\gie^{-1}\right) e^{B_{i,\eps}^2}\le C\gie e^{\gie^2} \left(1+\frac{\left\vert x-\xie\right\vert^2}{4\mie^2}\right)^{-\frac{7}{4}}$$
so that 
\begincal
&&\lambda_\eps\int_{\Omega_{0,\eps}} \left\vert x-\xie\right\vert \left(B_{i,\eps}+C\gie^{-1}\right) e^{B_{i,\eps}^2}\, dx\\
&&\quad \le C\mie^{-2}\gie^{-1} \int_{{\mathbb D}_{\xie}\left(\sie\right)}\left\vert x-\xie\right\vert \left(1+\frac{\left\vert x-\xie\right\vert^2}{4\mie^2}\right)^{-\frac{7}{4}}\,dx\\
&&\quad \le C\mie\gie^{-1}\hskip.1cm.
\fincal
In $\Omega_{1,\eps}$, we have that 
$$\left(B_{i,\eps}+C\gie^{-1}\right) e^{B_{i,\eps}(x)^2} \le \gie e^{\gie^2-\frac{1}{4}\gie} \left(1+\frac{\left\vert x-\xie\right\vert^2}{4\mie^2}\right)^{-1}$$
so that 
\begincal
&&\lambda_\eps\int_{\Omega_{1,\eps}} \left\vert x-\xie\right\vert \left(B_{i,\eps}+C\gie^{-1}\right) e^{B_{i,\eps}^2}\, dx\\
&&\quad \le C\mie^{-2}\gie^{-1}e^{-\frac{1}{4}\gie} \int_{{\mathbb D}_{\xie}\left(\sie\right)}\left\vert x-\xie\right\vert \left(1+\frac{\left\vert x-\xie\right\vert^2}{4\mie^2}\right)^{-1}\,dx\\
&&\quad \le \gie^{-1}e^{-\frac{1}{4}\gie}\sie\hskip.1cm.
\fincal
At last, in $\Omega_{2,\eps}$, we have that $B_{i,\eps}=O(1)$ so that 
$$\lambda_\eps\int_{\Omega_{2,\eps}} \left\vert x-\xie\right\vert \left(B_{i,\eps}+C\gie^{-1}\right) e^{B_{i,\eps}^2}\, dx \le \lambda_\eps s_{i,\eps}^3\hskip.1cm.$$
Combining the above estimates with Claim \ref{claim-MT-2}, we get the estimate of Step 3.\hfill $\spadesuit$

\medskip {\sc Step 4} - {\it Assume now that $d_\eps=d\left(x_\eps,\partial\Omega\right)\to 0$ as $\eps\to 0$ and that $\left\vert \xie-x_\eps\right\vert \le \sie+o\left(d_\eps\right)$. Then we have that 
$$\left\vert A_{i,\eps}\right\vert \le C \gie^{-1} \left(\frac{\mie}{\left\vert \xie-x_\eps\right\vert} + \gie^{-1} \ln\left(\frac{\sie}{\left\vert \xie-x_\eps\right\vert}+2\right)\right)\hskip.1cm.$$
}

\medskip {\sc Proof of Step 4} - Let us remark that in this case, we necessarily have that 
$$d_{\eps}\le \left\vert \xie-x_\eps\right\vert + d_{i,\eps} \le \sie +d_{i,\eps} +o\left(d_\eps\right)\le \frac{3}{2}d_{i,\eps}+o\left(d_\eps\right)$$ 
so that $d_\eps = O\left(d_{i,\eps}\right)$. This leads in turn to $\left\vert \xie-x_\eps\right\vert \le \sie+o\left(d_{i,\eps}\right)$. And then we can write that
$$d_{i,\eps}\le d_\eps + \left\vert \xie-x_\eps\right\vert \le \sie+o\left(d_{i,\eps}\right) + d_\eps \le \frac{1}{2}d_{i,\eps} + d_\eps + o\left(d_{i,\eps}\right)$$
so that $d_{i,\eps}=O\left(d_\eps\right)$. Thanks to \eqref{eq-appB-11}, we can write that 
$$\left\vert {\mathcal G}\left(x_\eps,x\right)- {\mathcal G}\left(x_\eps,\xie\right)\right\vert \le C\frac{\left\vert x-\xie\right\vert}{d_{i,\eps}} + C \left\vert \ln \frac{\left\vert \xie-x_\eps\right\vert}{\left\vert x_\eps-x\right\vert}\right\vert$$
so that the computations of Step 2 lead to the result of Step 4. \hfill $\spadesuit$

\medskip Of course, the combination of Steps 1 to 4 gives the estimate of the claim. \hfill $\diamondsuit$

\begin{claim}\label{claim-MT-7}
There exists $C>0$ such that 
$$\left\vert \nabla \left(u_\eps-\psi_\eps\right)(x)\right\vert \le C\sum_{i=1}^N \gie^{-1}\left(\mie+\left\vert x-\xie\right\vert\right)^{-1} $$
where $\psi_\eps$ is as in Claim \ref{claim-MT-6}.
\end{claim}

\medskip {\it Proof of Claim \ref{claim-MT-7}} - We use again the Green representation formula with equation \eqref{eq-MT-1} (together with the equation satisfied by $\psi_\eps$, see Claim \ref{claim-MT-6}) to write that 
$$\left\vert \nabla \left(u_\eps-\psi_\eps\right)(x)\right\vert \le \lambda_\eps \int_\Omega \left\vert \nabla {\mathcal G}(x,y)\right\vert f_\eps(y)\left(u_\eps(y)e^{u_\eps(y)^2}-w_\eps(y)e^{w_\eps(y)^2}\right)\, dy\hskip.1cm.$$
Thanks to standard estimates on the Green function and to the definition \eqref{eq-MT-11}, this leads to 
\begin{equation}\label{eq-claim-MT-7-1}
\left\vert \nabla \left(u_\eps-\psi_\eps\right)(x)\right\vert \le C\sum_{i=1}^N \lambda_\eps \int_{{\mathbb D}_{\xie}\left(\sie\right)} \left\vert x-y\right\vert^{-1} u_\eps(y)e^{u_\eps(y)^2}\, dy\hskip.1cm.
\end{equation}
Thanks to \eqref{eq-MT-9}, we have that 
\begin{eqnarray}\label{eq-claim-MT-7-3}
&&\lambda_\eps \int_{{\mathbb D}_{\xie}\left(\sie\right)} \left\vert x-y\right\vert^{-1} u_\eps(y)e^{u_\eps(y)^2}\, dy\\
&&\quad \le C\lambda_\eps  \sum_{k=0,1,2} \int_{\Omega_{k,\eps}} \left\vert x-y\right\vert^{-1} \left(B_{i,\eps}(y)+C_i\gie^{-1}\right)e^{B_{i,\eps}(y)^2}\, dy\nonumber
\end{eqnarray}
where the $\Omega_{\alpha,\eps}$'s are as in \eqref{eq-claim-MT-6-3}. In $\Omega_{0,\eps}$, we write that 
$$\left(B_{i,\eps}+C_i\gie^{-1}\right)\left(y\right)e^{B_{i,\eps}(y)^2}\le C \gie e^{\gie^2}\left(1+\frac{\left\vert y-\xie\right\vert^2}{4\mie^2}\right)^{-\frac{7}{4}}$$
so that 
\begincal
&&\lambda_\eps\int_{\Omega_{0,\eps}} \left\vert x-y\right\vert^{-1} \left(B_{i,\eps}(y)+C_i\gie^{-1}\right)e^{B_{i,\eps}(y)^2}\, dy\\
&&\quad \le C\mie^{-2}\gie^{-1}\int_{\Omega_{0,\eps}}\left\vert x-y\right\vert^{-1}\left(1+\frac{\left\vert y-\xie\right\vert^2}{4\mie^2}\right)^{-\frac{7}{4}}\, dy\hskip.1cm.
\fincal
Direct computations give that 
\begin{equation}\label{eq-claim-MT-7-5}
\lambda_\eps\int_{\Omega_{0,\eps}} \left\vert x-y\right\vert^{-1} \left(B_{i,\eps}(y)+C_i\gie^{-1}\right)e^{B_{i,\eps}(y)^2}\, dy\le C\gie^{-1}\left(\mie+\left\vert x-\xie\right\vert\right)^{-1}\hskip.1cm.
\end{equation}
In $\Omega_{1,\eps}$, we write that 
$$\left(B_{i,\eps}+C_i\gie^{-1}\right)\left(y\right)e^{B_{i,\eps}(y)^2}\le C\gie e^{\gie^2} e^{\frac{t_\eps(y)^2}{\gie^2}-2t_\eps(y)} $$
so that 
$$\lambda_\eps\int_{\Omega_{1,\eps}} \left\vert x-y\right\vert^{-1} \left(B_{i,\eps}(y)+C_i\gie^{-1}\right)e^{B_{i,\eps}(y)^2}\, dy\le C\gie^{-1}\mie^{-2} \int_{\Omega_{1,\eps}} \left\vert x-y\right\vert^{-1} e^{\frac{t_\eps(y)^2}{\gie^2}-2t_\eps(y)}\hskip.1cm.$$
In $\Omega_{1,\eps}$, we have that 
$$\frac{t_\eps(y)^2}{\gie^2}-2t_\eps(y) \le -t_\eps(y) -\frac{1}{4}\gie$$
so that 
\begincal
&&\lambda_\eps\int_{\Omega_{1,\eps}} \left\vert x-y\right\vert^{-1} \left(B_{i,\eps}(y)+C_i\gie^{-1}\right)e^{B_{i,\eps}(y)^2}\, dy\\
&&\quad \le C\gie^{-1}e^{-\frac{1}{4}\gie} \int_{\Omega_{1,\eps}} \left\vert x-y\right\vert^{-1} \left(\mie^2+\left\vert y-\xie\right\vert^2\right)^{-1}\, dy\hskip.1cm.
\fincal
In $\Omega_{1,\eps}$ we have that $\left\vert y-\xie\right\vert \ge \mie$ so that 
$$\lambda_\eps\int_{\Omega_{1,\eps}} \left\vert x-y\right\vert^{-1} \left(B_{i,\eps}(y)+C_i\gie^{-1}\right)e^{B_{i,\eps}(y)^2}\, dy\le C\gie^{-1}e^{-\frac{1}{4}\gie} \int_{\Omega_{1,\eps}}
\left\vert x-y\right\vert^{-1} \left\vert y-\xie\right\vert^{-2}\, dy\hskip.1cm.$$
Noting that ${\mathbb D}_{\xie}\left(r_{1,\eps}\right)\cap \Omega_{1,\eps}=\emptyset$ for $\eps$ small where 
$$r_{1,\eps} = \mie e^{\frac{1}{8}\gie^2}\hskip.1cm,$$
we get by direct computations that 
\begincal
&&\lambda_\eps\int_{\Omega_{1,\eps}} \left\vert x-y\right\vert^{-1} \left(B_{i,\eps}(y)+C_i\gie^{-1}\right)e^{B_{i,\eps}(y)^2}\, dy\\
&&\quad \le C\gie^{-1}e^{-\frac{1}{4}\gie}\left(\left\vert x-\xie\right\vert+r_{1,\eps}\right)^{-1}\ln \left(2+\frac{\left\vert x-\xie\right\vert}{r_{1,\eps}}\right)\hskip.1cm.
\fincal
Thanks to the value of $r_{1,\eps}$, this leads to 
\begin{equation}\label{eq-claim-MT-7-6}
\lambda_\eps\int_{\Omega_{1,\eps}} \left\vert x-y\right\vert^{-1} \left(B_{i,\eps}(y)+C_i\gie^{-1}\right)e^{B_{i,\eps}(y)^2}\, dy
=o\left(\gie^{-1} \left(\mie+\left\vert x-\xie\right\vert\right)^{-1}\right) \hskip.1cm.
\end{equation}
At last, in $\Omega_{2,\eps}$, we have that $B_\eps\left(y\right)=O(1)$ so that 
$$\lambda_\eps\int_{\Omega_{2,\eps}} \left\vert x-y\right\vert^{-1} \left(B_{i,\eps}(y)+C_i\gie^{-1}\right)e^{B_{i,\eps}(y)^2}\, dy\le C \lambda_\eps \int_{\Omega_{2,\eps}} \left\vert x-y\right\vert^{-1}\, dy \le C\lambda_\eps \frac{s_{i,\eps}^2}{\sie +\left\vert x-\xie\right\vert}\hskip.1cm.$$
Thanks to Claim \ref{claim-MT-2}, this leads to 
\begin{equation}\label{eq-claim-MT-7-7}
\lambda_\eps\int_{\Omega_{2,\eps}} \left\vert x-y\right\vert^{-1} \left(B_{i,\eps}(y)+C_i\gie^{-1}\right)e^{B_{i,\eps}(y)^2}\, dy\le C \gie^{-2} \left(\sie +\left\vert x-\xie\right\vert\right)^{-1}\hskip.1cm.
\end{equation}
Coming back to \eqref{eq-claim-MT-7-1} with \eqref{eq-claim-MT-7-3}, \eqref{eq-claim-MT-7-5}, \eqref{eq-claim-MT-7-6} and \eqref{eq-claim-MT-7-7}, we obtain the claim. \hfill $\diamondsuit$

\medskip Let us reorder the concentration points in a suitable way. For this purpose, we notice that, up to a subsequence, for any $i,j\in \left\{1,\dots,N\right\}$, there exists $C_{i,j}$, possibly $0$ or $+\infty$ (but nonnegative) such that 
\begin{equation}\label{eq-MT-12}
\lim_{\eps\to 0}\frac{\gie}{\gje} = C_{i,j}\hskip.1cm.
\end{equation}
Note that $C_{i,j}=C_{j,i}^{-1}$ (with obvious conventions when $C_{i,j}=0$ or $+\infty$). Then there exists $\tilde{C}\ge 1$ such that 
\begin{equation}\label{eq-MT-13}
\hbox{ for any }i,j\in\left\{1,\dots,N\right\}, \hbox{ either }C_{i,j}=0\hbox{ or }C_{i,j}=+\infty\hbox{ or }\frac{1}{\tilde{C}}\le C_{i,j}\le \tilde{C}\hskip.1cm.
\end{equation}
It is then easily checked that we can order the concentration points in such a way that 
\begin{equation}\label{eq-MT-14}
\hbox{ for any }i,j\in\left\{1,\dots,N\right\}, \, i<j \Rightarrow C_{i,j}<+\infty
\end{equation}
and 
\begin{equation}\label{eq-MT-15}
\hbox{ for any }i,j\in\left\{1,\dots,N\right\}, \, i<j \hbox{ and }C_{i,j}>0 \Rightarrow \rie\le \rje\hskip.1cm.
\end{equation}

\medskip Let us give some estimates on $\psi_\eps$, involved in Claims \ref{claim-MT-6} and \ref{claim-MT-7}. Using Claim \ref{claim-MT-5}, we clearly have that $\lambda_\eps^{-1}\Delta \psi_\eps$ is uniformly bounded in any $L^p\left(\Omega\right)$ thanks to Trudinger-Moser inequality. Thus we know that there exists $C>0$ such that 
\begin{equation}\label{eq-MT-16}
\left\Vert \psi_\eps\right\Vert_{{\mathcal C}^{1,\alpha}\left(\overline{\Omega}\right)} \le C\lambda_\eps
\end{equation}
for $0<\alpha<1$ by standard elliptic theory. Now, if $\lambda_\eps\to 0$, we know that $u_0\equiv 0$ and we can be a little bit more precise. Indeed, 
$$\left\Vert \Delta\psi_\eps\right\Vert_{L^p\left(\Omega\right)} \le \lambda_\eps \left\Vert f_\eps\right\Vert_{L^\infty\left(\Omega\right)} \left\Vert w_\eps e^{w_\eps^2}\right\Vert_{L^p\left(\Omega\right)} \le \lambda_\eps \left\Vert f_\eps\right\Vert_{L^\infty\left(\Omega\right)} \left\Vert u_\eps\right\Vert_{L^{2p}\left(\Omega\right)}\left\Vert e^{w_\eps^2}\right\Vert_{L^{2p}\left(\Omega\right)}\hskip.1cm.$$
Since $u_0\equiv 0$, we know thanks to Claim \ref{claim-MT-5} and to Trudinger-Moser inequality that 
$\left(e^{w_\eps^2}\right)$ is bounded in any $L^q$. Thus we have that 
$$\left\Vert \Delta\psi_\eps\right\Vert_{L^p\left(\Omega\right)} \le C\lambda_\eps \left\Vert u_\eps\right\Vert_{L^{2p}\left(\Omega\right)}$$
thanks to \eqref{eq-cv-heps}. Using Claim \eqref{claim-MT-7}, we get that 
$$\left\Vert \nabla \left(u_\eps-\psi_\eps\right)\right\Vert_{L^q\left(\Omega\right)} \le \frac{C_q}{\gamma_{1,\eps}}$$
for some $C_q>0$ for all $1\le q<2$. Remember that concentration points are ordered such that \eqref{eq-MT-14} holds. This gives that
$$\left\Vert u_\eps\right\Vert_{L^{2p}\left(\Omega\right)} \le C_p\left(\gamma_{1,\eps}^{-1} + \left\Vert \nabla \psi_\eps\right\Vert_{{\mathcal C}^1\left(\overline{\Omega}\right)}\right)$$
so that 
$$\left\Vert \Delta\psi_\eps\right\Vert_{L^p\left(\Omega\right)} \le C\lambda_\eps \left(\gamma_{1,\eps}^{-1} + \left\Vert \nabla \psi_\eps\right\Vert_{{\mathcal C}^1\left(\overline{\Omega}\right)}\right)\hskip.1cm.$$
By standard elliptic theory and since we assumed that $\lambda_\eps\to 0$, we finally obtain that 
\begin{equation}\label{eq-MT-17}
\hbox{if }\lambda_\eps\to 0\hbox{ as }\eps\to 0,\hbox{ then } \left\Vert \psi_\eps\right\Vert_{{\mathcal C}^{1,\alpha}\left(\overline{\Omega}\right)} \le C\frac{\lambda_\eps}{\gamma_{1,\eps}}\hskip.1cm.
\end{equation}

\begin{claim}\label{claim-MT-8}
We have that $r_{1,\eps}\ge \delta_0$ for some $\delta_0>0$.
\end{claim}

\medskip {\it Proof of Claim \ref{claim-MT-8}} - We assume by contradiction that $r_{1,\eps}\to 0$ as $\eps\to 0$. We let in the following 
\begin{equation}\label{eq-claim-MT-8-1}
{\mathcal D}_1^\star = \left\{i\in \left\{2,\dots,N\right\}\hbox{ s.t. }\left\vert \xie-x_{1,\eps}\right\vert = O\left(r_{1,\eps}\right)\right\}\hbox{ and }{\mathcal D}_1={\mathcal D}_1^\star \cup \left\{1\right\}\hskip.1cm.
\end{equation}
After passing to a subsequence, we let 
\begin{equation}\label{eq-claim-MT-8-2}
{\mathcal S}_1^\star = \left\{\tilde{x}_i=\lim_{\eps\to 0}\frac{\xie-x_{1,\eps}}{r_{1,\eps}},\, i\in {\mathcal D}_1^\star\right\}\hbox{ and }{\mathcal S}_1={\mathcal S}_1^\star\cup\left\{\tilde{x}_1=0\right\}\hskip.1cm.
\end{equation}
We also let 
\begin{equation}\label{eq-claim-MT-8-3}
\Omega_{1,\eps} = \left\{y\in {\mathbb R}^2\hbox{ s.t. }x_{1,\eps}+r_{1,\eps}y \in \Omega\right\}\hskip.1cm.
\end{equation}
Note that, after passing to a subsequence (and up to a harmless rotation if necessary), we have that 
\begin{equation}\label{eq-claim-MT-8-4}
\Omega_{1,\eps} \to \Omega_0\hbox{ as }\eps\to 0\hbox{ where }\left\{
\begin{array}{ll}
{\ds \Omega_0={\mathbb R}^2}&{\ds \hbox{ if }\frac{d_{1,\eps}}{r_{1,\eps}}\to +\infty\hbox{ as }\eps\to 0}\\
\,&\,\\
{\ds\Omega_0={\mathbb R}\times\left(-\infty,L\right)}&{\ds \hbox{ if } \frac{d_{1,\eps}}{r_{1,\eps}}\to L\hbox{ as }\eps\to 0}
\end{array}\right.
\end{equation}
Here $d_{1,\eps}= d\left(x_{1,\eps},\partial\Omega\right)$, as defined in \eqref{eq-MT-8}. For $R>0$, we shall also let 
\begin{equation}\label{eq-claim-MT-8-5}
\Omega_0^R = \left(\Omega_0\cap {\mathbb D}_0(R)\right)\setminus \bigcup_{i\in {\mathcal D}_1} {\mathbb D}_{\tilde{x}_i}\left(\frac{1}{R}\right)\hskip.1cm.
\end{equation}

We shall distinguish three cases, depending on the behaviour of $d_{1,\eps}=d\left(x_{1,\eps},\partial\Omega\right)$ and $r_{1,\eps}$. 

\smallskip {\sc Case 1} - We assume that $d_{1,\eps}\not \to 0$ as $\eps\to 0$, meaning that, after passing to a subsequence, $x_{1,\eps}\to x_1$ as $\eps\to 0$ with $x_1\in \Omega$. 

\smallskip\noindent We let $y\in \Omega_0^R$ for some $R>0$ and we set $x_\eps=x_{1,\eps}+r_{1,\eps}y$. Since $d_{1,\eps}\not\to 0$ and $r_{1,\eps}\to 0$, we are in situation a) of Claim \ref{claim-MT-6}. Note indeed that 
$$\frac{\left\vert x_\eps-\xie\right\vert}{\mie}\to +\infty\hbox{ as }\eps\to 0\hbox{ for all }i=1,\dots,N\hskip.1cm.$$
It is obvious if $i\in{\mathcal D}_1$ since we clearly have in this case 
$$\frac{\left\vert x_\eps-\xie\right\vert}{\mie} = \frac{\left\vert x_\eps-\xie\right\vert}{r_{1,\eps}}\frac{r_{1,\eps}}{\rie}\frac{\rie}{\mie}$$
with ${\ds \frac{\left\vert x_\eps-\xie\right\vert}{r_{1,\eps}}\ge R^{-1}+o(1)}$, ${\ds \frac{r_{1,\eps}}{\rie}\ge 2\left\vert \tilde{x}_i\right\vert^{-1}+o(1)}$ for $i\in {\mathcal D}_1^\star$ and equal to $1$ if $i=1$, and ${\ds \frac{\rie}{\mie}\to +\infty}$ as $\eps\to 0$ thanks to assertion c) of Proposition \ref{prop-DMJ}. While, if $i\not\in{\mathcal D}_1$, we can write that 
$$\frac{\left\vert x_\eps-\xie\right\vert}{\mie} \ge \left(1+o(1)\right) \frac{\left\vert \xie-x_{1,\eps}\right\vert}{\mie}\ge  \left(2+o(1)\right) \frac{\rie}{\mie}\to +\infty\hbox{ as }\eps\to 0\hskip.1cm.$$
Thus, applying a) of Claim \ref{claim-MT-6}, we can write that 
\begincal
u_\eps\left(x_\eps\right)&=&\psi_\eps\left(x_\eps\right)+\sum_{i=1}^N \left(4\pi+o(1)\right)\gie^{-1}{\mathcal G}\left(\xie,x_\eps\right)\\
&&+O\left(\sum_{i=1}^N \gie^{-1}\left(\frac{\mie}{\left\vert \xie-x_\eps\right\vert} +\gie^{-1} \ln\left(\frac{\sie}{\left\vert \xie-x_\eps\right\vert}+2\right)\right)\right)\hskip.1cm.
\fincal
Now, for any $i\in\left\{1,\dots,N\right\}$,
$$\gie^{-1}\frac{\mie}{\left\vert \xie-x_\eps\right\vert}  =o\left(\gie^{-1}\right)=o\left(\gamma_{1,\eps}^{-1}\right)$$
thanks to \eqref{eq-MT-14} and 
$$\gie^{-2} \ln\left(\frac{\sie}{\left\vert \xie-x_\eps\right\vert}+2\right) =o\left(\gamma_{1,\eps}^{-1}\right)$$
thanks to the fact that $\sie\le \rie =O\left(\left\vert \xie-x_\eps\right\vert\right)$. Note that Claim \ref{claim-MT-2} implies that $\lambda_\eps=O\left(\gamma_{1,\eps}^{-2}\right)$ in our case so that \eqref{eq-MT-17} gives that 
$$\psi_\eps\left(x_\eps\right)=O\left(\gamma_{1,\eps}^{-3}\right)=o\left(\gamma_{1,\eps}^{-1}\right)\hskip.1cm.$$
Thus we have that
\begin{equation}\label{eq-claim-MT-8-6}
u_\eps\left(x_\eps\right)=\sum_{i=1}^N \left(4\pi+o(1)\right)\gie^{-1}{\mathcal G}\left(\xie,x_\eps\right)+o\left(\gamma_{1,\eps}^{-1}\right)\hskip.1cm.
\end{equation}
We can now use \eqref{eq-appB-2} to write that 
$${\mathcal G}\left(\xie,x_\eps\right) = \frac{1}{2\pi}\ln \frac{1}{r_{1,\eps}} + O\left(1\right)$$
if $i\in {\mathcal D}_1$ and that 
$${\mathcal G}\left(\xie,x_\eps\right) = {\mathcal G}\left(\xie,x_{1,\eps}\right) + O(1)$$
if $i\not\in{\mathcal D}_1$. Thus we have that 
\begin{eqnarray}\label{eq-claim_MT-8-7}
u_\eps\left(x_\eps\right)&=&\left(2+o(1)\right)\gamma_{1,\eps}^{-1}\ln\frac{1}{r_{1,\eps}}\left(1+\sum_{i\in{\mathcal D}_1^\star} C_{1,i}\right)\\
&&+ \sum_{i\not\in {\mathcal D}_1} \left(4\pi+o(1)\right)\gie^{-1}{\mathcal G}\left(\xie,x_{1,\eps}\right)+O\left(\gamma_{1,\eps}^{-1}\right)\hskip.1cm.\nonumber
\end{eqnarray}
Note that $C_{1,i}\le \tilde{C}$ for all $i>1$ thanks to \eqref{eq-MT-13}. Thus we have in particular that 
\begin{equation}\label{eq-claim-MT-8-8}
\gamma_{1,\eps}^{-1}\ln\frac{1}{r_{1,\eps}}\left(1+\sum_{i\in{\mathcal D}_1^\star} C_{1,i}\right)\le \left(\frac{1}{2}+o(1)\right)u_\eps\left(x_\eps\right) \le 
\gamma_{1,\eps}^{-1}\ln\frac{1}{r_{1,\eps}}\left(1+(N-1)\tilde{C}\right) \hskip.1cm.
\end{equation}
Note that we also have thanks to Claim \ref{claim-MT-7} and to \eqref{eq-MT-17} that 
\begin{equation}\label{eq-claim-MT-8-9}
\left\vert \nabla u_\eps(x)\right\vert \le C\gamma_{1,\eps}^{-1}\left\vert x_{1,\eps}-x\right\vert^{-1} \hbox{ for all }x\in {\mathbb D}_{x_{1,\eps}}\left(r_{1,\eps}\right)\hskip.1cm.
\end{equation}
We are thus in position to apply Claim \ref{claim-MT-LBA} for $i=1$ to write that, if $\left\vert x\right\vert=\frac{1}{2}$,
$$u_\eps\left(x_{1,\eps}+r_{1,\eps} x\right)= B_{1,\eps}\left(r_{1,\eps}\right)+O\left(\gamma_{1,\eps}^{-1}\right)\hskip.1cm.$$
Combined with \eqref{eq-claim-MT-8-8}, this gives that 
\begin{equation}\label{eq-claim-MT-8-10}
\left(2+o(1)\right)\gamma_{1,\eps}^{-1}\ln\frac{1}{r_{1,\eps}}\left(1+\sum_{i\in{\mathcal D}_1^\star} C_{1,i}\right)\le B_{1,\eps}\left(r_{1,\eps}\right)\le \left(2+o(1)\right)\gamma_{1,\eps}^{-1}\ln\frac{1}{r_{1,\eps}}\left(1+(N-1)\tilde{C}\right) \hskip.1cm.
\end{equation}
We write now thanks to Claim \ref{claim-appA-4} of Appendix A that 
\begin{equation}\label{eq-claim-MT-8-10bis}
B_{1,\eps}\left(r_{1,\eps}\right)= 2\gamma_{1,\eps}^{-1} \ln \frac{1}{r_{1,\eps}} - \gamma_{1,\eps}^{-1} \ln \left(\lambda_\eps\gamma_{1,\eps}^2\right) + O\left(\gamma_{1,\eps}^{-1}\right)
\end{equation}
to deduce that 
\begin{equation}\label{eq-claim-MT-8-11}
\left(2+o(1)\right)\gamma_{1,\eps}^{-1}\ln\frac{1}{r_{1,\eps}}\sum_{i\in{\mathcal D}_1^\star} C_{1,i}\le - \gamma_{1,\eps}^{-1} \ln \left(\lambda_\eps\gamma_{1,\eps}^2\right)\le \left(2(N-1)\tilde{C}+o(1)\right)\gamma_{1,\eps}^{-1}\ln\frac{1}{r_{1,\eps}}\hskip.1cm.
\end{equation}
Fix now $i \in {\mathcal D}_1^\star$. It is clear that there exists $\delta>0$ such that $\partial {\mathbb D}_{\xie}\left(\delta r_{1,\eps}\right)\subset \left\{x_{1,\eps}+r_{1,\eps} y,\, y\in 
 \Omega_0^R\right\}$ for some $R>0$. Thus we can write that 
$$\inf_{\partial {\mathbb D}_{\xie}\left(\delta r_{1,\eps}\right)} u_\eps\ge \left(2+o(1)\right)\gamma_{1,\eps}^{-1}\ln\frac{1}{r_{1,\eps}}\left(1+\sum_{i\in{\mathcal D}_1^\star} C_{1,i}\right)$$
thanks to \eqref{eq-claim-MT-8-8}. We can also apply b) of Claim \ref{claim-MT-0} with $r_\eps=\delta r_{1,\eps}$ thanks to Claim \ref{claim-MT-1} and to the fact that ${\ds \frac{r_{1,\eps}}{\rie}\ge 2\left\vert \tilde{x}_i\right\vert^{-1}}+o(1)$. This leads to 
$$\left(2+o(1)\right)\gamma_{1,\eps}^{-1}\ln\frac{1}{r_{1,\eps}}\left(1+\sum_{i\in{\mathcal D}_1^\star} C_{1,i}\right)\le B_{i,\eps}\left(\delta r_{1,\eps}\right)+o\left(\gie^{-1}\right)\hskip.1cm.$$
We have that 
$$B_{i,\eps}\left(\delta r_{1,\eps}\right) = 2\gie^{-1} \ln \frac{1}{r_{1,\eps}} - \gie^{-1} \ln \left(\lambda_\eps\gie^2\right) + O\left(\gie^{-1}\right)\hskip.1cm.$$
This leads together with \eqref{eq-claim-MT-8-11} to 
$$\left(2+o(1)\right)\gamma_{1,\eps}^{-1}\ln\frac{1}{r_{1,\eps}}\left(1+\sum_{i\in{\mathcal D}_1^\star} C_{1,i}\right)\le  \left(2(N-1)\tilde{C}+2+o(1)\right)\gamma_{i,\eps}^{-1}\ln\frac{1}{r_{1,\eps}}-\gie^{-1}\ln\frac{\gie^2}{\gamma_{1,\eps}^2}\hskip.1cm.$$
This is clearly impossible if $C_{1,i}=0$. Thus we have proved that 
\begin{equation}\label{eq-claim-MT-8-12}
\hbox{ for any }i\in {\mathcal D}_1^\star,\, C_{1,i}>0\hskip.1cm.
\end{equation}
This implies thanks to \eqref{eq-MT-15} that $\rie\ge r_{1,\eps}$ for all $i\in {\mathcal D}_1^\star$. Then we can apply Claim \ref{claim-MT-LBA} to all $i\in {\mathcal D}_1^\star$ thanks to Claim \ref{claim-MT-7} and to what we just said to get that, for any $\tilde{x}_i\in {\mathcal D}_1$, 
\begin{equation}\label{eq-claim-MT-8-12bis}
\gie\left(u_\eps\left(x_{1,\eps}+r_{1,\eps}x\right)-B_{i,\eps}\left(r_{1,\eps}\right)\right) \to 2\ln \frac{1}{\left\vert x-\tilde{x}_i\right\vert} + {\mathcal H}_i\hbox{ in }{\mathcal C}^1_{loc}\left({\mathbb D}_{\tilde{x}_i}\left(1\right)\setminus\left\{\tilde{x}_i\right\}\right)\hbox{ as }\eps\to 0
\end{equation}
where ${\mathcal H}_i$ is some harmonic function in ${\mathbb D}_{\tilde{x}_i}\left(1\right)$ satisfying ${\mathcal H}_i\left(\tilde{x}_i\right)=0$ and $\nabla {\mathcal H}_i\left(\tilde{x}_i\right)=0$ (note here that we assumed that $r_{1,\eps}\to 0$ as $\eps\to 0$). Let us set now 
$$v_\eps\left(x\right)=\gamma_{1,\eps} \left(u_\eps\left(x_{1,\eps}+r_{1,\eps}x\right)-B_{1,\eps}\left(r_{1,\eps}\right)\right)\hskip.1cm.$$
Thanks to Claim \ref{claim-MT-7}, we have that 
$$\left\vert \nabla v_\eps\right\vert \le C_R \hbox{ in }\Omega_0^R$$
for all $R>0$. This clearly proves that $\left(v_\eps\right)$ is uniformly bounded in any $\Omega_0^R$. Since 
$$\Delta v_\eps = \lambda_\eps r_{1,\eps}^2 \gamma_{1,\eps} f_\eps\left(x_{1,\eps}+r_{1,\eps}x\right) u_\eps\left(x_{1,\eps}+r_{1,\eps}x\right)e^{u_\eps\left(x_{1,\eps}+r_{1,\eps}x\right)^2}$$
in $\Omega_0^R$, we have that 
$$\left\vert \Delta v_\eps\right\vert =O\left(\lambda_\eps r_{1,\eps}^2 \gamma_{1,\eps} \left(B_{1,\eps}\left(r_{1,\eps}\right)+\gamma_{1,\eps}^{-1}\right)e^{B_{1,\eps}\left(r_{1,\eps}\right)^2}\right)\hbox{ in }\Omega_0^R\hskip.1cm.  $$
Thanks to \eqref{eq-claim-MT-8-10bis}, we know that 
$$\lambda_\eps r_{1,\eps}^2 \le C\gamma_{1,\eps}^{-2} e^{-\gamma_{1,\eps}B_{1,\eps}\left(r_{1,\eps}\right)}$$
so that 
$$\left\vert \Delta v_\eps\right\vert =O\left(\gamma_{1,\eps}^{-1} \left(B_{1,\eps}\left(r_{1,\eps}\right)+\gamma_{1,\eps}^{-1}\right)e^{B_{1,\eps}\left(r_{1,\eps}\right)^2-\gamma_{1,\eps}B_{1,\eps}\left(r_{1,\eps}\right)}\right)=o(1)\hbox{ in }\Omega_0^R$$
thanks to Claim \ref{claim-MT-1}. Thus we have by standard elliptic theory that 
\begin{equation}\label{eq-claim-MT-8-13}
v_\eps\to v_0\hbox{ in }{\mathcal C}^1_{loc}\left({\mathbb R}^2\setminus {\mathcal S}_1\right)\hbox{ as }\eps\to 0
\end{equation}
where $v_0$ is some harmonic funtion in ${\mathbb R}^2\setminus {\mathcal S}_1$ which satisfies, thanks to Claim \ref{claim-MT-7},
\begin{equation}\label{eq-claim-MT-8-14}
\left\vert \nabla v_0\right\vert \le \frac{C}{\vert x\vert} \hbox{ for }\vert x\vert\hbox{ large.}
\end{equation}
Thanks to \eqref{eq-claim-MT-8-12bis}, we know that 
$$v_0(x)= 2C_{1,i}\ln\frac{1}{\left\vert x-\tilde{x}_i\right\vert} + C_{1,i} {\mathcal H}_i + B_i\hbox{ in }{\mathbb D}_{\tilde{x}_i}\left(1\right)$$
for all $i\in {\mathcal D}_1$ where $B_i$ is a constant given by 
$$B_i = \left(1-C_{1,i}\right)\left(1+\ln \frac{f_0\left(x_1\right)}{4}\right) + 2C_{1,i}\ln C_{1,i} +\lim_{\eps\to 0} \left( \left(1-\frac{\gamma_{1,\eps}}{\gie}\right)\ln \left(\lambda_\eps \gamma_{1,\eps}^2 r_{1,\eps}^2\right)\right)\hskip.1cm.$$ 
Thus we have that 
$$v_0(x)=2\ln\frac{1}{\vert x\vert}+2 \sum_{i\in {\mathcal D}_1^\star} C_{1,i}\ln\frac{1}{\left\vert x-\tilde{x}_i\right\vert}+ w_0$$
where $w_0$ is harmonic in ${\mathbb R}^2$ and satisfies thanks to \eqref{eq-claim-MT-8-14} that $\left\vert \nabla w_0\right\vert \le C\vert x\vert^{-1}$ for $\vert x\vert$ large. This implies that $w_0\equiv A_0$ for some constant $A_0$. Thus we have that 
\begin{equation}\label{eq-claim-MT-8-15}
v_0(x)=2\ln\frac{1}{\vert x\vert}+2 \sum_{i\in {\mathcal D}_1^\star} C_{1,i}\ln\frac{1}{\left\vert x-\tilde{x}_i\right\vert}+ A_0\hskip.1cm.
\end{equation}
Moreover, the ${\mathcal H}_i$'s of \eqref{eq-claim-MT-8-12bis} are given by 
$${\mathcal H}_i(x)=2\ln\frac{1}{\vert x\vert}+2 \sum_{j\in {\mathcal D}_1^\star, \, j\neq i} C_{1,j}\ln\frac{1}{\left\vert x-\tilde{x}_j\right\vert}+ A_0 $$
and they satisfy $\nabla {\mathcal H}_i\left(\tilde{x}_i\right)=0$ for all $i\in {\mathcal D}_1$. Note that, by the definition of $r_{1,\eps}$ and since we assumed that $\frac{r_{1,\eps}}{d_{1,\eps}}\to +\infty$ as $\eps\to 0$, we know that ${\mathcal D}_1^\star\neq\emptyset$. Let us pick up $i\in {\mathcal D}_1^\star$ such that $\left\vert \tilde{x}_i\right\vert \ge \left\vert \tilde{x}_j\right\vert$ for all $j\in {\mathcal D}_1^\star$. It is then clear that 
$$\langle \nabla {\mathcal H}_i\left(\tilde{x}_i\right), \tilde{x}_i\rangle =-2 - 2\sum_{j\in {\mathcal D}_1^\star, \, j\neq i} C_{1,i} \frac{\langle \tilde{x}_i-\tilde{x}_j,\tilde{x}_i\rangle}{\left\vert \tilde{x}_i-\tilde{x}_j\right\vert^2}\le -2\hskip.1cm,$$
which contradicts the fact that $\nabla {\mathcal H}_i\left(\tilde{x}_i\right)=0$. This is the contradiction we were looking for and this proves that, if $r_{1,\eps}\to 0$ as $\eps\to 0$, this first case can not happen, that is we must have $d_{1,\eps}\to 0$ as $\eps\to 0$. \hfill $\spadesuit$

\medskip {\sc Case 2} - We assume that $d_{1,\eps}\to 0$ and that $\frac{r_{1,\eps}}{d_\eps}\to 0$ as $\eps\to 0$. 

\smallskip\noindent We let $y\in \Omega_0^R$ for some $R>0$ and we set $x_\eps=x_{1,\eps}+r_{1,\eps}x$. Since $d_{1,\eps}\to 0$ and $r_{1,\eps}\to 0$, we are in situation b) of Claim \ref{claim-MT-6}. Indeed, as in Case 1, we have that 
$$\frac{\left\vert x_\eps-\xie\right\vert}{\mie}\to +\infty\hbox{ as }\eps\to 0\hbox{ for all }i=1,\dots,N\hskip.1cm.$$
\begincal
u_\eps\left(x_\eps\right)&=&\psi_\eps\left(x_\eps\right)+ \sum_{i=1}^N \frac{4\pi+o(1)}{\gie} {\mathcal G}\left(\xie,x_\eps\right)   \\
&&+O\left( \sum_{i\in {\mathcal A}} \gie^{-1}\left(\frac{\mie}{\left\vert \xie-x_\eps\right\vert} +\gie^{-1} \ln\left(\frac{\sie}{\left\vert \xie-x_\eps\right\vert}+2\right)\right)\right)\\
&&+O\left(\sum_{i\in {\mathcal B}} \frac{d_\eps}{d_\eps+d_{i,\eps}} \left(\gie^{-1}\mie + \gie^{-2}\sie\right)\right)
\fincal
where ${\mathcal A}$ is defined as the set of $i\in \left\{1,\dots,N\right\}$ such that $\left\vert \xie-x_\eps\right\vert \le \sie +o\left(d_\eps\right)$ and ${\mathcal B}$ as its complementary. Noting that $\left\vert \xie-x_\eps\right\vert \ge C\rie$ for all $i\in \left\{1,\dots,N\right\}$, we have that for any $i\in {\mathcal A}$,
$$\gie^{-1}\frac{\mie}{\left\vert \xie-x_\eps\right\vert} +\gie^{-2} \ln\left(\frac{\sie}{\left\vert \xie-x_\eps\right\vert}+2\right)=o\left(\gie^{-1}\right)$$
and, for any $i\in {\mathcal B}$, 
$$\frac{d_\eps}{d_\eps+d_{i,\eps}} \left(\gie^{-1}\mie + \gie^{-2}\sie\right) =o\left(\gie^{-1}\right)\hskip.1cm.$$
Thus we have that 
$$u_\eps\left(x_\eps\right)=\psi_\eps\left(x_\eps\right)+ \sum_{i=1}^N \frac{4\pi+o(1)}{\gie} {\mathcal G}\left(\xie,x_\eps\right) +o\left(\gamma_{1,\eps}^{-1}\right)$$
thanks to \eqref{eq-MT-14}. For $i\in {\mathcal D}_1$, we have that $\left\vert \xie-x_\eps\right\vert =o\left(d_{1,\eps}\right)$ so that, thanks to \eqref{eq-appB-11}, 
$${\mathcal G}\left(\xie,x_\eps\right) = \frac{1}{2\pi}\left(\ln \frac{2d_{1,\eps}}{r_{1,\eps}}\right) +O(1)\hskip.1cm.$$
For any $i\not\in {\mathcal D}_1$, we know that 
$${\mathcal G}\left(\xie,x_\eps\right) = {\mathcal G}\left(\xie,x_{1,\eps}\right)+o(1)$$
thanks to \eqref{eq-appB-11}. Thus we can write that 
$$u_\eps\left(x_\eps\right)=\psi_\eps\left(x_\eps\right)+ \sum_{i\in {\mathcal D}_1} \frac{2+o(1)}{\gie}\left(\ln \frac{2d_{1,\eps}}{r_{1,\eps}}\right) +\sum_{i\not\in {\mathcal D}_1}\gie^{-1}{\mathcal G}\left(\xie,x_{1,\eps}\right)+O\left(\gamma_{1,\eps}^{-1}\right)\hskip.1cm.$$
If $\lambda_0\neq 0$, then we can write thanks to the fact that $\psi_\eps=0$ on $\partial \Omega$ and to \eqref{eq-MT-16} that $\psi_\eps\left(x_\eps\right)=O\left(d_{1,\eps}\right)$. This leads with Claim \ref{claim-MT-2} to $\psi_\eps\left(x_\eps\right)=O\left(\gamma_{1,\eps}^{-1}\right)$. If $\lambda_0=0$, then we can use \eqref{eq-MT-17} to arrive to the same result. Thus we finally get that 
\begin{equation}\label{eq-claim-MT-8-16}
u_\eps\left(x_\eps\right)= \frac{2}{\gamma_{1,\eps}}\left(1+\sum_{i\in {\mathcal D}_1^\star}C_{1,i}\right) \left(\ln \frac{d_{1,\eps}}{r_{1,\eps}}\right)+\sum_{i\not\in {\mathcal D}_1}\gie^{-1}{\mathcal G}\left(\xie,x_{1,\eps}\right)+o\left(\gamma_{1,\eps}^{-1}\left(\ln \frac{d_{1,\eps}}{r_{1,\eps}}\right)\right)\hskip.1cm.
\end{equation}
Note that $C_{1,i}\le \tilde{C}$ for all $i>1$ thanks to \eqref{eq-MT-13}. Thus we have in particular that 
\begin{equation}\label{eq-claim-MT-8-17}
2\gamma_{1,\eps}^{-1}\ln\frac{d_{1,\eps}}{r_{1,\eps}}\left(1+\sum_{i\in{\mathcal D}_1^\star} C_{1,i}\right)\le u_\eps\left(x_\eps\right) +o\left(\gamma_{1,\eps}^{-1}\ln\frac{d_{1,\eps}}{r_{1,\eps}}\right)\le 
 2\gamma_{1,\eps}^{-1}\ln\frac{d_{1,\eps}}{r_{1,\eps}}\left(1+(N-1)\tilde{C}\right) \hskip.1cm.
\end{equation}
Here we used \eqref{eq-appB-11} to estimate ${\mathcal G}\left(\xie,x_{1,\eps}\right)$ for $i\not\in {\mathcal D}_1$. Note that we also have thanks to Claim \ref{claim-MT-7} and to \eqref{eq-MT-17} that 
\begin{equation}\label{eq-claim-MT-8-17b}
\left\vert \nabla u_\eps(x)\right\vert \le C\gamma_{1,\eps}^{-1}\left\vert x_{1,\eps}-x\right\vert \hbox{ for all }x\in {\mathbb D}_{x_{1,\eps}}\left(r_{1,\eps}\right)\hskip.1cm.
\end{equation}
The proof now follows exactly Case 1, from \eqref{eq-claim-MT-8-9} to the end. We will not repeat it here. \hfill $\spadesuit$

\medskip {\sc Case 3} - We assume that $d_{1,\eps}\to 0$ as $\eps\to 0$ and that ${\ds \frac{d_{1,\eps}}{r_{1,\eps}}\to L}$ as $\eps\to 0$ where $L\ge 2$. 

\smallskip\noindent We are thus in the case where, after some harmless rotation,
$$\Omega_0 = {\mathbb R}\times\left(-\infty,L\right)\hskip.1cm.$$
We let $y\in \Omega_0^R$ for some $R>0$ and we set $x_\eps=x_{1,\eps}+r_{1,\eps}y$. Since $d_{1,\eps}\to 0$ and $r_{1,\eps}\to 0$, we are in situation b) of Claim \ref{claim-MT-6}. Indeed, as in Case 1, we have that 
$$\frac{\left\vert x_\eps-\xie\right\vert}{\mie}\to +\infty\hbox{ as }\eps\to 0\hbox{ for all }i=1,\dots,N\hskip.1cm.$$
Thus we can write that 
\begincal
u_\eps\left(x_\eps\right)&=&\psi_\eps\left(x_\eps\right)+ \sum_{i=1}^N \frac{4\pi+o(1)}{\gie} {\mathcal G}\left(\xie,x_\eps\right)   \\
&&+O\left( \sum_{i\in {\mathcal A}} \gie^{-1}\left(\frac{\mie}{\left\vert \xie-x_\eps\right\vert} +\gie^{-1} \ln\left(\frac{\sie}{\left\vert \xie-x_\eps\right\vert}+2\right)\right)\right)\\
&&+O\left(\sum_{i\in {\mathcal B}} \frac{d_\eps}{d_\eps+d_{i,\eps}} \left(\gie^{-1}\mie + \gie^{-2}\sie\right)\right)
\fincal
where ${\mathcal A}$ is defined as the set of $i\in \left\{1,\dots,N\right\}$ such that $\left\vert \xie-x_\eps\right\vert \le \sie +o\left(d_\eps\right)$ and ${\mathcal B}$ as its complementary. As in Case 2, we have that 
$$\gie^{-1}\frac{\mie}{\left\vert \xie-x_\eps\right\vert} +\gie^{-2} \ln\left(\frac{\sie}{\left\vert \xie-x_\eps\right\vert}+2\right)=o\left(\gie^{-1}\right)$$
for all $i\in {\mathcal A}$ while 
$$\frac{d_\eps}{d_\eps+d_{i,\eps}} \left(\gie^{-1}\mie + \gie^{-2}\sie\right) = o\left(\gie^{-1}\right)$$
for all $i\in {\mathcal B}$. Thus we have that 
$$u_\eps\left(x_\eps\right)=\psi_\eps\left(x_\eps\right)+ \sum_{i=1}^N \frac{4\pi+o(1)}{\gie} {\mathcal G}\left(\xie,x_\eps\right)+o\left(\gamma_{1,\eps}^{-1}\right)\hskip.1cm.$$
If $i\in {\mathcal D}_1$, we have that 
$${\mathcal G}\left(\xie,x_\eps\right) = \frac{1}{2\pi} \ln \frac{\left\vert \tilde{y}_i-y\right\vert}{\left\vert \tilde{x}_i-y\right\vert} +o(1)$$
where 
$$\tilde{y}_i = {\mathcal R}\left(\tilde{x}_i\right)\hskip.1cm,$$
${\mathcal R}$ being the reflection with respect to the straight line ${\mathbb R}\times \left\{L\right\}$. Here we used \eqref{eq-appB-11}. If $i\not\in {\mathcal D}_1$, we have that 
$${\mathcal G}\left(\xie,x_\eps\right) = o(1)$$
thanks to \eqref{eq-appB-11}. Thus we can write, remembering \eqref{eq-MT-14}, that 
$$u_\eps\left(x_\eps\right)=\psi_\eps\left(x_\eps\right)+ 2\gamma_{1,\eps}^{-1}\sum_{i\in {\mathcal D}_1} C_{1,i} \ln \frac{\left\vert \tilde{y}_i-y\right\vert}{\left\vert \tilde{x}_i-y\right\vert} +o\left(\gamma_{1,\eps}^{-1}\right)\hskip.1cm.$$
Using \eqref{eq-MT-16}, we know that 
$$\frac{\psi_\eps\left(x_\eps\right)}{r_{1,\eps}} \to A\left(L-y_2\right)$$
where $y=\left(y_1,y_2\right)$ for some $A$ independent of $y$. Moreover, we have that $A\ge 0$ by the maximum principle since $\Delta\psi_\eps\ge 0$ in $\Omega$ and $\psi_\eps=0$ on $\partial \Omega$. If $\lambda_0\neq 0$, we can use Claim \ref{claim-MT-2} to deduce that 
$$\gamma_{1,\eps}\psi_\eps\left(x_\eps\right) \to B\left(L-y_2\right)$$
for some $B>0$, independent of $y$. If $\lambda_0=0$, then \eqref{eq-MT-17} implies that 
$$\gamma_{1,\eps}\psi_{1,\eps}\left(x_\eps\right) =O\left(\lambda_\eps r_{1,\eps}\right)=o(1)\hskip.1cm.$$
Thus, up to change the $B$ above, we can write that 
\begin{equation}\label{eq-claim-MT-8-26}
\gamma_{1,\eps}u_\eps\left(x_\eps\right) \to B\left(L-y_2\right) + 2\sum_{i\in {\mathcal D}_1} C_{1,i} \ln \frac{\left\vert \tilde{y}_i-y\right\vert}{\left\vert \tilde{x}_i-y\right\vert}\hbox{ as }\eps\to 0\hskip.1cm.
\end{equation}
Then, by the equation satisfied by $u_\eps$, it is clear that 
$$v_\eps(x)=\gamma_{1,\eps}u_\eps\left(x_{1,\eps}+r_{1,\eps} x\right)$$ 
has a Laplacian uniformly converging to $0$ in any $\Omega_0^R$. Thus, by standard elliptic theory, we can conclude that 
\begin{equation}\label{eq-claim-MT-8-27}
\gamma_{1,\eps}u_\eps\left(x_{1,\eps}+r_{1,\eps} y\right) \to B\left(L-y_2\right) + 2\sum_{i\in {\mathcal D}_1} C_{1,i} \ln \frac{\left\vert \tilde{y}_i-y\right\vert}{\left\vert \tilde{x}_i-y\right\vert}\hbox{ in }{\mathcal C}^1_{loc}\left(\Omega_0\setminus {\mathcal S}_1\right)\hbox{ as }\eps\to 0\hskip.1cm.
\end{equation}
Writing that 
$$\left\vert \nabla \psi_\eps\right\vert \le C\lambda_\eps \hbox{ in }{\mathbb D}_{x_{1,\eps}}\left(r_{1,\eps}\right)$$
thanks to \eqref{eq-MT-16}, we get with Claim \ref{claim-MT-2} that 
$$\left\vert \nabla \psi_\eps\right\vert \le C\sqrt{\lambda_\eps} d_{1,\eps}^{-1}\gamma_{1,\eps}^{-1} 
\hbox{ in }{\mathbb D}_{x_{1,\eps}}\left(r_{1,\eps}\right)$$
so that we can use Claim \ref{claim-MT-7} and \eqref{eq-MT-14} to obtain that 
$$\left\vert \nabla u_\eps\right\vert \le C\gamma_{1,\eps}^{-1}\left\vert x_{1,\eps}-x_\eps\right\vert^{-1}\hbox{ in }{\mathbb D}_{x_{1,\eps}}\left(r_{1,\eps}\right)\hskip.1cm.$$
We are thus in position to apply Claim \ref{claim-MT-LBA} to $i=1$. In particular, combined with \eqref{eq-claim-MT-8-27}, we get that 
$$\gamma_{1,\eps}B_{1,\eps}\left(r_{1,\eps}\right)=O\left(1\right)\hskip.1cm.$$
This leads with Claim \ref{claim-appA-4} of Appendix A to 
$$\ln \left(\lambda_\eps r_{1,\eps}^2\gamma_{1,\eps}^2\right) = O\left(1\right)\hskip.1cm.$$
Thus we have that, up to a subsequence, 
\begin{equation}\label{eq-claim-MT-8-28}
\lambda_\eps f_\eps\left(x_{1,\eps}\right)r_{1,\eps}^2\gamma_{1,\eps}^2\to \alpha_0\hbox{ as }\eps\to 0
\end{equation}
for some $\alpha_0>0$. Let now $i\in {\mathcal D}_1$ be such that the second coordinate of $\tilde{x}_i$ satisfies $\left(\tilde{x}_i\right)_2<L$ and 
$$\left(\tilde{x}_i\right)_2\ge \left(\tilde{x}_j\right)_2\hbox{ or }\left(\tilde{x}_j\right)_2=L\hbox{ for all }j\in{\mathcal D}_1\hskip.1cm.$$
Note that such a $i$ does exists since $1\in {\mathcal D}_1$. Moreover, we have that 
$$L>\left(\tilde{x}_i\right)_2\ge \left(\tilde{x}_1\right)_2=0\hskip.1cm.$$
Note also that $d_{i,\eps}\ge \left(L-\left(\tilde{x}_i\right)_2+o(1)\right)r_{1,\eps}$ so that Claim \ref{claim-MT-2} implies that 
$$O(1)= \lambda_\eps f_\eps\left(x_{1,\eps}\right)\gie^{2} d_{i,\eps}^2= \frac{\gie^2}{\gamma_{1,\eps}^2}\frac{d_{i,\eps}^2}{r_{1,\eps}^2} \lambda_\eps f_\eps\left(x_{1,\eps}\right)\gamma_{1,\eps}^2 r_{1,\eps}^2 \ge \left(\left(L-\left(\tilde{x}_i\right)_2\right)^2\alpha_0 +o(1)\right)\frac{\gie^2}{\gamma_{1,\eps}^2}$$
thanks to \eqref{eq-claim-MT-8-28}. This implies that $C_{1,i}\neq 0$. Thanks to \eqref{eq-MT-15}, we then have that $r_{i,\eps}\ge r_{1,\eps}$. Once again, thanks to Claim \ref{claim-MT-7} and to \eqref{eq-MT-14}, we see now that 
$$\left\vert \nabla u_\eps\right\vert \le C\gamma_{i,\eps}^{-1}\left\vert x_{i,\eps}-x_\eps\right\vert^{-1}\hbox{ in }{\mathbb D}_{x_{i,\eps}}\left(r_{1,\eps}\right)$$
and that we can apply Claim \ref{claim-MT-LBA}. In particular, using \eqref{eq-claim-MT-8-28}, we get that 
$$\gamma_{i,\eps} u_\eps\left(x_{1,\eps}+r_{1,\eps} x\right)\to 2\ln \frac{1}{\left\vert x-\tilde{x}_i\right\vert} + {\mathcal H}_i-\ln\left(\frac{\alpha_0}{4C_{1,i}^2}\right)\hbox{ in }{\mathcal C}^1_{loc}\left({\mathbb D}_{\tilde{x}_i}\left(1\right)\setminus \left\{\tilde{x}_i\right\}\right)\hbox{ as }\eps\to 0$$ 
where ${\mathcal H}_i$ is harmonic in ${\mathbb D}_{\tilde{x}_i}\left(1\right)$ and satisfies $\nabla {\mathcal H}_i\left(\tilde{x}_i\right)=0$ (since $r_{1,\eps}\to 0$ as $\eps\to 0$ by assumption). Now, combining this with \eqref{eq-claim-MT-8-27}, we know that 
$${\mathcal H}_i = \frac{B}{C_{1,i}}\left(L-x_2\right) + 2\sum_{j\in {\mathcal D}_1, j\neq i} \frac{C_{1,j}}{C_{1,i}} \ln \frac{\left\vert \tilde{y}_j-x\right\vert}{\left\vert \tilde{x}_j-x\right\vert}+2\ln \left\vert \tilde{y}_i-x\right\vert +\ln\left(\frac{\alpha_0}{4C_{1,i}^2}\right)\hskip.1cm.$$
The derivative of ${\mathcal H}_i$ with respect to the second coordinate at $\tilde{x}_i$ is 
$$\frac{\partial {\mathcal H}_i}{\partial x_2}\left(\tilde{x}_i\right) = -\frac{B}{C_{1,i}} +2\sum_{j\in {\mathcal D}_1, j\neq i} \frac{C_{1,j}}{C_{1,i}} \left(\frac{\left(\tilde{x}_i\right)_2 - \left(\tilde{y}_j\right)_2}{\left\vert \tilde{y}_j-\tilde{x}_i\right\vert^2} - \frac{\left(\tilde{x}_i\right)_2 - \left(\tilde{x}_j\right)_2}{\left\vert \tilde{x}_j-\tilde{x}_i\right\vert^2}\right) +2\frac{\left(\tilde{x}_i\right)_2 - \left(\tilde{y}_i\right)_2}{\left\vert \tilde{y}_i-\tilde{x}_i\right\vert^2}\hskip.1cm.$$
Note now that 
$$\left(\tilde{y}_j\right)_2=2L-\left(\tilde{x}_j\right)_2$$
so that 
$$\frac{\partial {\mathcal H}_i}{\partial x_2}\left(\tilde{x}_i\right) = -\frac{B}{C_{1,i}} +2\sum_{j\in {\mathcal D}_1, j\neq i} \frac{C_{1,j}}{C_{1,i}} \left(\frac{\left(\tilde{x}_i\right)_2+\left(\tilde{x}_j\right)_2 - 2L}{\left\vert \tilde{y}_j-\tilde{x}_i\right\vert^2} - \frac{\left(\tilde{x}_i\right)_2 - \left(\tilde{x}_j\right)_2}{\left\vert \tilde{x}_j-\tilde{x}_i\right\vert^2}\right) +4\frac{\left(\tilde{x}_i\right)_2 - L}{\left\vert \tilde{y}_i-\tilde{x}_i\right\vert^2}\hskip.1cm.$$
We claim that
\begin{equation}\label{eq-claim-MT-8-29}
\frac{\left(\tilde{x}_i\right)_2+\left(\tilde{x}_j\right)_2 - 2L}{\left\vert \tilde{y}_j-\tilde{x}_i\right\vert^2}\le \frac{\left(\tilde{x}_i\right)_2 - \left(\tilde{x}_j\right)_2}{\left\vert \tilde{x}_j-\tilde{x}_i\right\vert^2}
\end{equation}
for all $j\in {\mathcal D}_1$ with $j\neq i$. This will imply that 
$$\frac{\partial {\mathcal H}_i}{\partial x_2}\left(\tilde{x}_i\right) <0\hskip.1cm,$$
all the terms above being nonpositive, the last one being negative. This will give a contradiction with the fact that $\nabla {\mathcal H}_i\left(\tilde{x}_i\right)=0$, thus proving that this last case is not possible either. In order to prove \eqref{eq-claim-MT-8-29}, we first notice that 
$$\left\vert \tilde{x}_i-\tilde{y}_j\right\vert^2 = \left\vert \tilde{x}_i-\tilde{x}_j\right\vert^2 + 4\left(L-\left(\tilde{x}_i\right)_2\right) \left(L-\left(\tilde{x}_j\right)_2\right)\hskip.1cm.$$
Thus we can write that 
\begincal
&&\left(\left(\tilde{x}_i\right)_2+\left(\tilde{x}_j\right)_2 - 2L\right) \left\vert \tilde{x}_i-\tilde{x}_j\right\vert^2 - \left(\left(\tilde{x}_i\right)_2 - \left(\tilde{x}_j\right)_2\right) \left\vert \tilde{x}_i-\tilde{y}_j\right\vert^2 \\
&&\quad =2\left\vert \tilde{x}_i-\tilde{x}_j\right\vert^2 \left(\left(\tilde{x}_j\right)_2-L\right) -4 \left(L-\left(\tilde{x}_i\right)_2\right) \left(L-\left(\tilde{x}_j\right)_2\right)\left(\left(\tilde{x}_i\right)_2 - \left(\tilde{x}_j\right)_2\right)\\
&&\quad = 2\left(\left(\tilde{x}_j\right)_2-L\right) \left(\left\vert \tilde{x}_i-\tilde{x}_j\right\vert^2+ 2  \left(L-\left(\tilde{x}_i\right)_2\right)\left(\left(\tilde{x}_i\right)_2 - \left(\tilde{x}_j\right)_2\right)\right)\\
&&\quad \le 0
\fincal
since $\left(\tilde{x}_j\right)_2-L\le 0$ and if $\left(\tilde{x}_j\right)_2-L\neq 0$, $\left(\tilde{x}_i\right)_2 - \left(\tilde{x}_j\right)_2\ge 0$. This clearly proves \eqref{eq-claim-MT-8-29} and, as already said, proves that this last case is not possible. \hfill $\spadesuit$

\medskip The study of these three cases proves that the assumption $r_{1,\eps}\to 0$ is absurd and thus proves the claim. \hfill $\diamondsuit$

\medskip\noindent Note that this claim implies that 
\begin{equation}\label{eq-MT-19}
x_{1,\eps}\to x_1\hbox{ as }\eps\to 0 \hbox{ with }x_1\in \Omega\hskip.1cm.
\end{equation}
We also have thanks to Claims \ref{claim-MT-2} and \ref{claim-MT-8} that 
$$\lambda_\eps = O\left(\gamma_{1,\eps}^{-2}\right)$$
so that $\lambda_0= 0$ and $u_0\equiv 0$. Moreover, we can transform \eqref{eq-MT-17} into 
\begin{equation}\label{eq-MT-20}
\left\Vert \nabla \psi_\eps\right\Vert_{{\mathcal C}^{1,\alpha}\left(\overline{\Omega}\right)}=O\left(\gamma_{1,\eps}^{-3}\right)\hskip.1cm.
\end{equation}
Let us now give a simple consequence of the previous claim~:

\begin{claim}\label{claim-MT-8bis}
After passing to a subsequence,
$$\lambda_\eps \gamma_{1,\eps}^2\to \alpha_0\hbox{ as }\eps\to 0$$
for some 
$$0<\alpha_0 \le \frac{4}{f_0\left(x_1\right) d\left(x_1,\partial\Omega\right)^2}\hskip.1cm.$$
\end{claim}

\medskip {\it Proof of Claim \ref{claim-MT-8bis}} - We already said that $\lambda_\eps=O\left(\gamma_{1,\eps}^{-2}\right)$. Claim \ref{claim-MT-7} with \eqref{eq-MT-20} gives that 
$$\left\vert \nabla u_\eps\right\vert \le C \sum_{i=1}^N \gie^{-1} \left(\mie+\left\vert x-\xie\right\vert\right)^{-1}\hbox{ in }\Omega\hskip.1cm.$$
This gives in particular that 
$$\left\vert \nabla u_\eps\right\vert \le C\gamma_{1,\eps}^{-1}\left\vert x-x_{1,\eps}\right\vert^{-1}$$
in ${\mathbb D}_{x_{1,\eps}}\left(\delta_0\right)$ where $\delta_0$ is as in Claim \ref{claim-MT-8}. Thus we are in position to apply Claim \ref{claim-MT-LBA} to $i=1$. This gives in particular that 
$$\gamma_{1,\eps}\left(u_\eps(x)-B_{1,\eps}\left(\delta_0\right)\right)=O(1)$$
for all $\vert x-x_{1,\eps}\vert =\frac{\delta_0}{2}$. Now Claim \ref{claim-MT-6} combined with \eqref{eq-MT-20} gives that 
$$\gamma_{1,\eps}u_\eps\left(x\right)=O(1)\hbox{ on }\partial {\mathbb D}_{x_{1,\eps}}\left(\frac{\delta_0}{2}\right)$$ 
so that the above leads to 
$$\gamma_{1,\eps}B_{1,\eps}\left(\delta_0\right)= O(1)\hskip.1cm.$$
Since 
$$\gamma_{1,\eps} B_{1,\eps}\left(\delta_0\right) = -\ln\left(\lambda_\eps\gamma_{1,\eps}^2\right)+O(1)\hskip.1cm,$$
we obtain that 
$$\ln \left(\lambda_\eps\gamma_{1,\eps}^2\right)=O(1)\hskip.1cm.$$
This clearly permits to prove the claim. \hfill $\diamondsuit$

\begin{claim}\label{claim-MT-9}
We have that $\rie\ge \delta_1$ for some $\delta_1>0$ for all $i=1,\dots,N$.
\end{claim}

\medskip {\it Proof of Claim \ref{claim-MT-9}} - We shall prove it by induction on $i$. This is already proved for $i=1$ in the previous claim. Fix $2\le i\le N$ and assume that 
\begin{equation}\label{eq-claim-MT-9-1}
\rje \ge \delta_1>0 \hbox{ for all }j<i\hskip.1cm.
\end{equation}
In particular, after passing to a subsequence, we have that 
\begin{equation}\label{eq-claim-MT-9-1bis}
\xje\to x_j\hbox{ as }\eps\to 0\hbox{ with }x_j\in \Omega\hskip.1cm.
\end{equation}
Assume by contradiction that 
\begin{equation}\label{eq-claim-MT-9-2}
\rie\to 0\hbox{ as }\eps\to 0\hskip.1cm.
\end{equation}
By \eqref{eq-MT-15}, this implies that $C_{j,i}=0$ for all $j<i$ so that 
\begin{equation}\label{eq-claim-MT-9-3}
\gamma_{j,\eps} = o\left(\gie\right)\hbox{ for all }j<i\hskip.1cm.
\end{equation}
We shall now proceed as in the proof of Claim \ref{claim-MT-8} and distinguish three cases. 

\smallskip We let in the following 
\begin{equation}\label{eq-claim-MT-9-4}
{\mathcal D}_i^\star = \left\{j>i\hbox{ s.t. }\left\vert \xie-\xje\right\vert = O\left(\rie\right)\right\}\hbox{ and }{\mathcal D}_i={\mathcal D}_i^\star \cup \left\{i\right\}\hskip.1cm.
\end{equation}
After passing to a subsequence, we let 
\begin{equation}\label{eq-claim-MT-9-5}
{\mathcal S}_i^\star = \left\{\tilde{x}_j=\lim_{\eps\to 0}\frac{\xje-\xie}{\rie},\, j\in {\mathcal D}_i^\star\right\}\hbox{ and }{\mathcal S}_i={\mathcal S}_i^\star\cup\left\{\tilde{x}_i=0\right\}\hskip.1cm.
\end{equation}
We also let 
\begin{equation}\label{eq-claim-MT-9-6}
\Omega_{i,\eps} = \left\{y\in {\mathbb R}^2\hbox{ s.t. }\xie+\rie y \in \Omega\right\}\hskip.1cm.
\end{equation}
Note that, after passing to a subsequence (and up to a harmless rotation if necessary), we have that 
\begin{equation}\label{eq-claim-MT-9-7}
\Omega_{i,\eps} \to \Omega_0\hbox{ as }\eps\to 0\hbox{ where }\left\{
\begin{array}{ll}
{\ds \Omega_0={\mathbb R}^2}&{\ds \hbox{ if }\frac{\die}{\rie}\to +\infty\hbox{ as }\eps\to 0}\\
\,&\,\\
{\ds\Omega_0={\mathbb R}\times\left(-\infty,L\right)}&{\ds \hbox{ if } \frac{\die}{\rie}\to L\hbox{ as }\eps\to 0}
\end{array}\right.
\end{equation}
Here $d_{i,\eps}= d\left(x_{i,\eps},\partial\Omega\right)$, as defined in \eqref{eq-MT-8}. For $R>0$, we shall also let 
\begin{equation}\label{eq-claim-MT-9-8}
\Omega_0^R = \left(\Omega_0\cap {\mathbb D}_0(R)\right)\setminus \bigcup_{j\in {\mathcal D}_i} {\mathbb D}_{\tilde{x}_j}\left(\frac{1}{R}\right)\hskip.1cm.
\end{equation}

\smallskip {\sc Case 1} - We assume that $\die\not \to 0$ as $\eps\to 0$, meaning that, after passing to a subsequence, $\xie\to x_i$ as $\eps\to 0$ with $x_i\in \Omega$. 

\smallskip\noindent We let $y\in \Omega_0^R$ for some $R>0$ and we set $x_\eps=\xie+\rie y$. Since $\die\not\to 0$ and $\rie\to 0$ as $\eps\to 0$, we are in situation a) of Claim \ref{claim-MT-6}. Note indeed that 
$$\frac{\left\vert x_\eps-\xje\right\vert}{\mje}\to +\infty\hbox{ as }\eps\to 0\hbox{ for all }j=1,\dots,N\hskip.1cm.$$
It is obvious if $j<i$ since $r_{j,\eps}\ge \delta_1>0$ and $\rie\to 0$ as $\eps\to 0$. It is also obvious if $j\in{\mathcal D}_i$ since we clearly have in this case 
$$\frac{\left\vert x_\eps-\xje\right\vert}{\mje} = \frac{\left\vert x_\eps-\xje\right\vert}{\rie}\frac{\rie}{\rje}\frac{\rje}{\mje}$$
with ${\ds \frac{\left\vert x_\eps-\xje\right\vert}{\rie}\ge R^{-1}+o(1)}$, ${\ds \frac{\rie}{\rje}\ge 2\left\vert \tilde{x}_j\right\vert^{-1}+o(1)}$ for $j\in {\mathcal D}_i^\star$ and equal to $1$ if $j=i$, and ${\ds \frac{\rje}{\mje}\to +\infty}$ as $\eps\to 0$ thanks to assertion c) of Proposition \ref{prop-DMJ}. While, if $j>i$ and $j\not\in{\mathcal D}_i$, we can write that 
$$\frac{\left\vert x_\eps-\xje\right\vert}{\mie} \ge \left(1+o(1)\right) \frac{\left\vert \xie-\xje\right\vert}{\mje}\ge  \left(2+o(1)\right) \frac{\rje}{\mje}\to +\infty\hbox{ as }\eps\to 0\hskip.1cm.$$
Thus, applying a) of Claim \ref{claim-MT-6}, we can write that 
\begincal
u_\eps\left(x_\eps\right)&=&\psi_\eps\left(x_\eps\right)+\sum_{j=1}^N \left(4\pi+o(1)\right)\gje^{-1}{\mathcal G}\left(\xje,x_\eps\right)\\
&&+O\left(\sum_{i=1}^N \left(\gje^{-1}\frac{\mje}{\left\vert \xje-x_\eps\right\vert} +\gje^{-2} \ln\left(\frac{\sje}{\left\vert \xje-x_\eps\right\vert}+2\right)\right)\right)\hskip.1cm.
\fincal
For $j<i$, we have that 
$$\left(4\pi+o(1)\right)\gje^{-1}{\mathcal G}\left(\xje,x_\eps\right) = 4\pi \gamma_{1,\eps}^{-1} C_{1,j} {\mathcal G}\left(x_j,x_i\right)+o\left(\gamma_{1,\eps}^{-1}\right)$$
thanks to the assumption that $\xie\to x_i$ with $x_i\in \Omega$, to \eqref{eq-claim-MT-9-1bis} and to \eqref{eq-MT-12}, \eqref{eq-MT-14}. We also obviously have that 
$$\gje^{-1}\frac{\mje}{\left\vert \xje-x_\eps\right\vert} +\gje^{-2} \ln\left(\frac{\sje}{\left\vert \xje-x_\eps\right\vert}+2\right)= o\left(\gamma_{1,\eps}^{-1}\right)\hskip.1cm.$$
We also know thanks to \eqref{eq-MT-20} that $\psi_\eps\left(x_\eps\right)=o\left(\gamma_{1,\eps}^{-1}\right)$. For $j>i$, we can proceed exactly as in Case 1 of Claim \ref{claim-MT-8} to finally obtain that 
\begin{eqnarray}\label{eq-claim-MT-9-9}
u_\eps\left(x_\eps\right)&=& 4\pi \gamma_{1,\eps}^{-1} \sum_{j=1}^{i-1} C_{1,j} {\mathcal G}\left(x_j,x_i\right)+o\left(\gamma_{1,\eps}^{-1}\right)\nonumber\\
&&+\left(2+o(1)\right)\gie^{-1}\ln\frac{1}{\rie}\left(1+\sum_{j\in{\mathcal D}_i^\star} C_{i,j}\right)\\
&&+ \sum_{j>i,\, j\not\in {\mathcal D}_i} \left(4\pi+o(1)\right)\gje^{-1}{\mathcal G}\left(\xie,\xje\right)+O\left(\gie^{-1}\right)\nonumber
\end{eqnarray}
This gives in particular that $u_\eps\ge C\gamma_{1,\eps}^{-1}$ on $\partial{\mathbb D}_{\xie}\left(\rie\right)$ for some $C>0$. Using b) of Claim \ref{claim-MT-0}, we deduce that 
$$C\gamma_{1,\eps}^{-1}\le B_{i,\eps}\left(\rie\right)+O\left(\gie^{-1}\right)= -\gie^{-1} \ln \left(\lambda_\eps \gie^2\rie^2\right)+O\left(\gie^{-1}\right)\hskip.1cm.$$
Since $\gamma_{1,\eps}=o\left(\gie\right)$, see \eqref{eq-claim-MT-9-3}, we deduce that 
$$\lambda_\eps\gie^2\rie^2\to 0\hbox{ as }\eps\to 0\hskip.1cm.$$
Thanks to Claim \ref{claim-MT-8bis}, this gives that 
\begin{equation}\label{eq-claim-MT-9-10}
\gie\rie = o\left(\gamma_{1,\eps}\right)\hskip.1cm.
\end{equation}
We apply now Claim \ref{claim-MT-7} combined with \eqref{eq-MT-20} and this last estimate to write that 
$$\left\vert \nabla u_\eps\left(x\right)\right\vert \le C \gie^{-1}\left\vert \xie-x\right\vert^{-1} + C\gamma_{1,\eps}^{-1}\le C'  \gie^{-1}\left\vert \xie-x\right\vert^{-1}$$
in ${\mathbb D}_{\xie}\left(\rie\right)$. Thus we can apply Claim \ref{claim-MT-LBA} to $i$~: this gives that, if $\left\vert x\right\vert=\frac{1}{2}$,
$$u_\eps\left(\xie+\rie x\right)= B_{i,\eps}\left(\rie\right)+O\left(\gie^{-1}\right)\hskip.1cm.$$
Combined with \eqref{eq-claim-MT-9-9}, this leads to 
\begincal
B_{i,\eps}\left(\rie\right)  &=& 4\pi \gamma_{1,\eps}^{-1} \sum_{j=1}^{i-1} C_{1,j} {\mathcal G}\left(x_j,x_i\right)+o\left(\gamma_{1,\eps}^{-1}\right)\\
&&+\left(2+o(1)\right)\gie^{-1}\ln\frac{1}{\rie}\left(1+\sum_{j\in{\mathcal D}_i^\star} C_{i,j}\right)\\
&&+ \sum_{j>i,\, j\not\in {\mathcal D}_i} \left(4\pi+o(1)\right)\gje^{-1}{\mathcal G}\left(\xie,\xje\right)+O\left(\gie^{-1}\right)\hskip.1cm.
\fincal
Since 
$$B_{i,\eps}\left(\rie\right)= - \gie^{-1}\ln \left(\lambda_\eps \gie^2 \rie^2\right)+ O\left(\gie^{-1}\right)\hskip.1cm,$$
this leads to 
\begin{eqnarray}\label{eq-claim-MT-9-11}
- \ln \left(\lambda_\eps \gie^2 \rie^2\right)&=&\gie \gamma_{1,\eps}^{-1}\left(4\pi \sum_{j=1}^{i-1} C_{1,j} {\mathcal G}\left(x_j,x_i\right)+o(1)\right)\nonumber\\
&&+2\ln\frac{1}{\rie}\left(1+\sum_{j\in{\mathcal D}_i^\star} C_{i,j}\right)\\
&&+ \sum_{j>i,\, j\not\in {\mathcal D}_i} \left(4\pi+o(1)\right)C_{i,j}{\mathcal G}\left(\xie,\xje\right)+o\left(\ln\frac{1}{\rie}\right)\hskip.1cm.\nonumber
\end{eqnarray}
Thanks to Claim \ref{claim-MT-8bis} and \eqref{eq-claim-MT-9-3}, we deduce that 
\begin{equation}\label{eq-claim-MT-9-12}
4\pi \sum_{j=1}^{i-1} C_{1,j} {\mathcal G}\left(x_j,x_i\right)+o(1)+\frac{\gamma_{1,\eps}}{\gie}\ln\frac{1}{\rie}\left(2\sum_{j\in{\mathcal D}_1^\star} C_{i,j}+o(1)\right) \le 0\hskip.1cm.
\end{equation}
Let $k\in {\mathcal D}_i^\star$. It is clear that there exists $\delta>0$ such that $\partial {\mathbb D}_{\xke}\left(\delta \rie\right)\subset \left\{\xie+\rie y,\, y\in 
 \Omega_0^R\right\}$ for some $R>0$. Thus we can write that 
$$\inf_{\partial {\mathbb D}_{\xke}\left(\delta \rie\right)} u_\eps\ge C\gamma_{1,\eps}^{-1}+2\gie^{-1}\ln\frac{1}{\rie}\left(1+\sum_{j\in{\mathcal D}_i^\star} C_{i,j}+o(1)\right) $$
thanks to \eqref{eq-claim-MT-9-9}. We can also apply b) of Claim \ref{claim-MT-0} with $r_\eps=\delta \rie$ thanks to Claim \ref{claim-MT-1} and to the fact that ${\ds \frac{\rie}{\rke}\ge 2\left\vert \tilde{x}_k\right\vert^{-1}}+o(1)$. This leads to 
\begincal
&&C\gamma_{1,\eps}^{-1}+2\gie^{-1}\ln\frac{1}{\rie}\left(1+\sum_{j\in{\mathcal D}_i^\star} C_{i,j}+o(1)\right)\\
&&\quad  \le B_{k,\eps}\left(\delta \rie\right)+o\left(\gke^{-1}\right)=-\gke^{-1} \ln\left(\lambda_\eps \rie^2\gke^2\right) + O\left(\gke^{-1}\right) \hskip.1cm.
\fincal
Combined with \eqref{eq-claim-MT-9-11}, this gives that 
\begincal
C\gamma_{1,\eps}^{-1}+2\gie^{-1}\ln\frac{1}{\rie}\left(1+\sum_{j\in{\mathcal D}_i^\star} C_{i,j}+o(1)\right) &\le& \gke^{-1}\gie \gamma_{1,\eps}^{-1}\left(4\pi \sum_{j=1}^{i-1} C_{1,j} {\mathcal G}\left(x_j,x_i\right)+o(1)\right)\\
&&+O\left(\gke^{-1}\ln\frac{1}{\rie}\right) - \gke^{-1} \ln\left(\frac{\gke^2}{\gie^2}\right)\hskip.1cm.
\fincal
Assume by contradiction that $\gie=o\left(\gke\right)$. We then have that 
$$C\gamma_{1,\eps}^{-1}+2\gie^{-1}\ln\frac{1}{\rie}\left(1+\sum_{j\in{\mathcal D}_i^\star} C_{i,j}+o(1)\right) \le  o\left(\gamma_{1,\eps}^{-1}\right) + o\left(\gie^{-1}\ln\frac{1}{\rie}\right)\hskip.1cm,$$
which is absurd. Thus we have proved that 
\begin{equation}\label{eq-claim-MT-9-13} 
C_{i,j} >0 \hbox{ for all }j\in {\mathcal D}_i^\star\hskip.1cm.
\end{equation}
Since $\frac{\rie}{\die}\to +\infty$ as $\eps\to 0$ and since $\frac{\left\vert \xie-\xje\right\vert}{\rie}\to +\infty$ for all $j<i$, we are sure that ${\mathcal D}_i^\star\neq\emptyset$ and, with \eqref{eq-claim-MT-9-13}, that 
$$\sum_{j\in{\mathcal D}_i^\star} C_{i,j}>0\hskip.1cm.$$
Then \eqref{eq-claim-MT-9-12} leads to a contradiction. This proves that this first case is absurd. \hfill $\spadesuit$

\medskip {\sc Case 2} -  We assume that $\die\to 0$ and that $\frac{\rie}{\die}\to 0$ as $\eps\to 0$. 

\smallskip\noindent We let ${\ds x_i=\lim_{\eps\to 0} \xie}$. Note that $x_i\in \partial \Omega$. We let $y\in \Omega_0^R$ for some $R>0$ and we set $x_\eps=\xie+\rie y$. Since $\die\to 0$ and $\rie\to 0$, we are in situation b) of Claim \ref{claim-MT-6}. Indeed, as in Case 1, we have that 
$$\frac{\left\vert x_\eps-\xje\right\vert}{\mje}\to +\infty\hbox{ as }\eps\to 0\hbox{ for all }j=1,\dots,N\hskip.1cm.$$
Thus we have that 
\begincal
u_\eps\left(x_\eps\right)&=&\psi_\eps\left(x_\eps\right)+ \sum_{j=1}^N \frac{4\pi+o(1)}{\gje} {\mathcal G}\left(\xje,x_\eps\right)   \\
&&+O\left( \sum_{j\in {\mathcal A}} \left(\gje^{-1}\frac{\mje}{\left\vert \xje-x_\eps\right\vert} +\gje^{-2}\ln\left(\frac{\sje}{\left\vert \xje-x_\eps\right\vert}+2\right)\right)\right)\\
&&+O\left(\sum_{j\in {\mathcal B}} \frac{d_\eps}{d_\eps+d_{j,\eps}} \left(\gje^{-1}\mje + \gje^{-2}\sje\right)\right)
\fincal
where ${\mathcal A}$ is defined as the set of $j\in \left\{1,\dots,N\right\}$ such that $\left\vert \xje-x_\eps\right\vert \le \sje +o\left(d_\eps\right)$ and ${\mathcal B}$ as its complementary. Noting that $\left\vert \xje-x_\eps\right\vert \ge C\rje$ for all $j\in \left\{1,\dots,N\right\}$, we have that for any $j\in {\mathcal A}$,
$$\gje^{-1}\frac{\mje}{\left\vert \xje-x_\eps\right\vert} +\gje^{-2} \ln\left(\frac{\sje}{\left\vert \xje-x_\eps\right\vert}+2\right)=o\left(\gje^{-1}\right)\hskip.1cm.$$
And, for any $j\in {\mathcal B}$, 
$$\frac{d_\eps}{d_\eps+d_{j,\eps}} \left(\gje^{-1}\mje + \gje^{-2}\sje\right) =o\left(\gje^{-1}\right)\hskip.1cm.$$
Note also that, if $j<i$, we have that $j\in {\mathcal B}$ thanks to \eqref{eq-claim-MT-9-1} and that 
$$\frac{d_\eps}{d_\eps+d_{j,\eps}} \left(\gje^{-1}\mje + \gje^{-2}\sje\right) =o\left(d_\eps\gje^{-1}\right)\hskip.1cm.$$
Thus we have that 
$$u_\eps\left(x_\eps\right)=\psi_\eps\left(x_\eps\right)+ \sum_{j=1}^N \frac{4\pi+o(1)}{\gje} {\mathcal G}\left(\xje,x_\eps\right) +o\left(\gie^{-1}\right)+o\left(d_\eps\gamma_{1,\eps}^{-1}\right)\hskip.1cm.$$
We can write thanks to \eqref{eq-MT-20} and since $\psi_\eps=0$ on $\partial \Omega$ that 
$$\psi_\eps\left(x_\eps\right)=O\left(d_\eps \gamma_{1,\eps}^{-3}\right)\hskip.1cm.$$
Then we have that, for any $j<i$,
$${\mathcal G}\left(\xje,x_\eps\right)  = -d_\eps \partial_\nu {\mathcal G}\left(x_j,x_i\right)+o\left(d_\eps\right)\hskip.1cm.$$
And, for $j\ge i$, we have thanks to \eqref{eq-appB-11} that 
$${\mathcal G}\left(\xje,x_\eps\right) = \frac{1}{2\pi}\ln \frac{2\die}{\rie} +\frac{1}{2\pi}\ln\frac{\rie}{\left\vert \xje-x_\eps\right\vert} + O\left(\frac{\rie}{\die}\right) + O\left(\die\right)$$
if $j\in {\mathcal D}_i$ and that 
$${\mathcal G}\left(\xje,x_\eps\right) = {\mathcal G}\left(\xie,\xje\right) + o(1)$$
if $j\not\in {\mathcal D}_i$. We thus arrive to 
\begin{eqnarray}\label{eq-claim-MT-9-14}
u_\eps\left(x_\eps\right)&=&4\pi\frac{ d_{i,\eps}}{\gamma_{1,\eps}}\left(\sum_{j=1}^{i-1} C_{1,j}\left(-\partial_\nu {\mathcal G}\left(x_j,x_i\right)\right) \right)\nonumber\\
&&+2\left(\sum_{j\in {\mathcal D}_1}C_{i,j} +o(1)\right)\gie^{-1}\left(\ln \frac{2\die}{\rie}\right)\\
&&+\sum_{j>i,\, j\not\in{\mathcal D}_i} \left(4\pi+o(1)\right)\gje^{-1} {\mathcal G}\left(\xie,\xje\right) +o\left(\die\gamma_{1,\eps}^{-1}\right)\hskip.1cm.\nonumber
\end{eqnarray}
This gives in particular that $u_\eps\ge C\die\gamma_{1,\eps}^{-1}$ on $\partial{\mathbb D}_{\xie}\left(\rie\right)$ for some $C>0$. Using b) of Claim \ref{claim-MT-0}, we deduce that 
$$C\die\gamma_{1,\eps}^{-1}\le B_{i,\eps}\left(\rie\right)+O\left(\gie^{-1}\right)= -\gie^{-1} \ln \left(\lambda_\eps \gie^2\rie^2\right)+O\left(\gie^{-1}\right)\hskip.1cm.$$
Thanks to Claim \ref{claim-MT-8bis}, this gives that 
$$C\frac{\die\gie}{\gamma_{1,\eps}}\le -\ln \left(\frac{\gie^2\rie^2}{\gamma_{1,\eps}^2}\right)+O(1)\hskip.1cm.$$
Since $\frac{\die}{\rie}\to +\infty$ as $\eps\to 0$ in our case, this implies that 
\begin{equation}\label{eq-claim-MT-9-15}
\gie\rie = o\left(\gamma_{1,\eps}\right)\hskip.1cm.
\end{equation}
We apply now Claim \ref{claim-MT-7} combined with \eqref{eq-MT-20}, \eqref{eq-claim-MT-9-1} and this last estimate to write that 
$$\left\vert \nabla u_\eps\left(x\right)\right\vert \le C \gie^{-1}\left\vert \xie-x\right\vert^{-1} + C\gamma_{1,\eps}^{-1}\le C'  \gie^{-1}\left\vert \xie-x\right\vert^{-1}$$
in ${\mathbb D}_{\xie}\left(\rie\right)$. Thus we can apply Claim \ref{claim-MT-LBA} to $i$~: this gives that, if $\left\vert x\right\vert=\frac{1}{2}$,
$$u_\eps\left(\xie+\rie x\right)= B_{i,\eps}\left(\rie\right)+O\left(\gie^{-1}\right)\hskip.1cm.$$
Combined with \eqref{eq-claim-MT-9-14} and \eqref{eq-claim-MT-9-15}, this leads to 
\begin{eqnarray}\label{eq-claim-MT-9-15bis}
B_{i,\eps}\left(\rie\right)&=& 4\pi\frac{d_{i,\eps}}{\gamma_{1,\eps}}\left(\sum_{j=1}^{i-1} C_{1,j}\left(-\partial_\nu {\mathcal G}\left(x_j,x_i\right)\right)+o(1) \right)\nonumber\\
&&+2\left(\sum_{j\in {\mathcal D}_i}C_{i,j} +o(1)\right)\gie^{-1}\left(\ln \frac{2\die}{\rie}\right)\\
&&+\sum_{j>i,\, j\not\in{\mathcal D}_i} \left(4\pi+o(1)\right)\gje^{-1} {\mathcal G}\left(\xie,\xje\right) \hskip.1cm.\nonumber
\end{eqnarray}
Since 
$$B_{i,\eps}\left(\rie\right)= - \gie^{-1}\ln \left(\lambda_\eps \gie^2 \rie^2\right)+ O\left(\gie^{-1}\right)=-\gie^{-1}\ln \left(\frac{\gie^2 \die^2}{\gamma_{1,\eps}^2}\right)-2\gie^{-1}\ln \left(\frac{\rie}{\die}\right)+ O\left(\gie^{-1}\right)$$
thanks to Claim \ref{claim-MT-8bis}, this leads to 
\begin{eqnarray}\label{eq-claim-MT-9-16}
-\ln \left(\frac{\gie^2 \die^2}{\gamma_{1,\eps}^2}\right)  &=& 4\pi \frac{\die\gie}{\gamma_{1,\eps}}\left(\sum_{j=1}^{i-1} C_{1,j}\left(-\partial_\nu {\mathcal G}\left(x_j,x_i\right)\right)+o(1) \right)\nonumber\\
&&+2\left(\sum_{j\in {\mathcal D}_i^\star}C_{i,j} +o(1)\right)\left(\ln \frac{\die}{\rie}\right) \\
&&+\sum_{j>i,\, j\not\in{\mathcal D}_i} \left(4\pi +o(1)\right) \frac{\gie}{\gje}{\mathcal G}\left(\xie,\xje\right)\hskip.1cm,\nonumber
\end{eqnarray}
from which we can infer that, for $\eps$ small,
\begin{eqnarray}\label{eq-claim-MT-9-17}
2 \frac{\gamma_{1,\eps}}{\die\gie}\ln \left(\frac{\gamma_{1,\eps}}{\die\gie}\right) &\ge & 2\pi\sum_{j=1}^{i-1} C_{1,j}\left(-\partial_\nu {\mathcal G}\left(x_j,x_i\right)\right) \\
&&+2\left(\sum_{j\in {\mathcal D}_i^\star}C_{i,j} +o(1)\right)\frac{\gamma_{1,\eps}}{\die\gie}\left(\ln \frac{\die}{\rie}\right)\hskip.1cm.\nonumber
\end{eqnarray}
Let $j\in {\mathcal D}_i^\star$. Note that, since $\frac{\die}{\rie}\to +\infty$ as $\eps\to 0$, we know that ${\mathcal D}_i^\star\neq\emptyset$. There exists $\delta>0$ such that $\partial {\mathbb D}_{\xje}\left(\delta\rie\right)\subset \Omega_0^R$ for some $R>0$. Thus we can write that 
$$\inf_{\partial {\mathbb D}_{\xje}\left(\delta \rie\right)} u_\eps \ge \bigl(1+o(1)\bigr) B_{i,\eps}\left(\rie\right)$$
thanks to \eqref{eq-claim-MT-9-14} and \eqref{eq-claim-MT-9-15bis}. We can also apply b) of Claim \ref{claim-MT-0} with $r_\eps=\delta \rie$ thanks to Claim \ref{claim-MT-1} and to the fact that ${\ds \frac{\rie}{\rje}\ge 2\left\vert \tilde{x}_j\right\vert^{-1}}+o(1)$. This leads to 
$$B_{j,\eps}\left(\delta \rie\right) \ge \bigl(1+o(1)\bigr) B_{i,\eps}\left(\rie\right)\hskip.1cm.$$
Since 
$$B_{j,\eps}\left(\delta\rie\right)= - \gje^{-1}\ln \left(\lambda_\eps \gje^2 \rie^2\right)+ O\left(\gje^{-1}\right)$$
and 
$$B_{i,\eps}\left(\rie\right)= - \gie^{-1}\ln \left(\lambda_\eps \gie^2 \rie^2\right)+ O\left(\gie^{-1}\right)$$
thanks to Claim \ref{claim-MT-8bis}, we obtain that 
$$- \gje^{-1}\ln \left(\lambda_\eps \gje^2 \rie^2\right)+ O\left(\gje^{-1}\right)\ge - \bigl(1+o(1)\bigr)\gie^{-1}\ln \left(\lambda_\eps \gie^2 \rie^2\right)+O\left(\gie^{-1}\right)\hskip.1cm.$$
This implies since $\gie=O\left(\gje\right)$, see \eqref{eq-MT-14}, that 
$$\ln\left(\lambda_\eps\gie^2 \rie^2\right) \left(1+o(1)-\frac{\gie}{\gje}\right)\ge -C$$
for some $C>0$. Since $\frac{\rie}{\die}\to 0$ as $\eps\to 0$, we get with Claim \ref{claim-MT-2} that $\lambda_\eps\rie^2\gie^2\to 0$ as $\eps\to 0$ and the above implies that $C_{i,j}\ge 1$. Thus we have obtained that 
\begin{equation}\label{eq-claim-MT-9-18}
C_{i,j}\ge 1 \hbox{ for all }j\in {\mathcal D}_i^\star\hskip.1cm.
\end{equation}
Thanks to \eqref{eq-MT-15}, we know that $\rje\ge \rie$ for all $i\in{\mathcal D}_i^\star$. Using Claim \ref{claim-MT-7}, \eqref{eq-MT-20} and \eqref{eq-claim-MT-9-1}, we thus obtain that 
$$\left\vert \nabla u_\eps\right\vert \le C\left(\gamma_{1,\eps}^{-1} + \gie^{-1}\sum_{j\in {\mathcal D}_i} \left\vert \xje-x\right\vert^{-1}\right)$$
in ${\mathbb D}_{\xie}\left(R\rie\right)$ for all $R>0$. Thanks to \eqref{eq-claim-MT-9-15}, this leads to 
$$\left\vert \nabla u_\eps\right\vert \le C \gie^{-1}\sum_{j\in {\mathcal D}_i} \left\vert \xje-x\right\vert^{-1}$$
in ${\mathbb D}_{\xie}\left(R\rie\right)$ for all $R>0$. We are now in position to follow exactly the end of the proof of Case 2 of Claim \ref{claim-MT-8}. We can prove that 
$$\gie\left(u_\eps\left(\xie+\rie x\right)-B_{i,\eps}\left(\rie\right)\right) \to 2\ln\frac{1}{\vert x\vert}+2 \sum_{j\in {\mathcal D}_i^\star} C_{i,j}\ln\frac{1}{\left\vert x-\tilde{x}_j\right\vert}+ A_0$$
in ${\mathcal C}^1_{loc}\left({\mathbb R}^2\setminus {\mathcal S}_i\right)$ as $\eps\to 0$ for some constant $A_0$ and then get a contradiction with Claim \ref{claim-MT-LBA} for $j\in {\mathcal D}_i^\star$ (which is non-empty) such that $\left\vert \tilde{x}_j\right\vert\ge \left\vert \tilde{x}_k\right\vert$ for all $k\in {\mathcal D}_i^\star$. Note here that we assumed that $\rie \to 0$ as $\eps\to 0$, see \eqref{eq-claim-MT-9-2}. This proves that this second case can not happen either. \hfill $\spadesuit$

\medskip {\sc Case 3} - We assume that $\die\to 0$ as $\eps\to 0$ and that ${\ds \frac{\die}{\rie}\to L}$ as $\eps\to 0$ where $L\ge 2$. 

\smallskip\noindent We are thus in the case where, after some harmless rotation,
$$\Omega_0 = {\mathbb R}\times\left(-\infty,L\right)\hskip.1cm.$$
We let $y\in \Omega_0^R$ for some $R>0$ and we set $x_\eps=\xie+\rie x$. Since $\die\to 0$ and $\rie\to 0$, we are in situation b) of Claim \ref{claim-MT-6}. Indeed, as in Case 1, we have that 
$$\frac{\left\vert x_\eps-\xje\right\vert}{\mje}\to +\infty\hbox{ as }\eps\to 0\hbox{ for all }j=1,\dots,N\hskip.1cm.$$
Thus we can write that 
\begincal
u_\eps\left(x_\eps\right)&=&\psi_\eps\left(x_\eps\right)+ \sum_{j=1}^N \frac{4\pi+o(1)}{\gje} {\mathcal G}\left(\xje,x_\eps\right)   \\
&&+O\left( \sum_{j\in {\mathcal A}} \left(\gje^{-1}\frac{\mje}{\left\vert \xje-x_\eps\right\vert} +\gje^{-2} \ln\left(\frac{\sje}{\left\vert \xje-x_\eps\right\vert}+2\right)\right)\right)\\
&&+O\left(\sum_{j\in {\mathcal B}} \frac{d_\eps}{d_\eps+d_{j,\eps}} \left(\gje^{-1}\mje + \gje^{-2}\sje\right)\right)
\fincal
where ${\mathcal A}$ is defined as the set of $j\in \left\{1,\dots,N\right\}$ such that $\left\vert \xje-x_\eps\right\vert \le \sje +o\left(d_\eps\right)$ and ${\mathcal B}$ as its complementary. As in Case 2, we have that 
$$\gje^{-1}\frac{\mje}{\left\vert \xje-x_\eps\right\vert} +\gje^{-2}\ln\left(\frac{\sje}{\left\vert \xje-x_\eps\right\vert}+2\right)=o\left(\gje^{-1}\right)$$
for all $j\in {\mathcal A}$ while 
$$\frac{d_\eps}{d_\eps+d_{j,\eps}} \left(\gje^{-1}\mje + \gje^{-2}\sje\right) = o\left(\gje^{-1}\right)$$
for all $j\in {\mathcal B}$. Note also that, if $j<i$, we have that $j\in {\mathcal B}$ thanks to \eqref{eq-claim-MT-9-1} and that 
$$\frac{d_\eps}{d_\eps+d_{j,\eps}} \left(\gje^{-1}\mje + \gje^{-2}\sje\right) =o\left(d_\eps\gje^{-1}\right)=o\left(\rie \gje^{-1}\right)\hskip.1cm.$$
Thus we have that 
$$u_\eps\left(x_\eps\right)=\psi_\eps\left(x_\eps\right)+ \sum_{i=1}^N \frac{4\pi+o(1)}{\gie} {\mathcal G}\left(\xie,x_\eps\right)+o\left(\gie^{-1}\right)+o\left(\rie\gamma_{1,\eps}^{-1}\right)\hskip.1cm.$$
We can write thanks to \eqref{eq-MT-20} and since $\psi_\eps=0$ on $\partial \Omega$ that 
$$\psi_\eps\left(x_\eps\right)=O\left(\rie \gamma_{1,\eps}^{-3}\right)\hskip.1cm.$$
Then we have that, for any $j<i$,
$${\mathcal G}\left(\xje,x_\eps\right)  = -d_\eps \partial_\nu {\mathcal G}\left(x_j,x_i\right)+o\left(\rie\right)\hskip.1cm.$$
And, for $j\ge i$, we have that 
$${\mathcal G}\left(\xje,x_\eps\right) = \frac{1}{2\pi} \ln \frac{\left\vert \tilde{y}_j-y\right\vert}{\left\vert \tilde{x}_j-y\right\vert} +o(1)$$
if $j\in {\mathcal D}_i$ where 
$$\tilde{y}_j = {\mathcal R}\left(\tilde{x}_j\right)\hskip.1cm,$$
${\mathcal R}$ being the reflection with respect to the straight line ${\mathbb R}\times \left\{L\right\}$. Here we used \eqref{eq-appB-11}. At last, for $j\ge i$ and $j\not\in {\mathcal D}_i$, we have that 
$${\mathcal G}\left(\xje,x_\eps\right) = o(1)$$
thanks to \eqref{eq-appB-11}. This leads to 
\begincal
u_\eps\left(x_\eps\right)&=& 4\pi d_\eps\gamma_{1,\eps}^{-1} \sum_{j=1}^{i-1}\left(-C_{1,j}\partial_\nu {\mathcal G}\left(x_j,x_i\right)\right)\\
&&+\sum_{j\in {\mathcal D}_i} \frac{2+o(1)}{\gje} \ln \frac{\left\vert \tilde{y}_j-y\right\vert}{\left\vert \tilde{x}_j-y\right\vert}+o\left(\gie^{-1}\right)+o\left(\rie\gamma_{1,\eps}^{-1}\right)\hskip.1cm.
\fincal
This gives in particular that $u_\eps\ge C\rie\gamma_{1,\eps}^{-1}$ on $\partial{\mathbb D}_{\xie}\left(\rie\right)$ for some $C>0$. Using b) of Claim \ref{claim-MT-0}, we deduce that 
$$C\rie\gamma_{1,\eps}^{-1}\le B_{i,\eps}\left(\rie\right)+O\left(\gie^{-1}\right)= -\gie^{-1} \ln \left(\lambda_\eps \gie^2\rie^2\right)+O\left(\gie^{-1}\right)\hskip.1cm.$$
Thanks to Claim \ref{claim-MT-8bis}, this gives that 
$$C\frac{\rie\gie}{\gamma_{1,\eps}}\le -\ln \left(\frac{\gie^2\rie^2}{\gamma_{1,\eps}^2}\right)+O(1)\hskip.1cm.$$
This proves that 
$$\rie\gie =O\left(\gamma_{1,\eps}\right)$$
so that, up to a subsequence,
\begin{equation}\label{eq-claim-MT-9-19}
\frac{\rie\gie}{\gamma_{1,\eps}}\to B_0\hbox{ as }\eps\to 0\hskip.1cm.
\end{equation}
Then, by the equation satified by $u_\eps$, it is clear that 
$$v_\eps(x)=\gamma_{i,\eps}u_\eps\left(x_{1,\eps}+\rie x\right)$$ 
has a Laplacian uniformly converging to $0$ in any $\Omega_0^R$. Thus, by standard elliptic theory, we can conclude that 
\begin{equation}\label{eq-claim-MT-9-20}
\gie u_\eps\left(\xie+\rie x\right) \to B_1\left(L-y_2\right) + 2\sum_{j\in {\mathcal D}_i} C_{i,j} \ln \frac{\left\vert \tilde{y}_j-y\right\vert}{\left\vert \tilde{x}_j-y\right\vert}\hbox{ in }{\mathcal C}^1_{loc}\left(\Omega_0\setminus {\mathcal S}_1\right)\hbox{ as }\eps\to 0\hskip.1cm.
\end{equation}
Using \eqref{eq-MT-20}, \eqref{eq-claim-MT-9-19}, \eqref{eq-claim-MT-9-20} and Claim \ref{claim-MT-7}, we have that  
$$\left\vert \nabla u_\eps\right\vert \le C\gamma_{i,\eps}^{-1}\left\vert x_{i,\eps}-x\right\vert^{-1}\hbox{ in }{\mathbb D}_{\xie}\left(\rie\right)\hskip.1cm.$$
We are thus in position to apply the results of Section \ref{section-LBA} to $u_\eps\left(\xie+\cdot\right)$ in the disk ${\mathbb D}_0\left(\rie\right)$. In particular, applying c) of Proposition \ref{claim-LBA-1} and combining it with \eqref{eq-claim-MT-9-20}, we get that 
$$\gie B_{i,\eps}\left(\rie\right)=O\left(1\right)\hskip.1cm.$$
This leads with Claim \ref{claim-appA-4} of Appendix A to 
$$\ln \left(\lambda_\eps \rie^2\gie^2\right) = O\left(1\right)\hskip.1cm.$$
Thanks to Claim \ref{claim-MT-8bis}, we thus have that $B_0>0$ in \eqref{eq-claim-MT-9-19} and 
$B_1>0$ in \eqref{eq-claim-MT-9-20}. We can then proceed exactly as in Case 3 of Claim \ref{claim-MT-8} to get a contradiction in this last case. \hfill $\spadesuit$

\medskip The study of these three cases, all leading to a contradiction, proves that \eqref{eq-claim-MT-9-3} is absurd when we assume \eqref{eq-claim-MT-9-1}. As already said, this permits to prove the claim by induction on $i$. \hfill $\diamondsuit$

\medskip We are now in position to prove Theorem \ref{mainthm}. We know thanks to Claim \ref{claim-MT-9} that 
\begin{equation}\label{eq-MT-21}
\xie\to x_i \hbox{ as }\eps\to 0\hbox{ where }x_i\in \Omega\hskip.1cm.
\end{equation}
Claim \ref{claim-MT-2} then gives that $\lambda_\eps \gie^2=O(1)$ for all $i=1,\dots,N$. Thanks to Claim \ref{claim-MT-8bis} and \eqref{eq-MT-14}, this implies that, up to a subsequence 
\begin{equation}\label{eq-MT-22}
\frac{1}{\sqrt{\lambda_\eps} \gie} \to m_i\hbox{ as }\eps\to 0
\end{equation}
for all $i=1,\dots,N$ with $m_i>0$. Thanks to Claim \ref{claim-MT-6}, to \eqref{eq-MT-20} and to the equation satisfied by $u_\eps$, by standard elliptic theory, we obtain that 
\begin{equation}\label{eq-MT-23}
\frac{u_\eps}{\sqrt{\lambda_\eps}} \to 4\pi \sum_{i=1}^N m_i{\mathcal G}\left(x_i,x\right)\hbox{ in }{\mathcal C}^1_{loc}\left(\Omega\setminus {\mathcal S}\right)
\end{equation}
where ${\mathcal S}=\left\{x_i\right\}_{i=1,\dots,N}$. Moreover, using again \eqref{eq-MT-20} this time together with Claim \ref{claim-MT-7}, we know that 
$$\left\vert \nabla u_\eps\right\vert \le C\sqrt{\lambda_\eps}\sum_{i=1}^N \left\vert \xie-x\right\vert^{-1}$$
in $\Omega$. We are thus in position to apply Claim \ref{claim-MT-LBA} for all $i=1,\dots, N$. This gives that 
\begin{equation}\label{eq-MT-24}
\gie \left(u_\eps\left(\xie+\delta x\right)-B_{i,\eps}\left(\delta\right)\right) \to 2\ln \frac{1}{\left\vert x\right\vert} + {\mathcal H}_i(x)\hbox{ in }{\mathcal C}^1_{loc}\left({\mathbb D}_0(1)\setminus \left\{0\right\}\right)\hbox{ as }\eps\to 0
\end{equation}
where ${\mathcal H}_i\left(0\right)=0$ and $\nabla {\mathcal H}_i\left(0\right)=-\frac{1}{2}\delta \frac{\nabla f_0\left(x_i\right)}{f_0\left(x_i\right)}$. Let us write thanks to Claim \ref{claim-appA-4} that 
\begincal
B_{i,\eps}\left(\delta\right) &=& \gie-\gie^{-1}\left(1+\gie^{-2}\right)\ln \left(1+\frac{\delta^2}{4\mie^2}\right)+O\left(\gie^{-2}\right)\\
&=& \gie -\gie^{-1}\left(1+\gie^{-2}\right)\ln \frac{1}{\mie^2} -\gie^{-1}\ln \frac{\delta^2}{4} +o\left(\gie^{-1}\right)\\
&=& -\gie^{-1} -\gie^{-1} \ln \left(f_0\left(x_i\right)\lambda_\eps \gie^2\right) -\gie^{-1}\ln \frac{\delta^2}{4} +o\left(\gie^{-1}\right)
\fincal
so that, thanks to \eqref{eq-MT-22}, 
$$\gie B_{i,\eps}\left(\delta\right) \to -\ln \frac{ \delta^2 f_0\left(x_i\right)}{4m_i^2}-1 \hskip.1cm.$$
Coming back to \eqref{eq-MT-24} with this, we get that 
\begin{equation}\label{eq-MT-25}
\gie u_\eps\left(x\right) \to 2\ln \frac{1}{\left\vert x-x_i\right\vert} + {\mathcal H}_i\left(\frac{x-x_i}{\delta}\right)-\ln \frac{ f_0\left(x_i\right)}{4m_i^2}-1\hbox{ in }{\mathcal C}^1_{loc}\left({\mathbb D}_{x_i}\left(\delta\right)\setminus \left\{x_i\right\}\right)\hbox{ as }\eps\to 0\hskip.1cm.
\end{equation}
On the other hand, using \eqref{eq-MT-22} and \eqref{eq-MT-23}, we also have that 
\begin{equation}\label{eq-MT-26}
\gie u_\eps\left(x\right) \to \frac{4\pi}{m_i}\sum_{j=1}^N m_j{\mathcal G}\left(x_j,x\right)\hbox{ in }{\mathcal C}^1_{loc}\left({\mathbb D}_{x_i}\left(\delta\right)\setminus \left\{x_i\right\}\right)\hbox{ as }\eps\to 0\hskip.1cm.
\end{equation}
Combining \eqref{eq-MT-25} and \eqref{eq-MT-26}, we get that 
$$m_i{\mathcal H}_i\left(\frac{x-x_i}{\delta}\right) = 4\pi\sum_{j=1}^N m_j{\mathcal G}\left(x_j,x\right)-2m_i\ln \frac{1}{\left\vert x-x_i\right\vert}+m_i\ln \frac{ f_0\left(x_i\right)}{4m_i^2}+m_i\hskip.1cm.$$
Writing 
$${\mathcal G}\left(x,y\right)= \frac{1}{2\pi} \left(\ln \frac{1}{\left\vert x-y\right\vert}+{\mathcal H}\left(x,y\right)\right)\hskip.1cm,$$
this leads to 
$$m_i{\mathcal H}_i\left(\frac{x-x_i}{\delta}\right) = 4\pi\sum_{j\neq i} m_j{\mathcal G}\left(x_j,x\right)+2m_i{\mathcal H}\left(x_i,x\right)+m_i\ln \frac{ f_0\left(x_i\right)}{4m_i^2}+m_i\hskip.1cm.$$
The conditions that ${\mathcal H}_i(0)=0$ and $\nabla{\mathcal H}_i(0)=-\frac{1}{2}\delta \frac{\nabla f_0\left(x_i\right)}{f_0\left(x_i\right)}$ read as 
\begin{equation}\label{eq-MT-27}
 4\pi\sum_{j\neq i} m_j{\mathcal G}\left(x_j,x_i\right)+2m_i{\mathcal H}\left(x_i,x_i\right)+m_i\ln \frac{ f_0\left(x_i\right)}{4m_i^2}+m_i =0 
\end{equation}
and 
\begin{equation}\label{eq-MT-28}
4\pi\sum_{j\neq i} m_j\nabla_y {\mathcal G}\left(x_j,x_i\right)+2m_i\nabla_y {\mathcal H}\left(x_i,x_i\right)=-\frac{1}{2} m_i \frac{\nabla f_0\left(x_i\right)}{f_0\left(x_i\right)}\hskip.1cm.
\end{equation}
This ends the proof of Theorem \ref{mainthm}, up to change the $m_i$'s as in the statement of the theorem. \hfill $\diamondsuit$

\section{Appendix A - The standard bubble}

In this appendix, we develop the exact form of the standard bubble $B_\eps$ which is defined as the radial solution of 
\begin{equation}\label{eq-appA-1}
\Delta B_\eps = \mu_\eps^{-2} \gamma_\eps^{-2} B_\eps e^{B_\eps^2-\gamma_\eps^2} \hbox{ in }\rtwo\hbox{ with } B_\eps(0)=\gamma_\eps
\end{equation}
where $\gamma_\eps\to +\infty$ and $\mu_\eps\to 0$ as $\eps\to 0$. Note that, by standard ordinary differential equations theory, this function is defined on $[0,+\infty)$ and is decreasing. 

We perform the change of variables 
\begin{equation}\label{eq-appA-2}
t=\ln\left(1+\frac{r^2}{4\mu_\eps^2}\right)
\end{equation}
so that we can rewrite equation \eqref{eq-appA-1} as 
\begin{equation}\label{eq-appA-3}
e^{t}\left(\left(1-e^{-t}\right) B_\eps'\right)' = -\frac{B_\eps}{\gamma_\eps^2}e^{2t+B_\eps^2-\gamma_\eps^2}\hskip.1cm.
\end{equation}

We shall need the following lemma which can be proved by direct computations~:

\begin{lemma}\label{lemma-appA}
The solution $\varphi$ of 
$${\mathcal L}\left(\varphi\right)=e^t\left(\left(1-e^{-t}\right) \varphi'\right)' + 2\varphi=F$$
with $\varphi(0)=0$ and $F$ smooth is 
$$\varphi(t)=\int_0^t e^{-s}F(s)\left(\left(1-2e^{-t}\right)\left(1-2e^{-s}\right) \ln \frac{e^t-1}{e^s-1} + 4\left(e^{-s}-e^{-t}\right)\right)\, ds\hskip.1cm.$$
\end{lemma}

\medskip {\it Proof of Lemma \ref{lemma-appA}} - We clearly have that $\varphi(0)=0$ so that we just have to check that $\varphi$ satisfies the given differential equation. Let us differentiate to obtain that 
$$\varphi'(t)= \int_0^t e^{-s}F(s)\left(2e^{-t}\left(1-2e^{-s}\right) \ln \frac{e^t-1}{e^s-1} +\frac{e^t-2}{e^t-1}\left(1-2e^{-s}\right)+ 4e^{-t}\right)\, ds$$
so that 
$$
\left(1-e^{-t}\right)\varphi'(t)=\int_0^t e^{-s}F(s)\left(2\frac{e^t-1}{e^{2t}}\left(1-2e^{-s}\right) \ln \frac{e^t-1}{e^s-1} +\frac{e^t-2}{e^t}\left(1-2e^{-s}\right)+ 4\frac{e^t-1}{e^{2t}}\right)\, ds\hskip.1cm.$$
Differentiating again, we get that 
\begincal
\left(\left(1-e^{-t}\right)\varphi'\right)'(t)&=& e^{-t}F(t)\left(\frac{e^t-2}{e^t}\left(1-2e^{-t}\right)+ 4\frac{e^t-1}{e^{2t}}\right)\\
&&+\int_0^t e^{-s}F(s)\left(-2e^{-t}\left(1-2e^{-t}\right)\left(1-2e^{-s}\right) \ln \frac{e^t-1}{e^s-1}+8e^{-2t}-8e^{-t}e^{-s}\right)\, ds\\
&=& e^{-t}F(t)-2e^{-t}\varphi(t)\hskip.1cm,
\fincal
which proves the lemma. \hfill $\diamondsuit$

\medskip Let us define  
\begin{equation}\label{eq-appA-4}
\varphi_0(t)=\int_0^t e^{-s}\left(s-s^2\right)\left(\left(1-2e^{-t}\right)\left(1-2e^{-s}\right) \ln \frac{e^t-1}{e^s-1} + 4\left(e^{-s}-e^{-t}\right)\right)\, ds
\end{equation}
so that, by lemma \ref{lemma-appA},
\begin{equation}\label{eq-appA-5}
{\mathcal L}\left(\varphi_0\right)(t) = t-t^2\hskip.1cm.
\end{equation}
We claim now that 
\begin{equation}\label{eq-appA-6}
\left\vert \varphi_0(t)+t\right\vert \le C_0\hbox{ and }\varphi_0'(t)\to 1\hbox{ as }t\to +\infty
\end{equation}
for some $C_0>0$. Let us write that 
\begincal
\varphi_0(t)&=&\int_0^t e^{-s}\left(s-s^2\right)\left(2e^{-t}\left(1-2e^{-s}\right) \ln \frac{e^t-1}{e^s-1} +\frac{e^t-2}{e^t-1}\left(1-2e^{-s}\right)+ 4e^{-t}\right)\, ds\\
&=& \left(2 e^{-t}\ln \left(e^t-1\right)+ \frac{e^t-2}{e^t-1} \right)\int_0^t e^{-s}\left(s-s^2\right)\left(1-2e^{-s}\right) \, ds \\
&& + e^{-t} \int_0^t e^{-s}\left(s-s^2\right)\left(2\left(1-2e^{-s}\right) \ln \frac{1}{e^s-1} + 4\right)\, ds\\
&=&\left(2 e^{-t}\ln \left(e^t-1\right)+ \frac{e^t-2}{e^t-1} \right) \left(\left(1+t+t^2\right)e^{-t}-t^2e^{-2t}-1\right)\\
&& + O\left(e^{-t}\right)\\
&=& -1 + O\left(\left(1+t^2\right)e^{-t}\right)\hskip.1cm.
\fincal
This proves the second part of \eqref{eq-appA-6} by passing to the limit $t\to +\infty$ and the first part by integration. 

\medskip We set now 
\begin{equation}\label{eq-appA-7}
B_\eps(t) = \gamma_\eps-\frac{t}{\gamma_\eps} +\gamma_\eps^{-3}\varphi_0 +R_\eps\hskip.1cm.
\end{equation}

\begin{claim}\label{claim-appA-3}
There exists $D_0>0$ such that 
$$\left\vert R_\eps'(t)\right\vert \le D_0 \gamma_\eps^{-5} \hbox{ for all }0\le t\le \gamma_\eps^2-T_\eps$$
where $T_\eps$ is any sequence such that $T_\eps =o\left(\gamma_\eps\right)$ and $\gamma_\eps^ke^{-T_\eps}\to 0$ as $\eps\to 0$ for all $k$.
\end{claim}

\medskip {\it Proof of Claim \ref{claim-appA-3}} - Fix such a sequence $T_\eps$. Let $D_0>0$ that we shall choose later. Since $R_\eps'(0)=0$, there exists $0<t_\eps\le \gamma_\eps^2-T_\eps$ such that 
\begin{equation}\label{eq-claim-appA-3-1}
\left\vert R_\eps'(t)\right\vert \le D_0 \gamma_\eps^{-5} \hbox{ for all }0\le t\le t_\eps\hskip.1cm.
\end{equation}
Note that this implies since $R_\eps(0)=0$ that 
\begin{equation}\label{eq-claim-appA-3-2}
\left\vert R_\eps(t)\right\vert \le D_0 \gamma_{\eps}^{-5} t \hbox{ for all }0\le t\le t_\eps\hskip.1cm.
\end{equation}
We will prove that, for some choice of $D_0$, this $t_\eps$ may be chosen equal to $\gamma_\eps^2-T_\eps$, which will prove the claim. Now, assume this is not the case, then, for the maximal $t_\eps$ such that \eqref{eq-claim-appA-3-1} holds, we have that 
\begin{equation}\label{eq-claim-appA-3-2b}
\left\vert R_\eps'\left(t_\eps\right)\right\vert = D_0 \gamma_\eps^{-5}\hskip.1cm.
\end{equation}
This is the statement we will contradict by an appropriate choice of $D_0$. 
Let us use \eqref{eq-appA-3}, \eqref{eq-appA-5} and \eqref{eq-appA-7} to write that 
$${\mathcal L}\left(R_\eps\right) = F_\eps$$
where 
$$F_\eps = \frac{1}{\gamma_\eps} - \frac{B_\eps}{\gamma_\eps^2}e^{2t+B_\eps^2-\gamma_\eps^2} + 2 R_\eps - \gamma_\eps^{-3}\left(t-t^2-2\varphi_0\right)\hskip.1cm.$$
For $0\le t\le \min\left\{t_\eps,T_\eps\right\}$, we have that 
$$ 2t+B_\eps^2-\gamma_\eps^2 = \frac{t^2}{\gamma_\eps^2} +2\gamma_\eps R_\eps +2\gamma_\eps^{-2}\left(1-\frac{t}{\gamma_\eps^2}\right)\varphi_0 + o\left(\gamma_\eps^{-4}\right)$$
and that 
$$\frac{B_\eps}{\gamma_\eps^2} = \gamma_\eps^{-1} - \gamma_\eps^{-2} t +\gamma_\eps^{-5}\varphi_0 + o\left(\gamma_\eps^{-6}\right)$$
thanks to \eqref{eq-appA-6} and \eqref{eq-claim-appA-3-2}. Thus we have in particular that 
$$\left\vert 2t+B_\eps^2-\gamma_\eps^2 \right\vert \le \frac{2t^2}{\gamma_\eps^2} + 2D_0\gamma_\eps^{-4} t + 2 \gamma_\eps^{-2}\left(C_0+1\right)+o\left(\gamma_\eps^{-4}\right)=o(1)$$
again with \eqref{eq-appA-6} and \eqref{eq-claim-appA-3-2}. We can write that 
\begincal
\left\vert e^{2t+B_\eps^2-\gamma_\eps^2} - 1 - \left(2t+B_\eps^2-\gamma_\eps^2\right)\right\vert
&\le & 2 \left(2t+B_\eps^2-\gamma_\eps^2\right)^2\\
&\le & 20\gamma_\eps^{-4} \left(t^4+\left(C_0+1\right)^2\right)
\fincal
for all $0\le t\le \min\left\{t_\eps,T_\eps\right\}$ for $\eps$ small. Coming back to $F_\eps$, this leads to 
$$\left\vert F_\eps\right\vert \le D_1 \left(1+t^4\right)\gamma_\eps^{-5}$$
for all $0\le t\le \min\left\{t_\eps,T_\eps\right\}$ where $D_1$ depends on $C_0$ but not on $D_0$. We can use the representation formula of Lemma \ref{lemma-appA} to deduce that 
\begincal
\left\vert R_\eps'(t)\right\vert &\le& D_1\gamma_\eps^{-5} \int_0^t e^{-s}\left(1+s^4\right)\left\vert 2e^{-t}\left(1-2e^{-s}\right) \ln \frac{e^t-1}{e^s-1} +\frac{e^t-2}{e^t-1}\left(1-2e^{-s}\right)+ 4e^{-t}\right\vert\, ds \\
&\le & D_2 \gamma_\eps^{-5}
\fincal
for all $0\le t\le \min\left\{T_\eps,t_\eps\right\}$ where $D_2$ depends only on $C_0$, not on $D_0$. Up to choose $D_0>2D_2$, we get that $t_\eps>T_\eps$ thanks to \eqref{eq-claim-appA-3-2b}. Moreover we have that 
\begin{equation}\label{eq-claim-appA-3-3}
\left\vert R_\eps'\left(T_\eps\right)\right\vert \le D_2 \gamma_\eps^{-5}\hskip.1cm.
\end{equation}

\smallskip From now on, we assume that $t_\eps\ge T_\eps$. For all $T_\eps\le t\le \gamma_\eps^2-T_\eps$, we can write that 
$$\left\vert F_\eps(t)\right\vert \le C\gamma_\eps e^{\frac{t^2}{\gamma_\eps^2}} $$
for some $C>0$, depending on $D_0$ and $C_0$. Then we write that 
\begincal 
\left\vert R_\eps'(t)-R_\eps'\left(T_\eps\right)\right\vert &\le& C \gamma_\eps\int_{T_\eps}^t e^{\frac{s^2}{\gamma_\eps^2}-s}\left\vert 2e^{-t}\left(1-2e^{-s}\right) \ln \frac{e^t-1}{e^s-1} +\frac{e^t-2}{e^t-1}\left(1-2e^{-s}\right)+ 4e^{-t}\right\vert\, ds\\
&\le & C\gamma_\eps\int_{T_\eps}^t e^{\frac{s^2}{\gamma_\eps^2}-s}\, ds\\
&= &O\left(\gamma_\eps\int_{T_\eps}^{\gamma_\eps^2-T_\eps} e^{\frac{s^2}{\gamma_\eps^2}-s}\, ds\right)\\
&= & O\left( \gamma_\eps \int_{T_\eps}^{\frac{1}{2}\gamma_\eps^2} e^{\frac{s^2}{\gamma_\eps^2}-s}\, ds\right)\\
&= &O\left(\gamma_\eps \int_{T_\eps}^{\frac{1}{2}\gamma_\eps^2} e^{-\frac{1}{2}s}\, ds\right)\\
&=& O\left(\gamma_\eps e^{-\frac{1}{2}T_\eps}\right)=o\left(\gamma_\eps^{-5}\right)\hskip.1cm.
\fincal
Combined with \eqref{eq-claim-appA-3-3}, this gives that 
$$\left\vert R_\eps'(t)\right\vert \le D_2\gamma_\eps^{-5}+o\left(\gamma_\eps^{-5}\right)\hskip.1cm.$$
This proves that \eqref{eq-claim-appA-3-2b} is impossible, up to choose $D_0\ge 2D_2$. This ends the proof of this claim. \hfill $\diamondsuit$.

\medskip If we want to push a little bit further the estimates, we can get 

\begin{claim}\label{claim-appA-4}
There exists $C_0>0$ such that 
$$\left\vert B_\eps-\gamma_\eps+\frac{t}{\gamma_\eps}+\frac{t}{\gamma_\eps^3}\right\vert \le C_0\gamma_\eps^{-2}$$
for all $0\le t\le \gamma_\eps^2$.
\end{claim}

\medskip {\it Proof of Claim \ref{claim-appA-4}} - It is clear that it holds for any $0\le t\le \gamma_\eps^2-T_\eps$ for $T_\eps$ as in Claim \ref{claim-appA-3}. This is a consequence of Claim \ref{claim-appA-3} and of \eqref{eq-appA-6}. We also know that 
\begin{equation}\label{eq-claim-appA-4-1}
B_\eps\left(\gamma_\eps^2-T_\eps\right)= \frac{T_\eps}{\gamma_\eps} -\frac{1}{\gamma_\eps}+\frac{T_\eps}{\gamma_\eps^3}+O\left(\gamma_\eps^{-3}\right)\hskip.1cm.
\end{equation}
and that 
\begin{equation}\label{eq-claim-appA-4-2}
B_\eps'\left(\gamma_\eps^2-T_\eps\right)= -\frac{1}{\gamma_\eps} -\frac{1}{\gamma_\eps^3}+O\left(\gamma_\eps^{-5}\right)\hskip.1cm.
\end{equation}
Let us integrate twice the equation \eqref{eq-appA-3} between $\gamma_\eps^2-T_\eps$ and $t_\eps=\gamma_\eps^2-\alpha_\eps$ for $0\le \alpha_\eps\le T_\eps$ to write that 
\begin{eqnarray}\label{eq-claim-appA-4-3}
B_\eps\left(t_\eps\right)&=& B_\eps\left(\gamma_\eps^2-T_\eps\right)+B_\eps'\left(\gamma_\eps^2-T_\eps\right)\left(1-e^{T_\eps-\gamma_\eps^2}\right)\ln \left(\frac{e^{\gamma_\eps^2-\alpha_\eps}-1}{e^{\gamma_\eps^2-T_\eps}-1}\right)\\
&&-\frac{1}{\gamma_\eps^2}\int_{\gamma_\eps^2-T_\eps}^{\gamma_\eps^2-\alpha_\eps} 
\ln \left(\frac{e^{t_\eps}-1}{e^t-1}\right)B_\eps(t) e^{t+B_\eps(t)^2-\gamma_\eps^2}\, dt\hskip.1cm.\nonumber
\end{eqnarray}
Using \eqref{eq-claim-appA-4-1} and \eqref{eq-claim-appA-4-2}, and remembering that $\alpha_\eps\le T_\eps=o\left(\gamma_\eps\right)$, we obtain that 
\begin{equation}\label{eq-claim-appA-4-4}
B_\eps\left(t_\eps\right)= \gamma_\eps-\frac{t_\eps}{\gamma_\eps}-\frac{t_\eps}{\gamma_\eps^3} +O\left(\gamma_\eps^{-3}\right) -\frac{1}{\gamma_\eps^2}\int_{\gamma_\eps^2-T_\eps}^{\gamma_\eps^2-\alpha_\eps} 
\ln \left(\frac{e^{t_\eps}-1}{e^t-1}\right)B_\eps(t) e^{t+B_\eps(t)^2-\gamma_\eps^2}\, dt\hskip.1cm.
\end{equation}
Assume that the statement of the Claim holds up to $t_\eps$. If we are able to prove that, under this condition,
\begin{equation}\label{eq-claim-appA-4-5}
\int_{\gamma_\eps^2-T_\eps}^{\gamma_\eps^2-\alpha_\eps} 
\ln \left(\frac{e^{t_\eps}-1}{e^t-1}\right)B_\eps(t) e^{t+B_\eps(t)^2-\gamma_\eps^2}\, dt = o(1)\hskip.1cm,
\end{equation}
then the argument already used in the previous claim will conclude. 

If 
$$\left\vert B_\eps-\gamma_\eps+\frac{t}{\gamma_\eps}+\frac{t}{\gamma_\eps^3}\right\vert \le C_0\gamma_\eps^{-2}$$
for all $0\le t\le t_\eps$, then we can write that 
$$\ln \left(\frac{e^{t_\eps}-1}{e^t-1}\right)\left\vert B_\eps(t)\right\vert e^{t+B_\eps(t)^2-\gamma_\eps^2} = O\left(\gamma_\eps^{-1} \left(1+s^2\right)e^{-s}\right)$$
in $\left[\gamma_\eps^2-T_\eps,t_\eps\right]$ with $t=\gamma_\eps^2-s$ so that it is easily checked that 
$$\int_{\gamma_\eps^2-T_\eps}^{\gamma_\eps^2-\alpha_\eps} 
\ln \left(\frac{e^{t_\eps}-1}{e^t-1}\right)B_\eps(t) e^{t+B_\eps(t)^2-\gamma_\eps^2}\, dt = O\left(\gamma_\eps^{-1}\right)\hskip.1cm,$$
which ends the proof of this claim. \hfill $\diamondsuit$

\begin{claim}\label{claim-appA-5}
There exists $C_1>0$ such that 
$$\left\vert B_\eps'(t) +\gamma_\eps^{-1}\right\vert \le C_1\gamma_\eps^{-2}$$
for all $0\le t\le \gamma_\eps^2$.
\end{claim}

\medskip {\it Proof of Claim \ref{claim-appA-5}} - Let us start from the fact that 
$$ B_\eps'(t) = -\gamma_\eps^{-2} \frac{e^t}{e^t-1}\int_0^t B_\eps(s) e^{s+B_\eps(s)^2-\gamma_\eps^2}\, ds$$
obtained by integrating \eqref{eq-appA-3}. This leads to 
$$\left\vert B_\eps'(t)+\gamma_\eps^{-1}\right\vert \le \gamma_\eps^{-2} \frac{e^t}{e^t-1}\int_0^t \left\vert B_\eps(s) e^{2s+B_\eps(s)^2-\gamma_\eps^2}-\gamma_\eps\right\vert e^{-s}\, ds\hskip.1cm.$$
Let us use Claim \ref{claim-appA-4} to write that 
$$\left\vert B_\eps(s) e^{2s+B_\eps(s)^2-\gamma_\eps^2}-\gamma_\eps \right\vert \le C\gamma_\eps \left(e^{\frac{s^2}{\gamma_\eps^2}}-1\right)e^{-s} +C \frac{s+\gamma_\eps}{\gamma_\eps} e^{\frac{s^2}{\gamma_\eps^2}-s}$$
for some $C>0$ independent of $\eps$ and of $0\le s\le \gamma_\eps^2$. Thus we get that 
$$
\left\vert B_\eps'(t)+\gamma_\eps^{-1}\right\vert \le C\gamma_\eps^{-1} \frac{e^t}{e^t-1} \int_0^t\left(e^{\frac{s^2}{\gamma_\eps^2}}-1\right)e^{-s}\, ds +C\gamma_\eps^{-3}\frac{e^t}{e^t-1} \int_0^t \left(s+\gamma_\eps\right)e^{\frac{s^2}{\gamma_\eps^2}-s}\, ds\hskip.1cm.
$$
Arguing as above, one gets that 
$$\frac{e^t}{e^t-1} \int_0^t\left(e^{\frac{s^2}{\gamma_\eps^2}}-1\right)e^{-s}\, ds \le C\gamma_\eps^{-2}$$
and that 
$$\frac{e^t}{e^t-1} \int_0^t \left(s+\gamma_\eps\right)e^{\frac{s^2}{\gamma_\eps^2}-s}\, ds \le C\gamma_\eps$$
for all $0\le t\le \gamma_\eps^2$. This permits to end the proof of the claim. \hfill $\diamondsuit$

\section{Appendix B - Estimates on the Green function}\label{section-Green}

We list and prove some useful estimates on the Green function of the Laplacian with Dirichlet boundary condition in some smooth domain $\Omega$. We fix such a two-dimensional domain and we let ${\mathcal G}\left(x,y\right)$ be such that 
$$\Delta_x {\mathcal G}(x,y)=\delta_y\hbox{ with }{\mathcal G}\left(x,y\right)=0\hbox{ if }x\in \partial \Omega\hskip.1cm.$$
It is well known that ${\mathcal G}$ is symmetric and smooth outside of the diagonal. Except on the disk of radius $R$ where ${\mathcal G}$ is explicitly given by 
$${\mathcal G}\left(x,y\right)= \frac{1}{4\pi} \ln \frac{\left\vert \frac{\vert y\vert}{R}x-\frac{R y}{\vert y\vert}\right\vert^2}{\left\vert x-y\right\vert^2}$$
and so where all the estimates below follow from explicit computations, we need to be a little bit careful to estimate the Green function for various $x$ and $y$. 

We know that 
\begin{equation}\label{eq-appB-1}
{\mathcal G}\left(x,y\right)= \frac{1}{2\pi} \ln \frac{1}{\left\vert x-y\right\vert} + {\mathcal H}_y(x)
\end{equation}
where 
$$\Delta_x {\mathcal H}_y(x)=0\hbox{ in }\Omega \hbox{ and } {\mathcal H}_y(x)=-\frac{1}{2\pi} \ln \frac{1}{\left\vert x-y\right\vert }\hbox{ on }\partial \Omega\hskip.1cm.$$
First, if $y\in K$ for some compact subset $K$ of $\Omega$, we clearly have that 
\begin{equation}\label{eq-appB-1bis} 
\left\vert {\mathcal H}_y\left(x\right)\right\vert \le C_K \hbox{ and }\left\vert \nabla {\mathcal H}_y(x)\right\vert \le C_K
\end{equation}
for some $C_K>0$ for all $x\in \Omega$ so that 
\begin{eqnarray}\label{eq-appB-2} 
&&\left\vert {\mathcal G}\left(x,y\right)-\frac{1}{2\pi}\ln \frac{1}{\left\vert x-y\right\vert} \right\vert \le C_K\hskip.1cm,\nonumber\\
&&\left\vert \nabla_x {\mathcal G}\left(x,y\right)+\frac{1}{2\pi}\frac{x-y}{\left\vert x-y\right\vert^2}\right\vert \le C_K\hskip.1cm,\\
&& \left\vert \nabla_x {\mathcal G}\left(x,y\right)\right\vert \le C_K \left\vert x-y\right\vert^{-1}\hskip.1cm,\nonumber\\
&&\left\vert {\mathcal G}\left(x,y\right)-{\mathcal G}\left(z,y\right)-\frac{1}{2\pi}\ln \frac{\left\vert z-y\right\vert}{\left\vert x-y\right\vert}\right\vert \le C_K\left\vert x-z\right\vert \nonumber
\end{eqnarray}
for all $x,y,z\in K\subset\subset \Omega$ (distinct points). 

We let now $\left(y_\eps\right)$ be a sequence of points in $\Omega$ such that 
\begin{equation}\label{eq-appB-3}
d_\eps=d\left(y_\eps,\partial \Omega\right) \to 0\hbox{ as }\eps\to 0\hskip.1cm.
\end{equation}
We let now $\tilde{y}_\eps\in {\mathbb R}^2$ be such that 
\begin{equation}\label{eq-appB-4}
\tilde{y}_\eps = 2\pi\left(y_\eps\right) -y_\eps 
\end{equation}
where $\pi$ is the projection on the boundary of $\Omega$. Note that $\pi\left(y_\eps\right)$ is unique thanks to \eqref{eq-appB-3} and to the fact that $\Omega$ is smooth. Moreover, we have that 
\begin{equation}\label{eq-appB-5}
\tilde{y}_\eps = y_\eps + 2 d_\eps \nu_\eps
\end{equation}
where $\nu_\eps$ is the unit outer normal  of $\partial\Omega$ at $\pi\left(y_\eps\right)$. We let now 
\begin{equation}\label{eq-appB-6}
{\mathcal G}\left(x,y_\eps\right)=\frac{1}{2\pi}\ln\frac{\left\vert x-\tilde{y}_\eps\right\vert}{\left\vert x-y_\eps\right\vert} + \tilde{\mathcal H}_\eps\left(x\right)
\end{equation}
where $\tilde{\mathcal H}_\eps$ is harmonic in $\Omega$ and satisfies 
\begin{equation}\label{eq-appB-7}
\tilde{\mathcal H}_\eps\left(x\right)=-\frac{1}{2\pi}\ln\frac{\left\vert x-\tilde{y}_\eps\right\vert}{\left\vert x-y_\eps\right\vert}\hbox{ on }\partial \Omega\hskip.1cm.
\end{equation}
It is easily checked since $\Omega\in {\mathcal C}^2$ that 
$$\left\vert \tilde{\mathcal H}_\eps\left(x\right)\right\vert \le C_{\Omega} d_\eps$$
for some $C_\Omega>0$ independent of $\eps$ and for all $x\in \partial \Omega$. Thus we have that 
\begin{equation}\label{eq-appB-8}
\left\vert \tilde{\mathcal H}_\eps\left(x\right)\right\vert \le C_\Omega d_\eps \hbox{ in }\Omega\hskip.1cm.
\end{equation}
It is also easily checked that 
\begin{equation}\label{eq-appB-9}
\left\vert \nabla^T \tilde{\mathcal H}_\eps\left(x\right)\right\vert \le C_\Omega
\end{equation}
for all $x\in \partial \Omega$ where $\nabla^T$ denotes the tangential derivative. Thus we have that 
\begin{equation}\label{eq-appB-10}
\left\vert \nabla  \tilde{\mathcal H}_\eps\left(x\right)\right\vert \le C_\Omega\frac{d_\eps}{d_\eps+d\left(x,\partial \Omega\right)}\hbox{ in }\Omega\hskip.1cm.
\end{equation}
Let us give some useful consequences of \eqref{eq-appB-8} and \eqref{eq-appB-10}. Let $y_\eps$ be such that $d_\eps=d\left(y_\eps,\partial \Omega\right)\to 0$ as $\eps\to 0$, then we have that for any sequence $\left(x_\eps\right)$ in $\Omega$
\begin{equation}\label{eq-appB-11}
\begin{array}{ll}
{\ds {\mathcal G}\left(x_\eps,y_\eps\right) = O\left(\frac{d_\eps}{\left\vert x_\eps-y_\eps\right\vert}\right)} &{\ds \hbox{ if }d_\eps = O\left(\left\vert x_\eps-y_\eps\right\vert\right) }\\
\,&\,\\
{\ds {\mathcal G}\left(x_\eps,y_\eps\right) = \frac{1}{2\pi}\ln \frac{2d_\eps}{\left\vert x_\eps-y_\eps\right\vert}+ O\left(\frac{\left\vert x_\eps-y_\eps\right\vert}{d_\eps}\right)+O\left(d_\eps\right) }&{\ds \hbox{ if }\left\vert x_\eps-y_\eps\right\vert=o\left(d_\eps\right) }\\
\,&\,\\
{\ds \left\vert \nabla_x{\mathcal G}\left(x_\eps,y_\eps\right)\right\vert = O\left(\frac{d_\eps}{d_\eps+d\left(x_\eps,\partial \Omega\right)}\right)}&{\ds \hbox{ if }d_\eps = O\left(\left\vert x_\eps-y_\eps\right\vert\right) }\\
\,&\,\\
{\ds \left\vert \nabla_x{\mathcal G}\left(x_\eps,y_\eps\right)\right\vert = \frac{1}{2\pi \left\vert x_\eps-y_\eps\right\vert} + O\left(\frac{1}{d_\eps}\right)}&{\ds \hbox{ if }\left\vert x_\eps-y_\eps\right\vert=o\left(d_\eps\right) }
\end{array}
\end{equation}
These are the only estimates which were used in this paper. 

\section*{Acknowledgments} The authors wish to thank Luca Martinazzi to have drawn their attention again on this question and for the fruitful discussions they had during the preparation of this work.


\begin{thebibliography}{10}

\bibitem{AdiDruet}
Adimurthi and O.~Druet.
\newblock Blow up analysis in dimension 2 and a sharp form of
  {T}rudinger-{M}oser inequality.
\newblock {\em Comm. PDE's}, 29, 1-2:295--322, 2004.

\bibitem{AdiPrashanth1}
Adimurthi and S.~Prashanth.
\newblock Failure of {P}alais-{S}male condition and blow-up analysis for the
  critical exponent problem in {${\mathbb R}^2$}.
\newblock {\em Proc. Indian. Acad. Sci. (Math. Sci.)}, 107, 3:283--317, 1997.

\bibitem{AdiStruwe}
Adimurthi and M.~Struwe.
\newblock Global compactness properties of semilinear elliptic equations with
  critical exponential growth.
\newblock {\em J.F.A.}, 175:125--167, 2000.

\bibitem{Brendle}
S.~Brendle.
\newblock Blow-up phenomena for the {Y}amabe equation.
\newblock {\em J. Amer. Math. Soc.}, 21(4):951--979, 2008.

\bibitem{BrendleMarques}
S.~Brendle and F.C. Marques.
\newblock Blow-up phenomena for the {Y}amabe equation. {II}.
\newblock {\em J. Differential Geom.}, 81(2):225--250, 2009.

\bibitem{Carleson-Chang}
L.~Carleson and S.Y.A. Chang.
\newblock On the existence of an extremal function for an inequality of {J}.
  {M}oser.
\newblock {\em Bull. Sc. Math.}, 110:113--127, 1986.

\bibitem{Chen-Lin}
C.C. Chen and C.S. Lin.
\newblock Sharp estimates for solutions of multi-bubbles in compact {R}iemann
  surfaces.
\newblock {\em Comm. Pure Appl. Math.}, 55(6):728--771, 2002.

\bibitem{ChenLin-Liouville}
C.C. Chen and C.S. Lin.
\newblock Topological degree for a mean field equation on {R}iemann surfaces.
\newblock {\em Comm. Pure Appl. Math.}, 56(12):1667--1727, 2003.

\bibitem{ChenLi}
W.~Chen and C.~Li.
\newblock Classification of solutions of some nonlinear elliptic equations.
\newblock {\em Duke Math. J.}, 63:615--623, 1991.

\bibitem{Coron}
J.-M. Coron.
\newblock Topologie et cas limite des injections de {S}obolev.
\newblock {\em C. R. Acad. Sci. Paris S\'er. I Math.}, 299(7):209--212, 1984.

\bibitem{CostaTintarev}
D.G. Costa and C.~Tintarev.
\newblock Concentration profiles for the {T}rudinger-{M}oser functional are
  shaped like toy pyramids.
\newblock {\em J. Funct. Anal.}, 266(2):676--692, 2014.

\bibitem{DPMR1}
M.~del Pino, M.~Musso, and B.~Ruf.
\newblock New solutions for {T}rudinger-{M}oser critical equations in {$\Bbb
  R^2$}.
\newblock {\em J. Funct. Anal.}, 258(2):421--457, 2010.

\bibitem{DPMR2}
M.~del Pino, M.~Musso, and B.~Ruf.
\newblock Beyond the {T}rudinger-{M}oser supremum.
\newblock {\em Calc. Var. Partial Differential Equations}, 44(3-4):543--576,
  2012.

\bibitem{DruetJDG}
O.~Druet.
\newblock From one bubble to several bubbles~: the low-dimensional case.
\newblock {\em J. Diff. Geom.}, 63:399--473, 2003.

\bibitem{DruetDMJ}
O.~Druet.
\newblock Multibumps analysis in dimension 2 - {Q}uantification of blow up
  levels.
\newblock {\em Duke Math. J.}, 132 (2):217--269, 2006.

\bibitem{DruetENSAIOS}
O.~Druet.
\newblock La notion de stabilit{\'e} pour des {\'e}quations aux
  d{\'e}riv{\'e}es partielles elliptiques.
\newblock {\em ENSAIOS Matematicos}, 19:1--100, 2010.

\bibitem{DHRbook}
O.~Druet, E.~Hebey, and F.~Robert.
\newblock {\em Blow-up theory for elliptic {P}{D}{E}s in {R}iemannian
  geometry}, volume~45 of {\em Mathematical Notes}.
\newblock Princeton University Press, 2004.

\bibitem{Flucher}
M.~Flucher.
\newblock Extremal functions for the {T}rudinger-{M}oser inequality in $2$
  dimensions.
\newblock {\em Comm. Math. Helv.}, 67:471--497, 1992.

\bibitem{KhuriMarquesSchoen}
M.~A. Khuri, F.~C. Marques, and R.~M. Schoen.
\newblock A compactness theorem for the {Y}amabe problem.
\newblock {\em J. Differential Geom.}, 81(1):143--196, 2009.

\bibitem{LRS}
T.~Lamm, F.~Robert, and M.~Struwe.
\newblock The heat flow with a critical exponential nonlinearity.
\newblock {\em J. Funct. Anal.}, 257(9):2951--2998, 2009.

\bibitem{Laurain}
P.~Laurain.
\newblock Concentration of {$CMC$} surfaces in a 3-manifold.
\newblock {\em Int. Math. Res. Not.}, 24:5585--5649, 2012.

\bibitem{LiZhang}
Yan~Yan Li and Lei Zhang.
\newblock A {H}arnack type inequality for the {Y}amabe equation in low
  dimensions.
\newblock {\em Calc. Var. Partial Differential Equations}, 20(2):133--151,
  2004.

\bibitem{LiShafrir}
Y.Y. Li and I.~Shafrir.
\newblock Blow-up analysis for solutions of {$-\Delta u=Ve^u$} in dimension
  two.
\newblock {\em Indiana Univ. Math. J.}, 43(4):1255--1270, 1994.

\bibitem{LiZhang2}
Y.Y. Li and L.~Zhang.
\newblock Compactness of solutions to the {Y}amabe problem. {II}.
\newblock {\em Calc. Var. Partial Differential Equations}, 24(2):185--237,
  2005.

\bibitem{LiZhang3}
Y.Y. Li and L.~Zhang.
\newblock Compactness of solutions to the {Y}amabe problem. {III}.
\newblock {\em J. Funct. Anal.}, 245(2):438--474, 2007.

\bibitem{LiZhu}
Y.Y. Li and M.~Zhu.
\newblock Yamabe type equations on three dimensional {R}iemannian manifolds.
\newblock {\em Communications in Contemporary Mathematics}, 1:1--50, 1999.

\bibitem{Linsurvey}
C.-S. Lin.
\newblock An expository survey on the recent development of mean field
  equations.
\newblock {\em Discrete Contin. Dyn. Syst.}, 19(2):387--410, 2007.

\bibitem{Lions2}
P.L. Lions.
\newblock The concentration-compactness principle in the calculus of
  variations. {T}he limit case. {P}art {I}.
\newblock {\em Rev. Mat. Iberoamericano}, 1.1:145--201, 1985.

\bibitem{Malchiodi}
A.~Malchiodi.
\newblock Morse theory and a scalar field equation on compact surfaces.
\newblock {\em Adv. Differential Equations}, 13(11-12):1109--1129, 2008.

\bibitem{MalchiodiMartinazzi}
A.~Malchiodi and L.~Martinazzi.
\newblock Critical points of the {M}oser-{T}rudinger functional on a disk.
\newblock {\em J. Eur. Math. Soc. (JEMS)}, 16(5):893--908, 2014.

\bibitem{ManciniMartinazzi}
Gabriele Mancini and Luca Martinazzi.
\newblock The {M}oser-{T}rudinger inequality and its extremals on a disk via
  energy estimates.
\newblock {\em Calc. Var. Partial Differential Equations}, 56(4):Art. 94, 26,
  2017.

\bibitem{Marques}
F.C. Marques.
\newblock A priori estimates for the {Y}amabe problem in the non-locally
  conformally flat case.
\newblock {\em J. Differential Geom.}, 71(2):315--346, 2005.

\bibitem{Marques2}
F.C. Marques.
\newblock Blow-up examples for the {Y}amabe problem.
\newblock {\em Calc. Var. Partial Differential Equations}, 36(3):377--397,
  2009.

\bibitem{Martinazzi}
Luca Martinazzi.
\newblock A threshold phenomenon for embeddings of {$H^m_0$} into {O}rlicz
  spaces.
\newblock {\em Calc. Var. Partial Differential Equations}, 36(4):493--506,
  2009.

\bibitem{Moser}
J.~Moser.
\newblock A sharp form of an inequality by {N}. {T}rudinger.
\newblock {\em Indiana U. Math. J.}, 20, 11:1077--1092, 1971.

\bibitem{SchoenZhang}
R.~Schoen and D.~Zhang.
\newblock Prescribed scalar curvature on the $n$-sphere.
\newblock {\em Calculus of Variations and PDE's}, 4:1--25, 1996.

\bibitem{Struwe-MT}
M.~Struwe.
\newblock Critical points of embeddings of {$H^{1,n}_0$} into {O}rlicz spaces.
\newblock {\em Ann. Inst. H. Poincar\'e Anal. Non Lin\'eaire}, 5(5):425--464,
  1988.

\bibitem{Struwe-MT2}
M.~Struwe.
\newblock Positive solutions of critical semilinear elliptic equations on
  non-contractible planar domains.
\newblock {\em J. Eur. Math. Soc. (JEMS)}, 2(4):329--388, 2000.

\bibitem{Takahashi}
F.~Takahashi.
\newblock Blow up points and the {M}orse indices of solutions to the
  {L}iouville equation in two-dimension.
\newblock {\em Adv. Nonlinear Stud.}, 12(1):115--122, 2012.

\bibitem{Trudinger}
N.S. Trudinger.
\newblock On embedding into {O}rlicz spaces and some applications.
\newblock {\em J. Math. Mech.}, 17:473--483, 1967.

\end{thebibliography}
\end{document}